\documentclass[12pt]{report}

\input xy

\xyoption{all}

\usepackage{amssymb,amsbsy,amsthm,amsmath,graphicx,epsfig,times}

\newtheorem{thm}{Theorem}[section]
\newtheorem{lem}[thm]{Lemma}
\newtheorem{prop}[thm]{Proposition}

\newtheorem{conj}[thm]{Conjecture}
\newtheorem{dfn}[thm]{Definition}

\theoremstyle{remark}

\newcommand{\bs}[1]{\boldsymbol{#1}}
\renewcommand{\bf}[1]{\mathbf{#1}}
\renewcommand{\rm}[1]{\mathrm{#1}}
\renewcommand{\cal}[1]{\mathcal{#1}}

\newcommand{\bbN}{\mathbb{N}}

\newcommand{\bbT}{\mathbb{T}}
\newcommand{\bbZ}{\mathbb{Z}}

\newcommand{\bfW}{\mathbf{W}}
\newcommand{\bfX}{\mathbf{X}}
\newcommand{\bfY}{\mathbf{Y}}
\newcommand{\bfZ}{\mathbf{Z}}

\newcommand{\sfA}{\mathsf{A}}
\newcommand{\sfC}{\mathsf{C}}
\newcommand{\sfD}{\mathsf{D}}
\newcommand{\sfE}{\mathsf{E}}

\newcommand{\sfZ}{\mathsf{Z}}

\renewcommand{\d}{\mathrm{d}}

\newcommand{\B}{\mathcal{B}}

\newcommand{\I}{\mathcal{I}}

\newcommand{\frH}{\mathfrak{H}}

\newcommand{\G}{\Gamma}
\renewcommand{\L}{\Lambda}

\renewcommand{\S}{\Sigma}

\renewcommand{\a}{\alpha}
\renewcommand{\b}{\beta}

\newcommand{\g}{\gamma}
\renewcommand{\l}{\lambda}
\newcommand{\s}{\sigma}

\newcommand{\id}{\mathrm{id}}
\newcommand{\nil}{\mathrm{nil}}

\newcommand{\ol}[1]{\overline{#1}}

\renewcommand{\t}[1]{\tilde{#1}}
\newcommand{\into}{\hookrightarrow}
\newcommand{\fin}{\nolinebreak\hspace{\stretch{1}}$\lhd$\newline}
\renewcommand{\qed}{\nolinebreak\hspace{\stretch{1}}$\Box$\newline}

\newcommand{\actson}{\curvearrowright}
\renewcommand{\to}{\longrightarrow}
\newcommand{\uhr}{\upharpoonright}

\begin{document}

\title          {\textbf{Multiple recurrence and the structure of probability-preserving systems}}
\author         {Tim Austin\\ \\ \small{\emph{Department of Mathematics}}\\ \small{\emph{University of California, Los Angeles}}}
\date{}

\maketitle

\tableofcontents

\chapter*{Preface}

In 1975 Szemer\'edi proved the long-standing conjecture of Erd\H{o}s and Tur\'an that any subset of $\bbZ$ having positive upper Banach density contains arbitrarily long arithmetic progressions.  Szemer\'edi's proof was entirely combinatorial, but two years later Furstenberg gave a quite different proof of Szemer\'edi's Theorem by first showing its equivalence to an ergodic-theoretic assertion of multiple recurrence, and then bringing new machinery in ergodic theory to bear on proving that.  His ergodic-theoretic approach subsequently yielded several other results in extremal combinatorics, as well as revealing a range of new phenomena according to which the structures of probability-preserving systems can be described and classified.

In this work I survey some recent advances in understanding these ergodic-theoretic structures.  It contains proofs of the norm convergence of the `nonconventional' ergodic averages that underly Furstenberg's approach to variants of Szemer\'edi's Theorem, and of two of the recurrence theorems of Furstenberg and Katznelson: the Multidimensional Multiple Recurrence Theorem, which implies a multidimensional generalization of Szemer\'edi's Theorem; and a density version of the Hales-Jewett Theorem of Ramsey Theory.

\[\ast\quad\ast\quad\ast\]

The text below was originally submitted as my Ph.D. dissertation at UCLA, after being assembled from a number of earlier papers.  It seems worth repeating the acknowledgements from that dissertation as well.

Many people deserve my thanks for their part in my mathematical education.  Listing them in roughly the order we met, I must at least mention David Fremlin, Imre Leader, Tim Gowers, B\'ela Bollob\'as, Ben Garling, James Norris, Assaf Naor, Yuval Peres, Vitaly Bergelson, Christoph Thiele, Sorin Popa, David Aldous, Tamar Ziegler, Bryna Kra, Bernard Host, Mariusz Lema\'nczyk and Dan Rudolph.  I could have written a much longer list, but still the selection would have been slightly arbitrary: to make it complete would have required far more space than I have available.

During the same period, I have benefited from the financial support of Trinity College, Cambridge, the Shapiro and Huang Families through their UCLA graduate student fellowships, and Microsoft Corporation.  No less significant, I have been able to rely unquestioningly on the support of family and friends, for whom I can only hope to be so generous in turn should the need arise.

Terence Tao, who advised this dissertation, has certainly taught me more during the last four years than either of us fully appreciates, and his energy and enthusiasm for mathematics are a constant motivation for those around him.

\begin{flushright}
Venice Beach, California

May 2010
\end{flushright}

\chapter{Introduction}

The concerns of this work stem from the following remarkable result of
Szemer\'edi~(\cite{Sze75}), which confirmed an old conjecture of
Erd\H{o}s and Tur\'an~(\cite{ErdTur36}).

\vspace{7pt}

\noindent\textbf{Szemer\'edi's Theorem.}\quad
\emph{For any $\delta > 0$ and $k\geq 1$ there is some $N_0 \geq 1$ such
that if $N \geq N_0$ then any $A \subseteq \{1,2,3,\ldots,N\}$ with
$|A| \geq \delta N$ includes a nontrivial $k$-term arithmetic
progression: $A \supseteq \{a,a+n,\ldots,a+(k-1)n\}$ for some $a \in
\{1,2,\ldots,N\}$ and $n \geq 1$.}

\vspace{7pt}

This provides a considerable strengthening of a much older result of van der Waerden~\cite{vdW27}, according to which any colouring of $\bbN$ using a bounded number of colours witnesses arbitrarily long finite arithmetic progressions that are monochromatic.  Since any colouring with at most $c$ colours must have at least one colour class of upper Banach density at least $1/c$, van der Waerden's Theorem can be deduced by applying Szemer\'edi's Theorem to the intersection of that class with sufficiently long discrete intervals in $\bbN$.

Shortly after the appearance of Szemer\'edi's ingenious
combinatorial proof, Furstenberg gave a new proof of the above
theorem in~\cite{Fur77} using a superficially quite different
approach, relying on a conversion to a problem about
probability-preserving dynamical systems.

Such a system consists of a probability space $(X,\S,\mu)$ together
with an invertible, measurable, $\mu$-preserving transformation
$T:X\to X$.  Furstenberg proved that all such systems enjoy a
property of `multiple recurrence':

\vspace{7pt}

\noindent\textbf{Multiple Recurrence Theorem.}\quad
\emph{Whenever $(X,\S,\mu)$ and $T$ are as above, if $k\geq 1$ and $A \in \S$ has $\mu(A)
> 0$ then
\[\liminf_{N\to\infty}\frac{1}{N}\sum_{n=1}^N\mu(A\cap T^{-n}(A)\cap\cdots\cap T^{-(k-1)n}(A)) > 0.\]
In particular, there is some $n \geq 1$ such that
\[\mu(A\cap T^{-n}(A)\cap\cdots\cap T^{-(k-1)n}(A)) > 0.\]}

\vspace{7pt}

It is worth noting that analogously to this ergodic-theoretic proof of Szemer\'edi's Theorem, it is possible to deduce the colouring theorem of van der Waerden from a multiple recurrence result in topological dynamics.  We will not be concerned with this story here, but it is reported in detail in Furstenberg's book~\cite{Fur81}.

Shortly after the above result appeared, Furstenberg and Katznelson realized that the
same basic method could be modified to apply to collections of
commuting measure-preserving transformations, and proved the following in~\cite{FurKat78}.

\vspace{7pt}

\noindent\textbf{Theorem A} (Multidimensional Multiple Recurrence Theorem).\quad \emph{If $(X,\S,\mu)$ is a probability space, $T_1$, $T_2$, \ldots, $T_d$
are commuting measurable invertible $\mu$-preserving self-maps of
$X$ and $A \in \S$ has $\mu(A) > 0$, then
\[\liminf_{N\to\infty}\frac{1}{N}\sum_{n=1}^N\mu(T_1^{-n}(A)\cap\cdots\cap T_d^{-n}(A)) > 0.\]}

\vspace{7pt}

Of course this result implies one-dimensional multiple recurrence by setting $d := k$ and $T_i := T^i$ for $i=0,1,\ldots,k-1$.  In addition, Furstenberg and Katznelson were able to convert Theorem A back into a multidimensional
combinatorial result generalizing Szemer\'edi's Theorem.

\vspace{7pt}

\noindent\textbf{Multidimensional Szemer\'edi Theorem.}\quad \emph{For any $\delta > 0$ and $d\geq 1$ there is some $N_0 \geq 1$ such
that if $N \geq N_0$ then any $A \subseteq \{1,2,\ldots,N\}^d$ with
$|A| \geq \delta N^d$ includes the vertex set of the outer face of a nontrivial upright simplex:
\[A \supseteq \{\bf{a} + n\bf{e}_1,\bf{a}+n\bf{e}_2,\ldots,\bf{a}+n\bf{e}_d\}\]
for some $\bf{a} \in \{1,2,\ldots,N\}^d$ and $n \geq 1$, where
$\bf{e}_1$, $\bf{e}_2$, \ldots, $\bf{e}_d$ are the usual basis
vectors of $\bbZ^d$.}

\vspace{7pt}

This ergodic-theoretic approach to results in additive combinatorics
has since developed into a whole subdiscipline, sometimes termed
`Ergodic Ramsey Theory'; see, for instance, Bergelson's
survey~\cite{Ber96}. In particular, Furstenberg and Katznelson used
this approach to prove a number of further results concerning some
form of `recurrence', culminating in the following density version
of the classical Hales-Jewett Theorem~\cite{HalJew63} proved
in~\cite{FurKat91}:

\vspace{7pt}

\noindent\textbf{Theorem B} (Density version of the Hales-Jewett Theorem).\quad
\emph{For any $\delta > 0$ and $k\geq 1$ there is some $N_0 \geq 1$ such
that if $N \geq N_0$ then any $A \subseteq [k]^N$ with $|A| \geq
\delta k^N$ includes a \textbf{combinatorial line}: a subset
$L\subseteq [k]^N$ of the form
\[L = \{w\in [k]^N:\ w|_{[N]\setminus J} = w_0,\,\&\,w_j\ \hbox{is the same element of $[k]$ for all $j \in J$}\},\]
for some fixed nonempty $J \subseteq [N]$ and $w_0 \in
[k]^{[N]\setminus J}$.}

\vspace{7pt}

In fact, this result implies most of the other main results in
density Ramsey Theory, including Szemer\'edi's Theorem and its
multidimensional generalization.  This implication holds exactly as in the older setting of colouring Ramsey Theorems, which is well-treated in the book~\cite{GraRotSpe90} of Graham,
Rothschild and Spencer.

In addition to achieving some striking new combinatorial results,
Ergodic Ramsey Theory has also motivated new ergodic-theoretic
questions, and has witnessed an ongoing interplay between insights
into these two aspects of the subject.

One basic question that was resolved only recently is whether the
`multiple ergodic averages' studied in Theorems A and B above actually converge (that is, whether
`$\liminf$' can be replaced with `$\lim$').  In the case of
the original Multiple Recurrence Theorem, this was finally shown to be so by Host and
Kra in~\cite{HosKra05}, following the establishment of several
special cases and related results over two decades
in~\cite{ConLes84,ConLes88.1,ConLes88.2,Zha96,FurWei96,HosKra01}
(see also Ziegler's paper~\cite{Zie07} for another proof of the Host-Kra result). The more
general setting of Theorem A was then settled by
Tao in~\cite{Tao08(nonconv)}.

\vspace{7pt}

\noindent\textbf{Theorem C} (Norm convergence of nonconventional
averages).\quad \emph{For any commuting tuple of invertible
measurable $\mu$-preserving transformations $T_1$, $T_2$, \ldots,
$T_d\curvearrowright (X,\S,\mu)$ and any functions
$f_1,f_2,\ldots,f_d \in L^\infty(\mu)$, the multiple ergodic
averages
\[\frac{1}{|I_N|}\sum_{n\in I_N}\prod_{i=1}^df_i\circ T_i^n\]
converge in $L^2(\mu)$ as $N\to\infty$.}

\vspace{7pt}

While the sequence of works preceding the proof of convergence in
the one-dimensional setting of the Multiple Recurrence Theorem
develops a large body of ergodic-theoretic machinery for the
analysis of these averages, Tao departs quite markedly from those
approaches and effectively converts the problem of convergence into
a quantitative assertion concerning averages of $[-1,1]$-valued
functions on large finite grids $\{1,2,\ldots,N\}^d$.

A new proof of Tao's Theorem was given using classical
ergodic-theoretic machinery in~\cite{Aus--nonconv}. It turns out
that this convergence can be proved relatively quickly using a
version of the older approaches, with the one new twist that
starting from a system of commuting transformations of interest
$T_1,T_2,\ldots,T_d \curvearrowright (X,\S,\mu)$ one must first pass
to a carefully-chosen \emph{extended} system
$\tilde{T}_1,\tilde{T}_2,\ldots,\tilde{T}_d \curvearrowright
(\tilde{X},\tilde{\S},\tilde{\mu})$ (that is, a new system for which
the original one is isomorphic to the action of the $\tilde{T}_i$'s
on some globally invariant $\s$-subalgebra of $\tilde{\S}$: in
ergodic-theoretic terms, the original system is a `factor' of the
new one).  If the extension is constructed correctly then the
asymptotic behaviour of the multiple ergodic averages associated to
it admits a simplification allowing them to be compared with a
similar system of averages involving only $k-1$ transformations;
from this point convergence in $L^2$ follows quickly by induction on
$k$. The need for this extension also offers some explanation for
the advantage that Tao gains in his approach to Theorem C by
converting to the finitary, combinatorial world: during the course
of his proof he constructs new functions from the initial data of
the problem in ways that cannot be used to construct
\emph{measurable} functions in the ergodic-theoretic setting, but
suitable measurable functions are available using the larger
$\s$-algebra of the extended system.

Theorem C proves the convergence of the scalar averages appearing in
Theorem A because
\[\frac{1}{N}\sum_{n=1}^N\mu(T_1^{-n}(A)\cap T_2^{-n}(A)\cap \cdots\cap T_d^{-n}(A)) = \int_X \frac{1}{N}\sum_{n = 1}^N\prod_{i=1}^d(f_i\circ T_i^n)\,\d\mu\]
when $f_1 = f_2 = \ldots = f_d = 1_A$. Note that another re-proof of
Tao's theorem involving non-standard analysis has been given by
Towsner in~\cite{Tow09}, and that a different construction of some
extensions of probability-preserving systems that can be used as in
the proof of~\cite{Aus--nonconv} has since been given by Host
in~\cite{Hos09}.

Having found the extended systems appearing in the new proof of
Theorem C, it turns out that they also afford a somewhat simplified
description of the limiting value of the scalar averages appearing
in Theorem A.  These limiting values can always be expressed in
terms of a certain $(d+1)$-fold self-joining of the system
$(\tilde{X},\tilde{\S},\tilde{\mu},\tilde{T}_1,\tilde{T}_2,\ldots,\tilde{T}_d)$
(which appears already in the works of Furstenberg and Katznelson),
and one finds that for the extended system this self-joining takes a
special form.  Crucially, that special form is precisely the
hypothesis required to apply another result of Tao: the infinitary
analog of the hypergraph removal lemma from~\cite{Tao07}.  This
leads fairly quickly to a new proof of Theorem A (and hence also
one-dimensional multiple recurrence and their combinatorial
consequences), which appeared in~\cite{Aus--newmultiSzem}.

A similar story is now known in the setting of Theorem B.  For their
proof of that theorem, Furstenberg and Katznelson first provided a
correspondence with a class of stochastic processes enjoying
stationarity with respect to some semigroup of transformations. This
is broadly similar to Furstenberg's original correspondence between
Szemer\'edi's Theorem and the Multiple Recurrence Theorem, but
differs considerably in its details.  Having built this bridge to a
class of stochastic processes, Furstenberg and Katznelson then used
analogs of their earlier structural results from the setting of
probability-preserving $\bbZ^d$-actions to prove the `recurrence'
result that is the translation of Theorem B.  Here, too, it turns
out that the strategy of seeking extended systems in which the
behaviour of interest is simplified leads to a new proof of that
recurrence result, and so overall to a considerably shortened proof
of Theorem B, where again the punchline is an implementation of
Tao's infinitary hypergraph removal.  This new proof of Theorem B
appears in~\cite{Aus--DHJ}.  It was discovered simultaneously with
the work of the Polymath project~\cite{DHJ09}, which provided the
first finitary, effective proof of that theorem, and the proof
of~\cite{Aus--DHJ} used a key construction discovered by the members
of that project (again, suitably translated to apply to the
stochastic processes).

More recently still, in pursuit of some convergence results for `polynomial' analogs of the functional averages of Theorem C, it was found that a very abstract, unified approach could be given to the construction of the different extensions underlying the above-mentioned proofs of Theorems C, A and B.  This rests on the notion of a system that is `sated' relative to another class of systems.  In this dissertation, the new proofs of the above results are re-told using this unifying language, and some speculations offered concerning some further extensions of this machinery.

\vspace{7pt}

\subsubsection*{Outline of the following chapters}

In the next chapter we recall some basic definitions and conventions
from the study of measurable dynamical systems, and then introduce
the chief technical innovation on which most of the remaining
chapters will rest: a special property of certain dynamical systems
called `satedness'.  The main result of that chapter,
Theorem~\ref{thm:sateds-exist}, asserts that any
probability-preserving dynamical system admits extensions that enjoy
this `satedness' (where precisely what this means is relative to a
choice of another class of systems).

In Chapter~\ref{chap:conv} we use the existence of sated extensions to prove Theorem C.  After the introduction of another important technical device, the `Furstenberg self-joining', this follows by a quick induction once the strategy of passing to a sated extension has been decided.

Chapter~\ref{chap:multiMRT} is dedicated to Theorem A.  In this case
the use of sated extensions gives a relatively easy reduction of the
proof to a case in which the Furstenberg self-joining (which
describes the limiting averages of interest) admits a rather
detailed structural description; but the use of that description to
deduce the desired positivity of these averages is still rather
involved.  This requires an implementation of (a very slight
modification of) Tao's `infinitary hypergraph removal lemma', which
we will recall for completeness.

In Chapter~\ref{chap:DHJ} we prove Theorem B.  This proof follows very closely that of Theorem A, notwithstanding that the category of dynamical systems in which the proof takes place is very different.  However, the unusual features of this new category will require that we quickly re-examine the existence of sated extensions proved in Chapter 2 to check that a slightly modified version of that result holds here. After recalling Furstenberg and Katznelson's original reformulation of Theorem B in terms of a `recurrence' property of certain `strongly stationary' stochastic processes, we establish this new notion of `coordinatewise-satedness' and show that in this world it implies a similar structure for certain joint distributions to that obtained for the Furstenberg self-joining in Chapter~\ref{chap:multiMRT}.  The proof of Theorem B is then completed by another appeal to infinitary hypergraph removal, essential identical to that in Chapter~\ref{chap:multiMRT}.

Finally, Chapter~\ref{chap:spec} contains some speculations around an important question left open by our work.  In the case of $\bbZ^d$-actions treated by Chapters~\ref{chap:conv} and~\ref{chap:multiMRT}, one can discern in the background a very general ergodic-theoretic meta-question concerning the possible joinings among systems enjoying various additional invariances.  This is formulated precisely in Section~\ref{sec:meta}, but in that section it is answered only in a special case that suffices for the proof of Theorem A. A more general answer would be very interesting in its own right, as well as potentially offering new insights on other generalizations of nonconventional average convergence and multiple recurrence.  In Chapter~\ref{chap:spec} we will formulate a conjecture that would answer this question much more completely.

\chapter{Setting the stage}\label{chap:basics}

A handful of key technical ideas in ergodic theory will drive all of the proofs in the later chapters of this work.  After recalling some standard definitions and notation in the first section below, we introduce two such key ideas: that of a subclass of a class of dynamical systems that has the property of being `idempotent', and the constructions that this assumption of idempotence enables; and then the possibility of a system being `sated' relative to such an idempotent class, together with the result that all systems have extensions that are sated in this way.

These preliminary sections provide the necessary background for Chapters 3 and 4 (and also Chapter 6).  Unfortunately, the slightly unusual class of stochastic processes that appears in Chapter 5 is a little less willing to be analysed using this standard framework: the key ideas of idempotence and satedness will be central there too, but only after being modified to suit that class.  The modifications will be explained early in that chapter, together with those small changes that must accordingly be made to the proofs in Sections~\ref{sec:idem} and~\ref{sec:sateds}.  In principle one could give a unified treatment of all of these settings, but only at the expense of working with quite abstractly-defined categories of dynamical system and operations on them, in which our basic intuitions for the notions recalled in Section~\ref{sec:basics} may become obscured.  Although more unified, that route seems to pose too great a risk to the clarity of the other chapters, and so we shall only indicate it in passing during Chapter 5.

\section{Probability-preserving systems}\label{sec:basics}

Throughout this paper $(X,\S)$ will denote a measurable space.
Since our main results pertain only to the joint distribution of
countably many bounded real-valued functions on this space and their
shifts under some measurable transformations, by passing to the image
measure on a suitable product space we may always assume that $(X,\S)$
is standard Borel, and this will prove convenient for some of our
later constructions. In addition, $\mu$ will always denote a
probability measure on $\S$. We shall write $(X^S,\S^{\otimes S})$
for the usual product measurable structure indexed by a set $S$, and
$\mu^{\otimes S}$ for the product measure and $\mu^{\Delta S}$ for
the diagonal measure on this structure respectively.  Given a measurable map $\phi:(X,\S)\to (Y,\Phi)$
to another measurable space, we shall write $\phi_\#\mu$ for the
resulting pushforward probability measure on $(Y,\Phi)$.

Suppose now that $\G$ is a discrete semigroup, and consider the class of all probability-preserving actions $T:\G\actson (X,\S,\mu)$ on standard Borel probability spaces; these will be referred to as \textbf{$\G$-systems}, and will often be denoted by either the quadruple $(X,\S,\mu,T)$ or simply by a boldface letter such as $\bfX$.  If
$\L\leq \G$ is a subgroup we denote by $T^{\
\uhr\L}$ the $\L$-action on $(X,\S,\mu)$ defined by $(T^{\
\uhr\L})^\g := T^\g$ for $\g \in \L$, and refer to this as the
\textbf{$\L$-subaction}, and if $\bfX = (X,\S,\mu,T)$ is a $\G$-system
then we write similarly $\bfX^{\ \uhr\L}$ for the system $(X,\S,\mu,T^{\
\uhr\L})$ and refer to it as a \textbf{subaction system}.

A $\G$-system $(X,\S,\mu,T)$ is \textbf{trivial} if $\mu$ is supported on a single point.  Since any two such systems are measure-theoretically isomorphic simply by identifying these single points, we will usually refer to `the' trivial system.

We will make repeated use of a handful of standard constructions and properties of $\G$-systems.

\subsection*{Factors and joinings}

A \textbf{factor} of the $\G$-system $(X,\S,\mu,T)$ is a globally
$T$-invariant $\s$-subalgebra $\Phi \leq \S$. Relatedly, a
\textbf{factor map} from one $\G$-system $T:\G\curvearrowright
(X,\S,\mu)$ to another $S:\G\curvearrowright (Y,\Phi,\nu)$ is a
measurable map $\pi:X \to Y$ such that $\nu = \pi_\#\mu$ and
$S^\g\circ\pi = \pi\circ T^\g$ for all $\g\in\G$.  This situation is
often signified by writing $\pi:(X,\S,\mu,T) \to (Y,\Phi,\nu,S)$.
Factor maps comprise the natural morphisms between systems for a
fixed acting semigroup.

To any factor map $\pi$ is associated the factor $\{\pi^{-1}(A):\ A\in\Phi\}\leq \S$.  Two
factor maps $\pi$ and $\psi$ are \textbf{equivalent} if these
$\s$-subalgebras of $\S$ that they generate are equal up to
$\mu$-negligible sets, in which case we shall write $\pi \simeq
\psi$; this clearly defines an equivalence relation among factors.

It is a standard fact that in the category of standard Borel spaces
equivalence classes of factors are in bijective correspondence with
equivalence classes of globally invariant $\s$-subalgebras under the
relation of equality modulo negligible sets. A treatment of these
classical issues may be found, for example, in Chapter 2 of
Glasner~\cite{Gla03}.  Given a globally invariant $\s$-subaglebra in
$\bfX$, a choice of factor $\pi:\bfX\to \bfY$ generating that
$\s$-subalgebra will be referred to as \textbf{coordinatizing} the $\s$-subalgebra.

More generally, the factor map $\pi:(X,\S,\mu,T) \to (Y,\Phi,\nu,S)$
\textbf{contains} $\psi:(X,\S,\mu,T)\to (Z,\Psi,\theta,R)$ if
$\pi^{-1}(\Phi) \supseteq \psi^{-1}(\Psi)$ up to $\mu$-negligible
sets.  Another standard feature of standard Borel
spaces is that this inclusion is equivalent to the existence of a
\textbf{factorizing} factor map $\phi:(Y,\Phi,\nu,S) \to (Z,\Psi,\theta,R)$
with $\psi = \phi\circ\pi$ $\mu$-a.s., and that a measurable analog
of the Schroeder-Bernstein Theorem holds: $\pi \simeq \psi$ if and
only if a single such $\phi$ may be chosen that is invertible away
from some negligible subsets of the domain and target. If $\pi$
contains $\psi$ we shall write $\pi \succsim \psi$ or $\psi \precsim \pi$.

If $\pi:\bfX\to\bfY$ and $\psi:\bfX\to \bfZ$ are any two factor maps
as above (not necessarily ordered), then the $\s$-subalgebra
$\pi^{-1}(\Phi)\vee \psi^{-1}(\Psi)$ is another factor of $\bfX$. In
general we will write $\pi\vee \psi$ for an arbitrary choice of
factor map coordinatizing this factor, and similarly for larger
collections of factor maps.

Dual to the idea of a factor is that of an extension: if $\bfX$ is a $\G$-system, then an \textbf{extension} $\bfX$ is another $\G$-system $\t{\bfX}$ together with a factor map $\pi:\t{\bfX}\to \bfX$.

More general than the notion of a factor is that of a joining: if $\bfX_1$, $\bfX_2$, \ldots, $\bfX_k$ are $\G$-systems then a \textbf{joining} of them is another $\G$-system $\bfX$ together with factor maps $\pi_i:\bfX\to \bfX_i$ such that these $\pi_i$ together generate the whole $\s$-algebra of $\bfX$.  Since their introduction by Furstenberg in~\cite{Fur67}, joinings have become one of the most important concepts in the ergodic theorist's vocabulary, as is well-demonstrated in Glasner's book~\cite{Gla03}.

\subsection*{Partially invariant factors}

Given a $\G$-system $\bfX = (X,\S,\mu,T)$, the $\s$-algebra $\S^T$
of sets $A\in\S$ for which $\mu(A\triangle T^\g(A))=0$ for all $\g
\in \G$ is $T$-invariant, so defines a factor of $\bfX$. More
generally, if $\G$ is a group and $\L \unlhd\G$ then we can consider
the $\s$-algebra $\S^{T\uhr\L}$ generated by all $T^{\
\uhr\L}$-invariant sets: we refer to this as the
\textbf{$\L$-partially invariant factor}.  Note that in this case
the condition that $\L$ be normal is needed for this to be a
globally $T$-invariant factor. Similarly, if $S\subseteq \G$ and
$\L$ is the normal subgroup generated by $S$, we will sometimes
write $\S^{T\uhr S}$ for $\S^{T\uhr\L}$.

If moreover $\G$ is Abelian and $T_1$ and
$T_2$ are two commuting actions of $\G$ on $(X,\S,\mu)$,
then we can define a third action $T_1T_2^{-1}$ by setting
$(T_1T_2^{-1})^\g := T_1^\g T_2^{\g^{-1}}$. Given this we often write $\S^{T_1 = T_2}$ in place of $\S^{T_1^{-1}T_2}$, and
similarly for a larger number of actions of the same
group.

\subsection*{Relative independence}

If $\S_i \geq \Xi_i$ are factors of $(X,\S,\mu,T)$ for each $i\leq d$,
then the tuple of factors $(\S_1,\S_2,\ldots,\S_d)$ is
\textbf{relatively independent} over the tuple
$(\Xi_1,\Xi_2,\ldots,\Xi_d)$ if whenever $f_i\in
L^\infty(\mu)$ is $\S_i$-measurable for each $i\leq d$ we have
\[\int_X\prod_{i\leq d}f_i\,\d\mu = \int_X\prod_{i\leq d}\sfE_\mu(f_i\,|\,\Xi_i)\,\d\mu.\]
The information that various joint distributions are relatively independent will repeatedly prove pivotal in the following. Sometimes for brevity we will write that `$\S_1$ is relatively independent from $\S_2$, $\S_d$, \ldots, $\S_d$ over $\Xi_1$' if $(\S_1,\S_2,\ldots,\S_d)$ is relatively independent over $(\Xi_1,\S_2,\ldots,\S_d)$.

In case $\G$ is a group (not just a semigroup, so each $T^\g$ is invertible) we can construct examples of this situation as follows.  Suppose that $\bfY = (Y,\Phi,\nu,S)$ is a $\G$-system and \[\pi_i:\bfX_i  = (X_i,\S_i,\mu_i,T_i)\to \bfY\] are extensions of it for $i=1,2,\ldots,k$.  Then the \textbf{relatively independent product} of the systems $\bfX_i$ over their factor maps $\pi_i$ is the system
\[\prod_{\{\pi_1 = \ldots = \pi_k\}}\bfX_i = \Big(\prod_{\{\pi_1 = \ldots = \pi_k\}}X_i,\bigotimes_{\{\pi_1 = \ldots = \pi_k\}}\S_i,\bigotimes_{\{\pi_1 = \ldots = \pi_k\}}\mu_i,T_1\times \cdots\times T_k\Big)\]
where
\[\prod_{\{\pi_1 = \ldots = \pi_k\}}X_i :=
\{(x_1,\ldots,x_k)\in X_1\times \cdots\times X_k:\\
\pi_1(x_1) = \ldots = \pi_k(x_k)\},\]
$\bigotimes_{\{\pi_1 = \ldots = \pi_k\}}\S_i$ is the restriction of $\S_1\otimes\cdots\otimes\S_k$ to this subset of $X_1\times\cdots\times X_k$, and
\[\bigotimes_{\{\pi_1 = \ldots = \pi_k\}}\mu_i := \int_Y \bigotimes_{i=1}^k\mu_{i,y}\,\nu(\d y)\]
with $y\mapsto \mu_{i,y}$ an arbitrary choice of disintegration of $\mu_i$ over $\pi_i$. A quick check shows that the factors generated by the coordinate projections $\phi_j:\prod_{\{\pi_1 = \ldots = \pi_k\}}\bfX_i\to \bfX_j$ are relatively independent over the common further factor map
\[\pi_1\circ\phi_1\simeq \ldots \simeq \pi_k\circ\phi_k:\prod_{\{\pi_1 = \ldots = \pi_k\}}\bfX_i\to \bfY.\]

In case
$k=2$ we write the relatively independent product more simply as $\bfX_1\times_{\{\pi_1=
\pi_2\}}\bfX_2$, and in addition if $\bfX_1 = \bfX_2 = \bfX$ and
$\pi_1 = \pi_2 = \pi$ then we will abbreviate this further to
$\bfX\times_\pi\bfX$, and similarly for the individual spaces and
measures.

The need for the invertibility of $T$ in this construction arises in checking that $\bigotimes_{\{\pi_1 = \ldots = \pi_k\}}\mu_i$ is invariant under the product action.  For example, if $k=2$ then the invariance of $\mu_i$ under $T_i$ implies that for each $\g \in \G$ the disintegrations $\mu_{i,y}$ satisfy
\[\int_Y (T_i^\g)_\#\mu_{i,y}\,\nu(\d y) = \int_Y \mu_{i,y}\,\nu(\d y).\]
However, to argue from here to the invariance of $\mu_1\otimes_{\{\pi_1 = \pi_2\}}\mu_2$ we must know in addition that for $\nu$-almost every $y \in Y$ there is a unique point $S^{\g^{-1}}(y) \in Y$ such that $(T_i^\g)_\#\mu_{i,S^{\g^{-1}}(y)}$ is supported on the fibre over $y$.  Given this and the essential uniqueness of disintegrations, the above equation implies that $(T_i^\g)_\#\mu_{i,y} = \mu_{i,S^\g(y)}$ for $\nu$-almost every $y$, from which it also follows that
\[(T_1\times T_2)^\g_\#(\mu_{1,y}\otimes \mu_{2,y}) = (\mu_{1,S^\g(y)}\otimes \mu_{2,S^\g(y)})\]
$\nu$-almost surely, so that integrating again with respect to $y$ gives the desired invariance of $\mu_1\otimes_{\{\pi_1 = \pi_2\}}\mu_2$.  However, this latter argument is valid only if we can obtain the above equality pointwise in $y$, and this can fail if $T_i^\g$ is not invertible.

\subsection*{Inverse limits}

An \textbf{inverse sequence} of $\G$-systems is a family of
$\G$-systems $(X_m,\S_m,\mu_m,T_m)$
together with factor maps
\[\psi^m_k:(X_m,\S_m,\mu_m,T_m)\to
(X_k,\S_k,\mu_k,T_k)\quad\quad\hbox{for all}\ m\geq k\]
satisfying the compatibility property that $\psi^k_\ell\circ \psi^m_k = \psi^m_\ell$ whenever $m\geq k\geq \ell$. From such a family one can construct
an \textbf{inverse limit}
\[\lim_{m\leftarrow}\,\big((X_m,\S_m,\mu_m,T_m)_m,(\psi^m_k)_{m\geq k}\big) =: (X,\S,\mu,T)\]
together with a sequence of factor maps
\[\psi_m:(X,\S,\mu,T)\to (X_m,\S_m,\mu_m,T_m)\]
such that $\psi^m_k\circ \psi_m = \psi_k$ whenever $m\geq k$, and such that the lifted factors $\psi_m^{-1}(\S_m)$ together generate the whole of $\S$. Moreover, subject to these stipulations this inverse limit is unique up to isomorphisms that intertwine all the factor maps $\psi_m$. This construction is
described, for example, in Section 6.3 of Glasner~\cite{Gla03}.

\section{Idempotent classes}\label{sec:idem}

In much of the following we will be concerned with properties of one system that are defined relative to some other class of systems.

\begin{dfn}[Idempotent class] A subclass $\sfC$ of $\G$-systems is \textbf{idempotent} if it contains the trivial system and is closed under measure-theoretic isomorphism, inverse limits and joinings.
\end{dfn}

Note that our `classes' need not be sets in the sense of ZFC.  In all subsequent constructions involving these classes it will be clear that we need only some set-indexed family of members, and so we will not generally pass comment on this set-theoretic distinction.  Alternatively, we could circumvent this issue altogether by working only with probability-preserving systems modelled by some Borel transformations and invariant probability measure on, say, the Cantor space, since any standard Borel system admits such a model up to measure-theoretic isomorphism (see, for instance, Theorem 2.15 in~\cite{Gla03}).

\noindent\textbf{Examples}\quad Suppose that $\G$ is a group and that $\L \unlhd \G$.  Then the class of all $\G$-systems for which the subaction of $\L$ is trivial is easily seen to be idempotent.  This important example will usually be denoted by $\sfZ_0^{\L}$ in the following.

More generally, for $\L$ as above and any $n\in \bbN$
we let $\sfZ_n^\L$ denote the class of systems on which the
$\L$-subaction is a distal tower of height at most $n$, in the sense
of direct integrals of compact homogeneous space data introduced
in~\cite{Aus--ergdirint} to allow for the case of non-ergodic
systems.  Standard results on the possible joinings and inverse limits of isometric extensions show that this class is idempotent (see~\cite{Aus--ergdirint,Aus--lindeppleasant1}).  Those arguments also allow us to identify certain natural idempotent subclasses of $\sfZ_n^\L$, such as the class $\sfZ_{\rm{Ab},n}^\L$ of those systems with $\L$-subaction a
distal tower of height at most $n$ and in which each isometric
extension is Abelian. \fin

\begin{lem}
If $\sfC$ is an idempotent class of $\G$-systems then any $\G$-system $\bfX$ has an essentially unique maximal factor
in the class $\sfC$.
\end{lem}

\noindent\textbf{Proof}\quad It is clear that under the above assumption the
family of factors
\[\big\{\Xi \leq \S:\ \Xi\ \hbox{is generated by a factor map to a system in $\sfC$}\big\}\]
is nonempty (it contains $\{\emptyset,X\}$, which corresponds to the trivial system), upwards directed
(because $\sfC$ is closed under joinings) and closed under taking
$\s$-algebra completions of increasing unions (because $\sfC$ is
closed under inverse limits).  There is therefore a maximal
$\s$-subalgebra in this family. \qed

\begin{dfn}
If $\sfC$ is an idempotent class then $\bfX$ is a
\textbf{$\sfC$-system} if $\bfX \in \sfC$, and for any $\bfX$ we
write $\zeta_\sfC^\bfX:\bfX\to\sfC\bfX$ for an arbitrarily-chosen
coordinatization of its \textbf{maximal $\sfC$-factor} given by the
above lemma.

It is clear that if $\pi:\bfX\to\bfY$ then $\zeta_\sfC^\bfX \succsim
\zeta_\sfC^\bfY\circ\pi$, and so there is an essentially unique
factorizing map, which we denote by $\sfC\pi$, that makes the
following diagram commute:
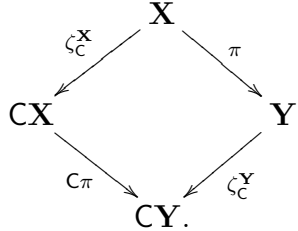
\begin{center}
$\phantom{i}$\xymatrix{
&\bfX\ar[dl]_{\zeta^\bfX_\sfC}\ar[dr]^{\pi}\\
\sfC\bfX\ar[dr]_{\sfC\pi} & & \bfY\ar[dl]^{\zeta^\bfY_\sfC}\\
&\sfC\bfY.}
\end{center}

In addition, we shall abbreviate $\bfX\times_{\zeta_\sfC^\bfX} \bfX$
to $\bfX\times_\sfC \bfX$, and similarly for the individual spaces
and measures defining this relatively independent product.
\end{dfn}

The above lemma and definition explain the choice of the term `idempotent', which is motivated by a more categorial viewpoint of such subclasses: if we identify such a class $\sfC$ with a full subcategory of the category of $\G$-systems with factor maps as morphisms, then the assignments $\bfX\mapsto \sfC\bfX$, $\pi \mapsto \sfC\pi$ define an autofunctor of this category which is idempotent.

The name we give for our next definition is also motivated by this
relationship with functors.

\begin{dfn}[Order continuity]
A class of $\G$-systems $\sfC$ is \textbf{order continuous} if
whenever $(\bfX_m)_{m\geq 0}$, $(\psi^m_k)_{m\geq k\geq
0}$ is an inverse sequence of $\G$-systems with inverse limit
$\bfX$, $(\psi_m)_{m\geq 0}$ we have
\[\zeta_\sfC^\bfX = \bigvee_{m\geq 0}\zeta_\sfC^{\bfX_m}\circ \psi_m:\]
that is, the maximal $\sfC$-factor of the inverse limit is simply
given by the (increasing) join of the maximal $\sfC$-factors of the
contributing systems.
\end{dfn}

\noindent\textbf{Example}\quad Although all the idempotent classes that will
matter to us later can be shown to be order continuous, it
may be instructive to exhibit one that is not.  In case $\G$ is an Abelian group, let us say that a
system $\bfX$ has a \textbf{finite-dimensional Kronecker factor} if
its Kronecker factor $\zeta_1^\bfX:X\to Z_1^\bfX$ can be
coordinatized as a direct integral (se Section 3
of~\cite{Aus--ergdirint}) of rotations on some measurably-varying
compact Abelian groups all of which can be isomorphically embedded
into a fibre repository $\bbT^D$ for some fixed $D \in \bbN$ (this
includes the possibility that the Kronecker factor is finite or
trivial). It is now easy to check that the class of $\bbZ$-systems
comprising all those that are either themselves finite-dimensional
Kronecker systems, or have a Kronecker factor that is \emph{not}
finite-dimensional (so we exclude just those systems that have a
finite-dimensional Kronecker factor but properly contain it), is
idempotent but not order continuous, since any infinite-dimensional
separable group rotation can be identified with an inverse limit of
finite-dimensional group rotations. \fin

\begin{dfn}[Hereditariness]
An idempotent class $\sfC$ is \textbf{hereditary} if it is also closed
under taking factors.
\end{dfn}

\begin{dfn}[Join]
If $\sfC_1$, $\sfC_2$ are idempotent classes, then the class
$\sfC_1\vee \sfC_2$ of all joinings of members of $\sfC_1$ and
$\sfC_2$ is clearly also idempotent. We call $\sfC_1\vee \sfC_2$ the
\textbf{join} of $\sfC_1$ and $\sfC_2$.
\end{dfn}

\begin{lem}[Join preserves order continuity]
If $\sfC_1$ and $\sfC_2$ are both order continuous then so is
$\sfC_1\vee \sfC_2$.
\end{lem}

\noindent\textbf{Proof}\quad Let $(\bfX_m)_{m\geq 0}$,
$(\psi^m_k)_{m\geq k \geq 0}$ be an inverse sequence with inverse
limit $\bfX$, $(\psi_m)_{m\geq 0}$.  Then $\zeta^\bfX_{\sfC_1\vee
\sfC_2}$ is the maximal factor of $\bfX$ that is a joining of a
$\sfC_1$-factor and a $\sfC_2$-factor (so, in particular, it must be
generated by its own $\sfC_1$- and $\sfC_2$-factors), and hence it
is equivalent to $\zeta^\bfX_{\sfC_1} \vee \zeta^\bfX_{\sfC_2}$.
Therefore any $f \in L^\infty(\mu)$ that is $\zeta^\bfX_{\sfC_1}
\vee \zeta^\bfX_{\sfC_2}$-measurable can be approximated in
$L^2(\mu)$ by some function of the finite-sum form $\sum_p
g_{p,1}\cdot g_{p,2}$ with each $g_{p,i} \in L^\infty(\mu)$ being
$\sfC_i$-measurable, and now since each $\sfC_i$ is order continuous
we may further approximate each $g_{p,i}$ by some $h_{p,i} \circ
\psi_m$ for a large integer $m$ and some $\sfC_i$-measurable
$h_{p,i} \in L^\infty(\mu_m)$. Combining these approximations
completes the proof. \qed

\noindent\textbf{Examples}\quad Of course, we can form the joins of
any of our earlier examples of idempotent classes: for example,
given a group $\G$ and subgroups $\L_1,\L_2,\ldots,\L_n \unlhd \G$
we can form $\sfZ_0^{\L_1}\vee \sfZ_0^{\L_2}\vee\cdots\vee
\sfZ_0^{\L_n}$. This particular example and several others like it
will appear frequently throughout the rest of this work. Clearly
each class $\sfZ_0^\L$ is hereditary, but in general joins of
several such classes are not; we will see this explicitly in the
first example of the next section. \fin

The following terminology will also prove useful.

\begin{dfn}[Joining to an idempotent class; adjoining]
If $\bfX$ is a system and $\sfC$ is an idempotent class then a
\textbf{joining of $\bfX$ to $\sfC$} or a \textbf{$\sfC$-adjoining
of $\bfX$} is a joining of $\bfX$ and $\bfY$ for some $\bfY \in
\sfC$.
\end{dfn}

\section{Sated systems}\label{sec:sateds}

The remainder of this dissertation concerns the consequences of one basic idea: that by extending a probability-preserving system, it is sometimes possible to impose on it some additional structure that makes its behaviour more transparent.  For our later applications, a notion of `additional structure' that is both useful and obtainable is best summarized by demanding that the system does not admit a nontrivial joining to systems drawn from various other special classes.  We will soon show that all systems admit extensions for which some version of this is true.  This idea, although very abstract and very simple, will repeatedly
prove surprisingly powerful.

\begin{dfn}[Sated system]\label{dfn:sated}
Given an idempotent class $\sfC$, a system $\bfX$ is
\textbf{$\sfC$-sated} if whenever $\pi:\t{\bfX} =
(\t{X},\t{\S},\t{\mu},\t{T}) \to \bfX$ is an extension, the factor maps
$\pi$ and $\zeta^{\t{\bfX}}_\sfC$ on $\t{X}$ are relatively
independent over $\zeta^\bfX_\sfC\circ\pi =
\sfC\pi\circ\zeta^{\t{\bfX}}_\sfC$ under $\t{\mu}$.  Phrased more pictorially, the two systems in the middle row of the commutative diagram
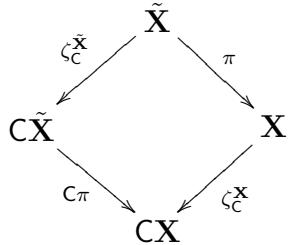
\begin{center}
$\phantom{i}$\xymatrix{
&\t{\bfX}\ar[dl]_{\zeta^{\t{\bfX}}_\sfC}\ar[dr]^{\pi}\\
\sfC\t{\bfX}\ar[dr]_{\sfC\pi} & & \bfX\ar[dl]^{\zeta^\bfX_\sfC}\\
&\sfC\bfX}
\end{center}
are relatively independent over their common factor copy of the system $\sfC\bfX$.

An inverse sequence is \textbf{$\sfC$-sated} if it has a cofinal
subsequence all of whose systems are $\sfC$-sated.
\end{dfn}

\noindent\textbf{Remark}\quad This definition has an important precedent in
Furstenberg and Weiss' notion of a `pair homomorphism' betwen
extensions elaborated in Section 8 of~\cite{FurWei96}. \fin

\noindent\textbf{Example}\quad If $\bfX = (U,\rm{Borel},\rm{Haar},R_\phi)$ with $U$ a compact metrizable Abelian group, $\phi:\bbZ^2\to U$ a dense homomorphism and $R_\phi$ the corresponding action of $\bbZ^2$ by rotations (so $R^\bf{n}(z) := z + \phi(\bf{n})$), then $\sfZ_0^{\bf{e}_i}\bfX$ is coordinatized by the quotient homomorphism
\[U\to U/\ol{\phi(\bbZ\bf{e}_i)}, \]
and so $\bfX$ is a member of $\sfZ_0^{\bf{e}_1}\vee \sfZ_0^{\bf{e}_2}$ if and only if these quotients together generate the whole of $U$, hence if and only if $\ol{\phi(\bbZ\bf{e}_1)}\cap\ol{\phi(\bbZ\bf{e}_2)}=\{0\}$.

On the other hand, any ergodic action $\bfX$ of $\bbZ^2$ by compact group rotations can be extended to a member of $\sfZ_0^{\bf{e}_1}\vee \sfZ_0^{\bf{e}_2}$. To see this we first note that ergodicity is equivalent to the denseness of $\phi(\bbZ^2)$ in $U$, and so in particular that $\ol{\phi(\bbZ\bf{e}_1)} + \ol{\phi(\bbZ\bf{e}_2)} = U$.  It follows that the `larger' group rotation system
\[\t{\bfX} = (\t{U},\rm{Borel},\rm{Haar},R_{\t{\phi}}),\]
where $\t{U} := \ol{\phi(\bbZ\bf{e}_1)} \oplus \ol{\phi(\bbZ\bf{e}_2)}$ and the homomorphism $\t{\phi}:\bbZ^2\to U^2$ is defined by
\[\t{\phi}(\bf{e}_1) := (\phi(\bf{e}_1),0)\quad\hbox{and}\quad \t{\phi}(\bf{e}_2) := (0,\phi(\bf{e}_2)),\] is an extension of $\bfX$ through the factor map
\[\t{U}\to U:(x,y)\mapsto x+y.\]
Now $\t{\bfX}$ clearly satisfies the above condition for membership of $\sfZ_0^{\bf{e}_1}\vee \sfZ_0^{\bf{e}_2}$, since the quotients by $\ol{\t{\phi}(\bbZ\bf{e}_i)}$ for $i=1,2$ are respectively the second and first coordinate projections.  It follows that every such $\bfX$ admits a $(\sfZ_0^{\bf{e}_1}\vee \sfZ_0^{\bf{e}_2})$-adjoining that generates the whole of $\bfX$, and which is therefore not relatively independent over any proper factor of $\bfX$, and hence that $\bfX$ itself is $(\sfZ_0^{\bf{e}_1}\vee \sfZ_0^{\bf{e}_2})$-sated if and only if it is already in the class $\sfZ_0^{\bf{e}_1}\vee \sfZ_0^{\bf{e}_2}$.  This reasoning also shows that the class $\sfZ_0^{\bf{e}_1}\vee \sfZ_0^{\bf{e}_2}$ is not hereditary.

A little more generally, if $\bfX$ is a totally weakly mixing extension of an ergodic action $\bfY$ of $\bbZ^2$ by compact group rotations, then routine arguments show that $\bfX$ is $(\sfZ_0^{\bf{e}_1}\vee \sfZ_0^{\bf{e}_2})$-sated if and only if this is true of $\bfY$ (since a totally weakly mixing extension is relatively disjoint from any $\sfZ_0^{\bf{e}_1}$-system, and given this the Furstenberg-Zimmer Inverse Theorem implies that the $\bf{e}_2$-invariant factor of any $\sfZ_0^{\bf{e}_1}$-adjoining of $\bfX$ is also relatively independent from $\bfX$ over its factor map to $\bfY$; see, for instance, Chapters 9 and 10 of~\cite{Gla03}).  Therefore such an $\bfX$ is $(\sfZ_0^{\bf{e}_1}\vee \sfZ_0^{\bf{e}_2})$-sated if and only if $\bfY \in \sfZ_0^{\bf{e}_1}\vee \sfZ_0^{\bf{e}_2}$.
\fin

The crucial technical fact that turns satedness into a useful tool
is the ability to construct sated extensions of arbitrary systems.
This can be seen as a natural abstraction from Propositions 4.6
of~\cite{Aus--nonconv} and 4.3 of~\cite{Aus--newmultiSzem}, and appears in its full strength as Theorem 3.11 in~\cite{Aus--lindeppleasant1}.

\begin{thm}[Idempotent classes admit multiply sated
extensions]\label{thm:sateds-exist} If $(\sfC_i)_{i\in I}$ is a
countable family of idempotent classes then any system $\bfX_0$
admits an extension $\pi:\bfX \to \bfX_0$ such that
\begin{itemize}
\item $\bfX$ is $\sfC_i$-sated for every $i\in I$;
\item the factors $\pi$ and $\bigvee_{i\in I}\zeta^{\bfX}_{\sfC_i}$ generate the whole of
$\bfX$.
\end{itemize}
\end{thm}

We shall prove this result after a preliminary lemma.

\begin{lem}\label{lem:inverse-lim-sated}
If $\sfC$ is an idempotent class then the inverse limit of any
$\sfC$-sated inverse sequence is $\sfC$-sated.
\end{lem}

\noindent\textbf{Proof}\quad By passing to a subsequence if necessary, it
suffices to suppose that $(\bfX_m)_{m\geq 0}$,
$(\psi^m_k)_{m\geq k\geq 0}$ is an inverse sequence of
$\sfC$-sated systems with inverse limit $\bfX_\infty$,
$(\psi_m)_{m\geq 1}$, and let $\pi:\t{\bfX} \to \bfX_\infty$
be any further extension and $f \in L^\infty(\mu_\infty)$.  We
will commit the abuse of identifying such a function with its lift
to any given extension when the extension in question is obvious.
With this in mind, we need to show that
\[\sfE(f\,|\,\zeta^{\t{\bfX}}_\sfC) = \sfE(f\,|\,\zeta^{\bfX_\infty}_\sfC).\]

However, by the $\sfC$-satedness of each $\bfX_m$, we certainly
have
\[\sfE(\sfE(f\,|\,\psi_m)\,|\,\zeta^{\t{\bfX}}_\sfC) = \sfE(f\,|\,\zeta^{\bfX_m}_\sfC),\]
and now as $m\to\infty$ this equation converges in $L^2(\mu)$ to
\[\sfE(f\,|\,\zeta^{\t{\bfX}}_\sfC) = \sfE\Big(f\,\Big|\,\bigvee_{m\geq 1}\,(\zeta^{\bfX_m}_\sfC\circ\psi_m)\Big).\]
By monotonicity we have
\[\zeta^{\t{\bfX}}_\sfC\succsim\zeta^{\bfX_\infty}_\sfC \succsim \bigvee_{m\geq 1}\,(\zeta^{\bfX_m}_\sfC\circ \psi_m),\]
and so by sandwiching the desired equality of conditional
expectations must also hold. \qed

\noindent\textbf{Proof of Theorem~\ref{thm:sateds-exist}}\quad We first prove
this for $I$ a singleton, and then in the general case.

\textbf{Step 1}\quad Suppose that $I = \{i\}$ and $\sfC_i =
\sfC$. This case will follow from a simple `energy increment'
argument.

Let $(f_r)_{r\geq 1}$ be a countable subset of the $L^\infty$-unit
ball $\{f\in L^\infty(\mu):\ \|f\|_\infty\leq 1\}$ that is dense in
this ball for the $L^2$-norm, and let $(r_i)_{i\geq 1}$ be a member
of $\bbN^\bbN$ in which every non-negative integer appears
infinitely often.

We will construct an inverse sequence $(\bfX_m)_{m\geq 0}$,
$(\psi^m_k)_{m\geq k \geq 0}$ starting from $\bfX_0$ such that each $\bfX_{m+1}$ is a $\sfC$-adjoining of
$\bfX_m$.  Suppose that for some $m_1 \geq 0$ we have already
obtained $(\bfX_m)_{m=0}^{m_1}$, $(\psi^m_k)_{m_1 \geq
m\geq k\geq 0}$ such that $\id_{X_{m_1}} \simeq
\zeta^{\bfX_{m_1}}_\sfC\vee \psi^{m_1}_{0}$.  We consider two
separate cases:
\begin{itemize}
\item If there is some further extension $\pi:\t{\bfX}\to
\bfX_{m_1}$ such that
\[\|\sfE_{\t{\mu}}(f_{r_{m_1}}\circ \psi^{m_1}_{0}\circ\pi\,|\,\zeta_\sfC^{\t{\bfX}})\|_2^2 > \|\sfE_{\mu_{m_1}}(f_{r_{m_1}}\circ \psi^{m_1}_0\,|\,\zeta_\sfC^{\bfX_{m_1}})\|_2^2 + 2^{-m_1},\]
then choose a particular $\pi:\t{\bfX}\to \bfX_{m_1}$ such that
the increase
\[\|\sfE_{\t{\mu}}(f_{r_{m_1}}\circ \psi^{m_1}_0\circ\pi\,|\,\zeta_\sfC^{\t{\bfX}})\|_2^2 - \|\sfE_{\mu_{m_1}}(f_{r_{m_1}}\circ \psi^{m_1}_0\,|\,\zeta_\sfC^{\bfX_{m_1}})\|_2^2\]
is at least half its supremal possible value over all extensions. By
restricting to the possibly smaller subextension of
$\t{\bfX}\to\bfX_{m_1}$ generated by $\pi$ and
$\zeta_\sfC^{\t{\bfX}}$ we may assume that $\t{\bfX}$ is itself a
$\sfC$-adjoining of $\bfX_{m_1}$ and hence of $\bfX_0$, and now we
let $\bfX_{m_1+1} := \t{\bfX}$ and $\psi^{m_1 + 1}_{m_1} :=
\pi$ (the other connecting factor maps being determined by this
one).
\item If, on the other hand, for every further extension $\pi:\t{\bfX}\to
\bfX_{m_1}$ we have \[\|\sfE_{\t{\mu}}(f_{r_{m_1}}\circ
\psi^{m_1}_0\circ\pi\,|\,\zeta_\sfC^{\t{\bfX}})\|_2^2 \leq
\|\sfE_{\mu_{m_1}}(f_{r_{m_1}}\circ
\psi^{m_1}_0\,|\,\zeta_\sfC^{\bfX_{m_1}})\|_2^2 + 2^{-m_1}\]
then we simply set $\bfX_{m_1+1} := \bfX_{m_1}$ and
$\psi^{m_1+1}_{m_1} := \id_{X_{m_1}}$.
\end{itemize}

Finally, let $\bfX_\infty$, $(\psi_m)_{m \geq 0}$ be the
inverse limit of this sequence.  We have
\begin{multline*}
\id_{X_\infty} \simeq \bigvee_{m\geq 0}\psi_m \simeq \bigvee_{m\geq
0}(\zeta_\sfC^{\bfX_m}\vee \psi^m_0)\circ\psi_m\\ \simeq
\bigvee_{m\geq 0}(\zeta_\sfC^{\bfX_m}\circ\psi_m)\vee\bigvee_{m\geq
0}( \psi^m_0\circ\psi_m) \precsim \zeta_\sfC^{\bfX_\infty}\vee
\psi_0,
\end{multline*}
so $\bfX_\infty$ is still a $\sfC$-adjoining of $\bfX_0$.  To
show that it is $\sfC$-sated, let $\pi:\t{\bfX}\to\bfX_\infty$
be any further extension, and suppose that $f \in
L^\infty(\mu_\infty)$. We will complete the proof for Step 1 by
showing that
\[\sfE_{\t{\mu}}(f\circ\pi\,|\,\zeta_\sfC^{\t{\bfX}}) =
\sfE_{\mu_\infty}(f\,|\,\zeta_\sfC^{\bfX_\infty})\circ\pi.\]

Since $\bfX_\infty$ is a $\sfC$-adjoining of $\bfX$, this $f$
may be approximated arbitrarily well in $L^2(\mu_\infty)$ by
finite sums of the form $\sum_p g_p\cdot h_p$ with $g_p$ being
bounded and $\zeta_\sfC^{\bfX_\infty}$-measurable and $h_p$
being bounded and $\psi_0$-measurable, and now by density we may
also restrict to using $h_p$ that are each a scalar multiple of some
$f_{r_p}\circ\psi_0$, so by continuity and multilinearity it
suffices to prove the above equality for just one such product $g\cdot
(f_r\circ\psi_0)$.  Since $g$ is
$\zeta_\sfC^{\bfX_\infty}$-measurable, this requirement now
reduces to
\[\sfE_{\t{\mu}}(f_r\circ\psi_0\circ\pi\,|\,\zeta_\sfC^{\t{\bfX}}) =
\sfE_{\mu_\infty}(f_r\circ\psi_0\,|\,\zeta_\sfC^{\bfX_\infty})\circ\pi.\]
Since $\zeta_\sfC^{\t{\bfX}} \succsim
\zeta_\sfC^{\bfX_\infty}\circ\pi$, this will follow if we only
show that
\[\|\sfE_{\t{\mu}}(f_r\circ\psi_0\circ\pi\,|\,\zeta_\sfC^{\t{\bfX}})\|_2^2 =
\|\sfE_{\mu_\infty}(f_r\circ\psi_0\,|\,\zeta_\sfC^{\bfX_\infty})\|_2^2.\]

Now, by the martingale convergence theorem we have
\[\|\sfE_{\mu_m}(f_r\circ\psi^m_0\,|\,\zeta_\sfC^{\bfX_m})\|_2^2 \to \|\sfE_{\mu_\infty}(f_r\circ\psi_0\,|\,\zeta_\sfC^{\bfX_\infty})\|_2^2\]
as $m\to\infty$.  It follows that if
\[\|\sfE_{\t{\mu}}(f_r\circ\psi_0\circ\pi\,|\,\zeta_\sfC^{\t{\bfX}})\|_2^2
>
\|\sfE_{\mu_\infty}(f_r\circ\psi_0\,|\,\zeta_\sfC^{\bfX_\infty})\|_2^2\]
then for some sufficiently large $m$ we would have $r_m = r$ (since
each integer appears infinitely often as some $r_m$) but also
\begin{eqnarray*}
&&\|\sfE_{\mu_{m+1}}(f_r\circ\psi^{m+1}_0\,|\,\zeta_\sfC^{\bfX_{m+1}})\|_2^2
-
\|\sfE_{\mu_m}(f_r\circ\psi^m_0\,|\,\zeta_\sfC^{\bfX_m})\|_2^2\\
&& \leq
\|\sfE_{\mu_\infty}(f_r\circ\psi_0\,|\,\zeta_\sfC^{\bfX_\infty})\|_2^2
-
\|\sfE_{\mu_m}(f_r\circ\psi^m_0\,|\,\zeta_\sfC^{\bfX_m})\|_2^2\\
&& <
\frac{1}{2}\Big(\|\sfE_{\t{\mu}}(f_r\circ\psi_0\circ\pi\,|\,\zeta_\sfC^{\t{\bfX}})\|_2^2
-
\|\sfE_{\mu_m}(f\circ\psi^m_0\,|\,\zeta_\sfC^{\bfX_m})\|_2^2\Big)
\end{eqnarray*}
and
\[\|\sfE_{\t{\mu}}(f_r\circ\psi_0\circ\pi\,|\,\zeta_\sfC^{\t{\bfX}})\|_2^2\geq
\|\sfE_{\mu_m}(f\circ\psi^m_0\,|\,\zeta_\sfC^{\bfX_m})\|_2^2
+ 2^{-m},\] so contradicting our choice of $\bfX_{m+1}\to
\bfX_m$ in the first alternative in our construction above. This
contradiction shows that we must actually have the equality of
$L^2$-norms required.

\textbf{Step 2}\quad The general case follows easily from Step 1 and
a second inverse limit construction: choose a sequence $(i_m)_{m\geq
1} \in I^\bbN$ in which each member of $I$ appears infinitely often,
and form an inverse sequence $(\bfX_m)_{m\geq 0}$,
$(\psi^m_k)_{m\geq k\geq 0}$ starting from $\bfX_0$ such that each
$\bfX_m$ is $\sfC_{i_m}$-sated for $m\geq 1$. The inverse limit
$\bfX$ is now sated for every $\sfC_i$, by
Lemma~\ref{lem:inverse-lim-sated}. \qed

\noindent\textbf{Remark}\quad Thierry de la Rue has shown me another proof of
Theorem~\ref{thm:sateds-exist} in case $\G$ is a group that follows very quickly from ideas
contained in his paper~\cite{LesRitdelaRue03} with Lesigne and
Rittaud, and which has now received a nice separate writeup
in~\cite{delaRue09}. The key observation is that
\begin{center}
\emph{An idempotent class $\sfC$ is hereditary if and only if
every system is $\sfC$-sated.}
\end{center}
This in turn follows from a striking result of Lema\'nczyk, Parreau
and Thouvenot~\cite{LemParTho00} that if two systems $\bfX$ and
$\bfY$ are not disjoint then $\bfX$ shares a nontrivial factor with
the infinite Cartesian power $\bfY^{\times \infty}$.  Given now an
idempotent class $\sfC$ and a system $\bfX$, let $\sfC^\ast$ be the
hereditary idempotent class of all factors of members of $\sfC$, and
let $\bfY$ be any $\sfC$-system admitting a factor map $\pi:\bfY\to
\sfC^\ast\bfX$ (such exists because by definition $\sfC^\ast\bfX$ is
a factor of some $\sfC$-system).  Now forming $\t{\bfX}:=
\bfX\times_{\{\zeta_{\sfC^\ast}^\bfX = \pi\}}\bfY$ (so here is where we need $\G$ to be a group), a quick check
using the above fact shows that $\sfC\t{\bfX} = \sfC^\ast\t{\bfX}$,
and that this is equivalent to the $\sfC$-satedness of $\t{\bfX}$. \fin

\chapter{The convergence of nonconventional averages}\label{chap:conv}

In this chapter Theorem C will be deduced from Theorem~\ref{thm:sateds-exist}.  This amounts to a rather simpler outing for many of the same ideas that will go into proving recurrence in the next chapter.

We first recall the Hilbert space version of a classical estimate due to van der Corput, which has long been a workhorse of Ergodic Ramsey Theory.  After giving this its own section, the Furstenberg self-joining for a tuple of transformations is introduced, and then in the last section we show how the right instance of satedness implies that these enjoy some additional structure from which a proof of Theorem C follows quite quickly.

\subsection*{Notation}

Before commencing with any of these proofs, we make a slight modification to the notation of the Introduction to be more in keeping with that of Chapter~\ref{chap:basics}: rather than letting $T_1$, $T_2$, \ldots, $T_d$ denote a tuple of commuting individual transformations on $(X,\S,\mu)$, we henceforth regard these as the subactions of the basis vectors $\bf{e}_1$, $\bf{e}_2$, \ldots, $\bf{e}_d$ for a single $\bbZ^d$-action $T$. Theorem C is accordingly re-phrased as asserting that the averages
\[\frac{1}{N}\sum_{n=1}^N(f_1\circ T^{n\bf{e}_1})\cdot(f_2\circ T^{n\bf{e}_2})\cdot\cdots\cdot (f_d\circ T^{n\bf{e}_d})\]
converge in $L^2(\mu)$ for any $\bbZ^d$-system $(X,\S,\mu,T)$.  This slight increase in abstraction will prove worth tolerating when we come to various constructions of new actions from old during our later arguments, in which we will need to keep efficient track of how the action of one vector in $\bbZ^d$ may have been re-assigned to that of another.  It follows that in the remainder of this work, a list such as `$T_1$, $T_2$, \ldots, $T_d$' will denote a tuple of \emph{whole actions} of some previously-decided group, rather than individual transformations.

\section{The van der Corput estimate}

This result and a related discussion can be found, for example, as Theorem 2.2 of Bergelson~\cite{Ber96}.

\begin{prop}[Van der Corput estimate]\label{prop:vdC}
Suppose that $(u_n)_{n\geq 1}$ is a bounded sequence in a Hilbert space $\frH$.  If the vector-valued averages
\[\frac{1}{N}\sum_{n=1}^Nu_n\]
do not converge to $0$ in norm as $N\to\infty$, then also the scalar-valued averages
\[\frac{1}{M}\sum_{m=1}^M\frac{1}{N}\sum_{n=1}^N\langle u_n,u_{n+m}\rangle\]
do not converge to $0$ as $N\to\infty$ and then $M\to\infty$.
\end{prop}

\noindent\textbf{Proof}\quad For any fixed $H\geq 1$ we have
\[\frac{1}{N}\sum_{n=1}^Nu_n \sim \frac{1}{N}\sum_{n=1}^N \frac{1}{H}\sum_{h=1}^H u_{n+h}\]
as $N\to\infty$, where the notation $w_N\sim v_N$ denotes that $w_N
- v_N \to 0$ in $\frH$.  However, the squared norm of the right-hand
double average may be estimated by
\[\Big\|\frac{1}{N}\sum_{n=1}^N \frac{1}{H}\sum_{h=1}^H u_{n+h}\Big\|^2 \leq \frac{1}{N}\sum_{n=1}^N\Big\|\frac{1}{H}\sum_{h=1}^H u_{n+h}\Big\|^2\]
(using the triangle and Cauchy-Schwartz inequalities), and this right-hand side is equal to
\[\frac{1}{H^2}\sum_{h_1,h_2=1}^M\frac{1}{N}\sum_{n=1}^N \langle u_{n+h_1},u_{n+h_2}\rangle.\]
It follows that these averages must also not converge to $0$ as $N\to\infty$ and then $H\to\infty$; but for large $H$ these can be expressed as averages of the averages
\[\frac{1}{M}\sum_{m=1}^M\frac{1}{N}\sum_{n=1}^N\langle u_n,u_{n+m}\rangle\]
for correspondingly large values of $M$, and so these also cannot converge to $0$ as $N\to\infty$ and then $M\to\infty$, as required. \qed

\section{The Furstenberg self-joining}\label{sec:Fberg}

Theorem C is proved by induction on $d$.  In the first instance,
this induction is enabled by a construction that is made possible
once convergence is known for a smaller number of transformations,
and which will also be central to the proof of Theorem A in the next
chapter.

Thus, suppose now that for some $d\geq 1$ the convergence of Theorem C is known for all tuples of at most $d-1$ commuting transformations (so this assumption is vacuous if $d=1$).  Let $\bfX = (X,\S,\mu,T)$ be a $\bbZ^d$-system, and let $A_1$, $A_2$, \ldots, $A_d \in \S$.  By integrating and using the invariance of $\mu$ under $T^{\bf{e}_1}$, our assumption applied to the transformations $T^{\bf{e}_2 - \bf{e}_1}$, \ldots, $T^{\bf{e}_d - \bf{e}_1}$ implies that the scalar averages
\begin{multline*}
\frac{1}{N}\sum_{n=1}^N\mu(T^{-n\bf{e}_1}(A_1)\cap T^{-n\bf{e}_2}(A_2)\cap\cdots\cap T^{-n\bf{e}_d}(A_d))\\
 = \int_X 1_{A_1}\cdot\Big(\frac{1}{N}\sum_{n=1}^N (1_{A_2}\circ T^{n(\bf{e}_2 - \bf{e}_1)})\cdot\cdots\cdot (1_{A_d}\circ T^{n(\bf{e}_d - \bf{e}_1)}) \Big)\,\d\mu
\end{multline*}
converge as $N\to\infty$.  Moreover, the limit takes the form $\mu^{\rm{F}}(A_1\times A_2\times \cdots\times A_d)$ for some probability $\mu^{\rm{F}}$ on $X^d$ that is invariant under the diagonal $\bbZ^d$-action defined by $(T^{\times d})^{\bf{n}} := T^\bf{n}\times T^\bf{n}\times \cdots\times T^\bf{n}$, simply because it is a limit of averages of the off-diagonal joinings
\[\int_X \delta_{(T^{n\bf{e}_1}(x),T^{n\bf{e}_2}(x),\ldots,T^{n\bf{e}_d}(x))}\,\mu(\d x)\quad\quad \hbox{for}\ n \in \bbN.\]

The $\bbZ^d$-system $\bfX^{\rm{F}} := (X^d,\S^{\otimes
d},\mu^{\rm{F}},T^{\times d})$ is therefore a $d$-fold self-joining
of $\bfX$ through the $d$ coordinate projections
$\pi_i:\bfX^{\rm{F}} \to \bfX$.  We refer to either $\mu^\rm{F}$ or
$\bfX^\rm{F}$ as the \textbf{Furstenberg self-joining} of $\bfX$.
Given functions $f_1$, $f_2$, \ldots, $f_d \in L^\infty(\mu)$, by
approximating each of them in $L^\infty$ using step functions we may
extend the above definition of $\mu^\rm{F}$ to the convergence
\[\frac{1}{N}\sum_{n=1}^N\int_X (f_1\circ T^{n\bf{e}_1})\cdot (f_2\circ T^{n\bf{e}_2})\cdot\cdots \cdot(f_d\circ T^{n\bf{e}_d})\,\d\mu\to \int_{X^d} f_1\otimes f_2\otimes \cdots\otimes f_d\,\d\mu^\rm{F}\]
as $N\to\infty$.

In addition to its invariance under $T^{\times d}$, the definition of $\mu^\rm{F}$ gives an additional invariance that will shortly prove crucial.

\begin{lem}\label{lem:Fberg-offdiag-invar}
Provided the limiting self-joining $\mu^{\rm{F}}$ exists, it is also invariant under the transformation $T^{\bf{e}_1}\times T^{\bf{e}_2}\times \cdots\times T^{\bf{e}_d}$.
\end{lem}

\noindent\textbf{Proof}\quad For any $A_1$, $A_2$, \ldots, $A_d \in \S$ we have
\begin{eqnarray*}
&&\mu^{\rm{F}}((T^{\bf{e}_1}\times T^{\bf{e}_2}\times \cdots\times T^{\bf{e}_d})^{-1}(A_1\times A_2\times \cdots\times A_d))\\ &&\quad = \lim_{n\to\infty}\frac{1}{N}\sum_{n=1}^N\mu(T^{-n\bf{e}_1}(T^{-\bf{e}_1}(A_1))\cap \cdots\cap T^{-n\bf{e}_d}(T^{-\bf{e}_d}(A_d)))\\
&&\quad= \lim_{n\to\infty}\frac{1}{N}\sum_{n=2}^{N+1}\mu(T^{-n\bf{e}_1}(A_1)\cap \cdots\cap T^{-n\bf{e}_d}(A_d))\\ &&\quad= \mu^{\rm{F}}(A_1\times A_2\times \cdots\times A_d),
\end{eqnarray*}
where the last equality follows because the discrete intervals  $\{1,2,\ldots,N\}$ and $\{2,3,\ldots,N+1\}$ asymptotically overlap in $1 - \rm{o}(1)$ of their lengths. \qed

It will be important to know that Furstenberg self-joinings behave
well under inverse limits.  The following is another immediate
consequence of the definition, and we omit the proof.

\begin{lem}\label{lem:Fberg-inv-lim}
If $(\bfX_m)_{m\geq 0}$, $(\psi^m_k)_{m\geq k\geq 0}$ is an inverse sequence with inverse limit $\bfX$, $(\psi_m)_{m\geq 0}$, then the Furstenberg self-joinings $\bfX_m^\rm{F}$ form an inverse sequence under the factor maps $(\psi^m_k)^{\times d}$ with inverse limit $\bfX^\rm{F}$, $(\psi^{\times d}_m)_{m\geq 0}$. \qed
\end{lem}

\section{The proof of convergence}

The final observation needed before we prove Theorem C is that satedness implies a certain inverse result for the situation in which the functional averages
\[S_N(f_1,f_2,\ldots,f_d) := \frac{1}{N}\sum_{n=1}^N (f_1\circ T^{n\bf{e}_1})\cdot (f_2\circ T^{n\bf{e}_2})\cdot\cdots \cdot (f_d\circ T^{n\bf{e}_d})\]
do not converge to $0$.

\begin{prop}\label{prop:sated-implies-pleasant}
Suppose that $\bfX$ is $\sfC$-sated for the idempotent class
\[\sfC := \sfZ_0^{\bf{e}_1}\vee \bigvee_{j=2}^d\sfZ_0^{\bf{e}_1 - \bf{e}_j}\]
and that $f_i \in L^\infty(\mu)$ for $i=1,2,\ldots,d$.  In addition, let $\Phi:= \S^{T^{\bf{e}_1}}\vee \bigvee_{j=2}^d\S^{T^{\bf{e}_1} = T^{\bf{e}_j}}$, so this is a factor of $\bfX$. If
\[S_N(f_1,f_2,\ldots,f_d)\not\to 0\]
as $N\to\infty$, then also $\sfE(f_1\,|\,\Phi) \neq 0$.
\end{prop}

\noindent\textbf{Remark}\quad In the terminology of~\cite{FurWei96}, which has since become standard in this area (and is roughly followed in~\cite{Aus--nonconv}), this asserts that for a $\sfC$-sated system $\bfX$ the factor $\Phi$ is \textbf{partially characteristic}. \fin

\noindent\textbf{Proof}\quad This rests on an appeal to the van der
Corput estimate followed by a re-interpretation of what it tells us.
Letting $u_n := (f_1\circ T^{n\bf{e}_1})\cdot (f_2\circ
T^{n\bf{e}_2})\cdot\cdots \cdot (f_d\circ T^{n\bf{e}_d})$,
Proposition~\ref{prop:vdC} and our assumption imply that the double
averages
\begin{multline*}
\frac{1}{M}\sum_{m=1}^M\frac{1}{N}\sum_{n=1}^N\langle u_n,u_{n+m}\rangle\\ = \frac{1}{M}\sum_{m=1}^M\frac{1}{N}\sum_{n=1}^N\int_X \big((f_1\circ T^{n\bf{e}_1})\cdot \cdots \cdot (f_d\circ T^{n\bf{e}_d})\big)\quad\quad\quad\quad\quad\quad\quad\quad\\ \quad\quad\quad\quad\quad\quad\quad\quad\quad\quad\quad\quad\cdot \big((\ol{f_1}\circ T^{(n+m)\bf{e}_1})\cdot \cdots \cdot (\ol{f_d}\circ T^{(n+m)\bf{e}_d})\big)\,\d\mu
\end{multline*}
do not tend to $0$ as $N\to\infty$ and then $M\to\infty$.  However, simply by re-arranging the individual functions and recalling the definition of $\mu^\rm{F}$, the limit in $N$ behaves as
\begin{multline*}
\frac{1}{M}\sum_{m=1}^M\frac{1}{N}\sum_{n=1}^N\int_X ((f_1\cdot(\ol{f_1}\circ T^{m\bf{e}_1}))\circ T^{n\bf{e}_1})\cdot\cdots \cdot ((f_d\cdot (\ol{f_d}\circ T^{m\bf{e}_d}))\circ T^{n\bf{e}_d})\,\d\mu\\ \to \frac{1}{M}\sum_{m=1}^M\int_{X^d}(f_1\cdot(\ol{f_1}\circ T^{m\bf{e}_d}))\otimes\cdots\otimes (f_d\cdot(\ol{f_d}\circ T^{m\bf{e}_d}))\,\d\mu^\rm{F}\\
= \frac{1}{M}\sum_{m=1}^M\int_{X^d}(f_1\otimes \cdots\otimes f_d)\cdot (\ol{f_1\otimes\cdots\otimes f_d}\circ (T^{\bf{e}_1}\times \cdots\times T^{\bf{e}_d})^m)\,\d\mu^\rm{F}.
\end{multline*}

Now, since Lemma~\ref{lem:Fberg-offdiag-invar} gives that $\mu^\rm{F}$ is invariant under $T^{\bf{e}_1}\times T^{\bf{e}_2}\times \cdots\times T^{\bf{e}_d}$, the classical mean ergodic theorem allows us to take the limit in $M$ to obtain
\[\int_{X^d}(f_1\otimes \cdots\otimes f_d)\cdot \sfE_{\mu^\rm{F}}\big(\ol{f_1\otimes\cdots\otimes f_d}\,\big|\,(\S^{\otimes d})^{T^{\bf{e}_1}\times \cdots\times T^{\bf{e}_d}}\big)\,\d\mu^\rm{F}.\]

Thus the van der Corput estimate tells us that this integral is non-zero.  The proof is completed simply by re-phrasing this conclusion slightly.  We have previously used $\mu^\rm{F}$ to define a $\bbZ^d$-system $\bfX^\rm{F}$, but in light of Lemma~\ref{lem:Fberg-offdiag-invar} we may alternatively use it to define a $\bbZ^d$-system $\t{\bfX}$ by setting
\[(\t{X},\t{\S},\t{\mu}) := (X^d,\S^{\otimes d},\mu^\rm{F}),\]
\[\t{T}^{\bf{e}_1} := T^{\bf{e}_1}\times T^{\bf{e}_2}\times \cdots\times T^{\bf{e}_d}\]
and
\[\t{T}^{\bf{e}_i} := (T^{\times d})^{\bf{e}_i}\quad\quad\hbox{for}\ i=2,3,\ldots,d\]
(thus, the basis direction $\bf{e}_1$ is treated differently from the others). With this definition the first coordinate projection $\pi_1:X^d\to X$ still defines a factor map of $\bbZ^d$-systems $\t{\bfX}\to \bfX$, because $\t{T}^\bf{n}$ does agree with $T^\bf{n}$ on the first coordinate in $X^d$ for every $\bf{n}$.  On the other hand, for $i=2,3,\ldots,d$ the function $f_i\circ\pi_i \in L^\infty(\mu^\rm{F})$ depends only on the $i^{\rm{th}}$ coordinate in $X^d$, and on this coordinate the transformations $\t{T}^{\bf{e}_1}$ and $\t{T}^{\bf{e}_i}$ agree, so that $f_i\circ \pi_i$ is $\t{T}^{\bf{e}_i - \bf{e}_1}$-invariant.  Thus the nonvanishing
\[\int_{X^d}(f_1\otimes \cdots\otimes f_d)\cdot \sfE_{\mu^\rm{F}}\big(\ol{f_1\otimes\cdots\otimes f_d}\,\big|\,\t{\S}^{\t{T}^{\bf{e}_1}}\big)\,\d\mu^\rm{F}\neq 0\]
asserts that the lifted function $f_1\circ \pi_1$ has a nontrivial inner product with a function that is a pointwise product of $\t{\S}^{\t{T}^{\bf{e}_1} = \t{T}^{\bf{e}_i}}$-measurable functions for $i=2,3,\ldots,d$ and the function $\sfE_{\mu^\rm{F}}\big(\ol{f_1\otimes\cdots\otimes f_d}\,\big|\,\t{\S}^{\t{T}^{\bf{e}_1}}\big)$, which is manifestly $\t{\S}^{\t{T}^{\bf{e}_1}}$-measurable.  Therefore $f_1\circ \pi_1$ has a nontrivial conditional expectation onto $\t{\S}^{\t{T}^{\bf{e}_1}}\vee\bigvee_{j=2}^d\t{\S}^{\t{T}^{\bf{e}_1} = \t{T}^{\bf{e}_j}}$, which is the $\s$-algebra generated by the factor map $\t{\bfX}\to \sfC\t{\bfX}$.  On the other hand, by $\sfC$-satedness $f_1\circ \pi_1$ must be relatively independent from this $\s$-algebra over $\Phi$, and so we also have $\sfE_\mu(f_1\,|\,\Phi) \neq 0$, as required. \qed

\noindent\textbf{Proof of Theorem C}\quad This proceeds by induction on $d$.  The case $d=1$ is the classical mean ergodic theorem, so suppose now that $d \geq 2$, that we know the result for all tuples of at most $d-1$ transformations and that we are given $T:\bbZ^d\actson (X,\S,\mu)$.

Let $\sfC$ be the class in Proposition~\ref{prop:sated-implies-pleasant}.  By Theorem~\ref{thm:sateds-exist} we may choose a $\sfC$-sated extension $\pi:\t{\bfX}\to \bfX$, and now since the corresponding inclusion $L^\infty(\mu)\subseteq L^\infty(\t{\mu})$ is an embedding of algebras that preserves the norms $\|\cdot\|_2$ it will suffice to prove convergence for the analogs of the averages $S_N$ associated to $\t{\bfX}$.  To lighten notation we henceforth assume that $\bfX$ itself is $\sfC$-sated.

Suppose that $f_1,f_2,\ldots,f_{d+1} \in L^\infty(\mu)$.  Letting $\Phi:= \S^{T^{\bf{e}_1}}\vee \bigvee_{j = 2}^d\S^{T^{\bf{e}_1} = T^{\bf{e}_j}}$, we see that the function $f_1 - \sfE(f_1\,|\,\Phi)$ has zero conditional expectation onto $\Phi$, and so by the multilinearity of $S_N$ and Proposition~\ref{prop:sated-implies-pleasant} we have that
\begin{multline*}
S_N(f_1,f_2,\ldots,f_d) - S_N(\sfE(f_1\,|\,\Phi),f_2,\ldots,f_d) \\ = S_N(f_1 - \sfE(f_1\,|\,\Phi),f_2,\ldots,f_d)\to 0
\end{multline*}
in $L^2(\mu)$ as $N\to\infty$. It therefore suffices to prove convergence with $f_1$ replaced by $\sfE(f_1\,|\,\Phi)$, or equivalently under the assumption that $f_1$ is $\Phi$-measurable.

However, this implies that $f_1$ may be approximated in $\|\cdot\|_2$ by finite sums of the form $\sum_p g_p\cdot h_{2,p} \cdot h_{3,p} \cdot\cdots \cdot h_{d,p}$ in which each $g_p$ is $T^{\bf{e}_1}$-invariant and each $h_{j,p}$ is $T^{\bf{e}_j - \bf{e}_1}$-invariant.  Since the operator
\[f_1\mapsto S_N(f_1,f_2,\ldots,f_d)\]
is linear and uniformly continuous in $L^2(\mu)$ for fixed bounded $f_2$, $f_3$, \ldots, $f_d$, it therefore suffices to prove convergence in case  $f_1$ is simply one such product, say $gh_2h_3\cdots h_d$.  For this function, however, we can re-arrange our averages as
\begin{multline*}
S_N(f_1,f_2,\ldots,f_d) = \frac{1}{N}\sum_{n=1}^N((gh_2h_3\cdots h_d)\circ T^{n\bf{e}_1})\cdot (f_2\circ T^{n\bf{e}_2})\cdot\cdots\cdot (f_d\circ T^{n\bf{e}_d})\\
= g\cdot \frac{1}{N}\sum_{n=1}^N ((f_2h_2)\circ T^{n\bf{e}_2})\cdot\cdots\cdot ((f_dh_d)\circ T^{n\bf{e}_d}) = gS_N(1_X,f_2h_2,\ldots,f_dh_d),
\end{multline*}
since $g\circ T^{n\bf{e}_1} = g$ and $h_j\circ T^{n\bf{e}_1} =
h_j\circ T^{n\bf{e}_j}$ for each $j=2,3,\ldots,d$. Now the averages
appearing on the right are uniformly bounded in $\|\cdot\|_\infty$
and involve only the $d-1$ transformations $T^{\bf{e}_2}$,
$T^{\bf{e}_3}$, \ldots, $T^{\bf{e}_d}$, and so the inductive
hypothesis gives their convergence in $\|\cdot\|_2$. Since
$\|g\|_\infty < \infty$ this gives also the convergence of the
left-hand averages in $\|\cdot\|_2$, as required. \qed

\noindent\textbf{Remark}\quad In fact the above proof gives a slight strengthening of Theorem C, in that the convergence is uniform
in the location of the interval of averaging: that is, the averages
\[\frac{1}{|I_N|}\sum_{n\in I_N}\prod_{i=1}^d f_i\circ T^{n\bf{e}_i}\]
converge in $L^2(\mu)$ for any sequence of increasingly long finite intervals $I_N \subset \bbZ$, and the limit does not depend on the choice of these intervals.  This result is treated in full in~\cite{Aus--nonconv}. \fin

\chapter{Multiple recurrence for commuting transformations}\label{chap:multiMRT}

In this chapter we deduce Theorem A from Theorem~\ref{thm:sateds-exist}.  Coupled with Furstenberg and Katznelson's correspondence principle from~\cite{FurKat78}, this gives a new proof of the Multidimensional Szemer\'edi Theorem, but we will not recount that correspondence here since it is already well-known from that paper and several subsequent accounts, such as those in the books~\cite{Fur81} of Furstenberg and~\cite{TaoVu06} of Tao and Vu.

After introducing a more convenient reformulation of Theorem A below, we first introduce a very general meta-question that covers most of the ergodic-theory we need. We then show how it specializes to give quite detailed information on the Furstenberg self-joining corresponding to a tuple of commuting transformations.  From this the proof of Theorem A follows by appealing to a version of Tao's infinitary hypergraph removal lemma.

We will continue the practice begun in the previous chapter of writing a tuple of commuting transformations as $T^{\bf{e}_1}$, $T^{\bf{e}_2}$, \ldots, $T^{\bf{e}_d}$ for some $\bbZ^d$-action $T$. The convergence result of the previous chapter implies that for any such $T^{\bf{e}_1}$, $T^{\bf{e}_2}$, \ldots, $T^{\bf{e}_d}$ the Furstenberg self-joining $\mu^\rm{F}$ of Section~\ref{sec:Fberg} exists.  Knowing this, Theorem A about the limit infima of scalar averages is a consequence of the following more general result:

\begin{thm}\label{thm:multirec2}
If $T:\bbZ^d\actson(X,\S,\mu)$, $\mu^\rm{F}$ denotes the Furstenberg self-joining of the transformations $T^{\bf{e}_1}$, $T^{\bf{e}_2}$, \ldots, $T^{\bf{e}_d}$ and $A_1$, $A_2$, \ldots, $A_d\in \S$ then
\[\mu^\rm{F}(A_1\times A_2\times \cdots\times A_d) = 0 \quad\quad\Rightarrow\quad\quad \mu(A_1\cap A_2\cap\cdots\cap A_d) =0.\]
\end{thm}

Inded, in case $A_i = A$ for each $A$ this assertion is precisely
the contrapositive of Theorem A.  However, the formulation of
Theorem~\ref{thm:multirec2} has the great advantage of allowing us
to manipulate the sets $A_i$ separately in setting up a proof by
induction.

\section{The question in the background}\label{sec:meta}

Having reformulated our goal in this chapter as Theorem~\ref{thm:multirec2}, it becomes clear that it is really an assertion about the joint distribution of the coordinate projections $\pi_i:X^d\to X$, $i=1,2,\ldots,d$ under $\mu^\rm{F}$.

By Lemma~\ref{lem:Fberg-offdiag-invar} $\mu^\rm{F}$ is an invariant measure for the action $\vec{T}$ of the larger group $\bbZ^{d+1}$ defined by setting
\[\vec{T}^{\,\uhr\, \bbZ^d\oplus \{\bs{0}\}} := T^{\times d}\quad\quad\hbox{and}\quad\quad \vec{T}^{\bf{e}_{d+1}} := T^{\bf{e}_1}\times T^{\bf{e}_2}\times \cdots\times T^{\bf{e}_d}.\]
Thus this defines a $\bbZ^{d+1}$-system $\vec{\bfX}$ in which the Furstenberg self-joining $\bfX^\rm{F}$ corresponds to the subaction of $\bbZ^d\oplus\{\bs{0}\}$.  The key to our proof is the observation that the coordinate projections $\pi_i$ now define factor maps of $\vec{\bfX}$ onto a collection of $\bbZ^{d+1}$-systems $\bfX_1$, $\bfX_2$, \ldots, $\bfX_d$ for each of which some one-dimensional subgroup of $\bbZ^{d+1}$ acts trivially: specifically, this is so with $\bfX_i = (X_i,\S_i,\mu_i,T_i)$ defined simply by `doubling up' the $\bbZ\bf{e}_i$-subaction of $T$:
\[(X_i,\S_i,\mu_i) := (X,\S,\mu),\quad T_i^{\,\uhr\,\bbZ^d\oplus\{\bs{0}\}} := T\quad\hbox{and}\quad T_i^{\bf{e}_{d+1}} := T^{\bf{e}_i}.\]
It follows immediately from these specifications that $\pi_i\circ \vec{T} = T_i\circ\pi_i$ and that $\bfX_i \in \sfZ_0^{\bf{e}_{d+1} - \bf{e}_i}$.

Having made these observations, our principal results on $\mu^\rm{F}$ will fall within the pattern of the following:
\begin{quote}
\textbf{Meta-question:}

Given subgroups $\G_1$, $\G_2$, \ldots, $\G_r \leq \bbZ^D$ and $\bbZ^D$-systems $(X_i,\S_i,\mu_i,T_i)$ for $i=1,2,\ldots,r$ such that $T_i^{\uhr\,\G_i}=\id$, what do these partial invariances imply about the possible joinings of these $\bbZ^D$-systems?
\end{quote}

The first stage in proving Theorem~\ref{thm:multirec2} will boil down to a handful of special cases of this question.  In this section we show that a partial answer covering all of the cases we need can be given quite easily, subject to an algebraic constraint on the subgroups $\G_i$ and an allowance to pass to extended systems.

First, it is instructive to understand the simple case $r=2$:

\begin{lem}\label{lem:two-fold-joinings}
If the systems $\bfX_i$ are $\G_i$-partially invariant for $i=1,2$, then any joining of them is relatively independent over their factors $\S_i^{T_i\uhr(\G_1 + \G_2)}$.
\end{lem}

\noindent\textbf{Proof}\quad Suppose $\pi_i:(Y,\Phi,\nu,S)\to (X_i,\S_i,\mu_i,T_i)$ is a joining of the two systems and consider subsets $A_i \in \S_i$.  In addition let $(F_N)_{N\geq 1}$ be a F\o lner sequence of subsets of $\G_1$.  Then the invariance of $\nu$ and the Mean Ergodic Theorem give
\begin{eqnarray*}
&&\nu(\pi_1^{-1}(A_1)\cap \pi_2^{-1}(A_2))\\ &&\quad = \lim_{N\to\infty}\frac{1}{|F_N|}\sum_{\g\in F_N}\int_Y (1_{A_1}\circ \pi_1)(1_{A_2}\circ T_2^{\g}\circ \pi_2)\,\d\nu\\
&&\quad = \lim_{N\to\infty}\int_Y (1_{A_1}\circ \pi_1)\Big(\Big(\frac{1}{|F_N|}\sum_{\g\in F_N}1_{A_2}\circ T_2^{\g}\Big)\circ \pi_2\Big)\,\d\nu\\ &&\quad = \int_Y (1_{A_1}\circ \pi_1)(\sfE_{\mu_2}(A_2\,|\,\S_2^{T_2\uhr\G_1})\circ\pi_2)\,\d\nu.
\end{eqnarray*}
Since $T_2^{\uhr\,\G_2} = \id$ the factor $\S_2^{T_2\uhr\G_1}$ consists of sets that are invariant under the whole group $\G_1 + \G_2$, and hence agrees with $\S_2^{T_2\uhr(\G_1 + \G_2)}$.  Arguing similarly with the roles of $\bfX_1$ and $\bfX_2$ reversed, this shows that that above is equal to
\[\int_Y (\sfE_{\mu_1}(A_1\,|\,\S_1^{T_1\uhr(\G_1 + \G_2)})\circ \pi_1)(\sfE_{\mu_2}(A_2\,|\,\S_2^{T_2\uhr(\G_1 + \G_2)})\circ\pi_2)\,\d\nu,\]
as required. \qed

For $r\geq 3$ we will not obtain an answer as complete as the above.
However, a natural generalization is available for certain special
tuples of subgroups, subject to the further provision that we may
replace the originally-given systems $\bfX_i$ with some extensions
of them.  The extensions, of course, will be sated extensions, and
for them the picture is given by the following.

\begin{thm}\label{thm:sateds-joint-dist}
Suppose that
\[\bbZ^D \cong \G_1 \oplus \G_2 \oplus \cdots \oplus \G_r \oplus \L\]
is a direct sum decomposition of $\bbZ^D$ into the subgroups $\G_i$ and some auxiliary subgroup $\L$, and that $\bfX_i \in \sfZ_0^{\G_i}$ for $i=1,2,\ldots,r$ are systems such that each $\bfX_i$ is $\sfC_i$-sated for
\[\sfC_i := \bigvee_{j\leq r,\,j\neq i}\sfZ_0^{\G_i + \G_j}.\]
Then for any joining $\pi_i:\bfY \to \bfX_i$, $i=1,2,\ldots,r$, the factors $\pi_i^{-1}(\S_i)$ are relatively independent over their further factors
\[\pi_i^{-1}\Big(\bigvee_{j\leq r,\,j\neq i}\S_i^{T_i\uhr(\G_i + \G_j)}\Big).\]
\end{thm}

\noindent\textbf{Proof}\quad This is a simple appeal to the definition of satedness.  We will show that $\pi_1^{-1}(\S_1)$ is relatively independent from $\bigvee_{j=2}^r\pi_j^{-1}(\S_j)$ over\\ $\pi_1^{-1}\big(\bigvee_{j=2}^r\S_1^{T_1\uhr(\G_1 + \G_j)}\big)$, the cases of the other factors being similar.

Let $\G := \G_2 \oplus \cdots \oplus \G_r \oplus \L\leq \bbZ^D$, so this complements $\G_1$ in $\bbZ^D$, and let $\bfY = (Y,\Phi,\nu,S)$.  From $S$ we may construct a new $\nu$-preserving $\bbZ^D$-action $S'$ by defining
\[(S')^{\bf{m} + \bf{n}} := S^{\bf{n}}\quad\quad\hbox{for all}\ \bf{m}\in \G_1,\ \bf{n} \in \G.\]

Let $\bfY' := (Y,\Phi,\nu,S')$, so manifestly $\bfY \in \sfZ_0^{\G_1}$. Similarly define the systems $\bfX'_i = (X_i,\S_i,\mu_i,T_i')$ for $i=2,3,\ldots,r$, so these also have trivial $\G_1$-subactions and hence in fact lie in the classes $\sfZ_0^{\G_1 + \G_i}$. Since $T_1^{\bf{m}} = \id_{X_1}$ for all $\bf{m} \in \G_1$ by assumption, we see that $\pi_1\circ S' = T_1\circ \pi_1$, so $\pi_1$ still defines a factor map $\bfY' \to \bfX_1$.  On the other hand, we also have
\[\pi_i\circ (S')^{\bf{m} + \bf{n} + \bf{p}} = \pi_i\circ S^{\bf{n} + \bf{p}} = T_i^{\bf{n} + \bf{p}}\circ \pi_i = T_i^{\bf{p}}\circ \pi_i = (T'_i)^{\bf{m} + \bf{n} + \bf{p}}\circ \pi_i\]
whenever $i=2,3,\ldots,r$ and $\bf{m} \in \G_1$, $\bf{n} \in \G_i$ and $\bf{p} \in \bigoplus_{j\neq 1,i}\G_j \oplus \L$.

Therefore $\pi_i$ is a factor map $\bfY' \to \bfX_i'$ for $i=2,3,\ldots,d$, and so $\bfY'$ is a joining of $\bfX_1$ with members of the classes $\sfZ_0^{\G_1 + \G_i}$ for $i=2,3,\ldots,r$: that is, $\bfY$ is a $\sfC_1$-adjoining of $\bfX_1$.  By the assumption of $\sfC_1$-satedness, it follows that this adjoining is relatively independent over the maximal $\sfC_1$-factor of $\bfX_1$, which equals $\bigvee_{j=2}^r\S_1^{T_1\uhr(\G_1 + \G_j)}$, as required. \qed

\noindent\textbf{Example}\quad Without the assumption of satedness,
more complicated phenomena can appear in the joint distribution of
three partially-invariant systems.  For example, let $(X,\S,\mu,T)$
be the $\bbZ^3$-system on the two-torus $\bbT^2$ with its Borel
$\s$-algebra and Haar measure defined by $T^{\bf{e}_1} :=
R_{(\a,0)}$, $T^{\bf{e}_2} := R_{(0,\a)}$ and $T^{\bf{e}_3} :=
R_{(\a,\a)}$, where $R_q$ denotes the rotation of $\bbT^2$ by an
element $q \in \bbT^2$ and we choose $\a \in \bbT$ irrational.  In
this case we have natural coordinatizations of the partially
invariant factors $\zeta_0^{T^{\bf{e}_i}}:X \to \bbT$ given by
\[\zeta_0^{T^{\bf{e}_1}}(t_1,t_2) = t_2,\quad \zeta_0^{T^{\bf{e}_2}}(t_1,t_2) = t_1\quad\hbox{and}\quad \zeta_0^{T^{\bf{e}_3}}(t_1,t_2) = t_1 - t_2.\]
It follows that in this example any two of $\S^{T^{\bf{e}_1}}$,
$\S^{T^{\bf{e}_2}}$ and $\S^{T^{\bf{e}_3}}$ are independent, but also
that any two of them generate the whole system (and so overall
independence fails).

In fact, it is possible to give a fairly complete answer to our meta-question in the case of any three $\bbZ$-subactions of some $\bbZ^D$-action, without the simplifying power of extending our systems.  However, that answer in general requires the handling of extensions of non-ergodic systems by measurably-varying compact homogeneous space data: it is contained in Theorem 1.1 of~\cite{Aus--ergdirint}, in which such extensions are studied in suitable generality.  The full formulation of that Theorem 1.1 is rather lengthy, and will not be repeated here; and it seems clear that matters will only become more convoluted for larger $r$. \fin

Theorem~\ref{thm:sateds-joint-dist} already suffices for the coming
applications, but it is natural to ask about more general
collections of subgroups $\G_i \leq \bbZ^D$.  In fact it is possible
to do slightly better than Theorem~\ref{thm:sateds-joint-dist} with
just a little extra effort: the same conclusion holds given only
that these subgroups are \textbf{linearly independent}, in the sense
that for any $\bf{n}_i \in \G_i$ we have
\[\bf{n}_1 + \bf{n}_2 + \cdots + \bf{n}_r = \bs{0}\quad \Rightarrow\quad \bf{n}_i = \bs{0}\ \forall i\leq r.\]
Indeed, given this linear independence, one can let $\Delta := \G_1 + \G_2 + \ldots + \G_r$ and now argue as in the above proof to deduce that the conclusion holds provided that $\bfX_1$ is $\sfC_1$-sated \emph{among all $\Delta$-systems}.  However, it is not quite obvious that this is the same as being $\sfC_1$-sated among $\bbZ^D$-systems. This turns out to be true, but it requires the key additional result that whenever $\Delta \leq \L$ are discrete Abelian groups, $\bfX$ is a $\L$-system and $\a:\bfY \to \bfX^{\uhr\Delta}$ is an extension of the $\Delta$-subaction, there is an extension of $\L$-systems $\b:\bfZ\to \bfX$ that fits into a commutative diagram
\begin{center}
$\phantom{i}$\xymatrix{ \bfZ^{\uhr \Delta}\ar[dr]\ar[rr]^\b & & \bfX^{\uhr \Delta}\\
& \bfY\ar[ur]_{\a}. }
\end{center}

The elementary but slightly messy proof of this can be found in Subsection 3.2 of~\cite{Aus--lindeppleasant1}.

What happens when there are linearly dependences among the subgroups $\G_1$, $\G_2$, \ldots, $\G_r$?  An answer to this question could have several applications to understanding multiple recurrence, but it is also clearly of broader interest in ergodic theory.  At present the picture remains unclear, but a number of recent works have provided answers in several further special cases, and in moments of optimism it now seems possible that a quite general extension of Theorem~\ref{thm:sateds-joint-dist} (using satedness relative to a much larger list of classes of system) may be available.  A more precise conjecture in this vein will be formulated in Chapter~\ref{chap:spec}.

\noindent\textbf{Remark}\quad Before leaving this section, it is
worth contrasting the feature seen above that linear independence is
helpful with previous works in this area. In the early study of
special cases of Theorems B or C it was generally found that the
analysis of powers of a single transformation (or correspondingly of
arithmetic progressions in $\bbZ$) revealed more usable structure
and was thus more tractable than the general case.  Of course,
Furstenberg's original Multiple Recurrence Theorem preceded Theorem
B; and the conclusion of Theorem C was known in many such
`one-dimensional' cases long before the general case was treated
(see~\cite{ConLes84,ConLes88.1,ConLes88.2,FurWei96,HosKra07,Zie07},
although we note that Conze and Lesigne did also treat a
two-dimensional case of Theorem C, and that in~\cite{Zha96} Zhang
extended this result to three dimensions subject to some additional
assumptions).

The same phenomenon is apparent in the search for finitary, quantitative approaches to Szemer\'edi's Theorem and its relatives. Indeed, a purely finitary proof of the Multidimensional Szemer\'edi Theorem appeared only recently in works of R\"odl and Skokan~\cite{RodSko04}, Nagle, R\"odl and Schacht~\cite{NagRodSch06} and Gowers~\cite{Gow07}, building on the development by those authors of sufficiently powerful hypergraph variants of Szemer\'edi's Regularity Lemma in graph theory.  Furthermore, the known bounds for how large $N_0$ must be taken in terms of $\delta$ and $k$ are far better for Szemer\'edi's Theorem than for its multidimensional generalization, owing to the powerful methods developed by Gowers in~\cite{Gow98,Gow01}, which extend Roth's proof for $k=3$ from~\cite{Rot53} and are much more efficient than the hypergraph regularity proofs. As yet these methods have resisted extension to the multidimensional setting, except in one two-dimensional case recently treated by Shkredov~\cite{Shk05}.  This story is discussed in much greater depth in Chapters 10 and 11 of~\cite{TaoVu06}.

Running counter to this trend, the value of linear independence for the present work is a consequence of our strategy of passing to extensions of probability-preserving systems.  Although such extensions can lose any a priori algebraic structure (such as being a $\bbZ^D$-action in which the transformations $T^{\bf{e}_i}$ are actually all powers of one fixed transformation), the various instances of satedness that it allows us to assume will furnish enough power to drive all of our subsequent proofs.  These instances of satedness will all be relative to joins of different classes of partially invariant systems, and, as illustrated by the above proof of Theorem~\ref{thm:sateds-joint-dist}, the usefulness of this kind of satedness will rely on the ability to construct new systems for which the corresponding subgroups behave in specified ways.  With this in mind it is natural that having those subgroups linearly independent removes a potential obstacle from these arguments, and that answering our meta-question for sated systems will be more difficult when the subgroups exhibit some linear dependences. \fin

\section{More on the Furstenberg self-joining}\label{sec:moreFberg}

We now return to the study of the Furstenberg self-joining $\mu^\rm{F}$ introduced in the previous chapter, with the goal of deriving a structure theorem for it as a consequence of Theorem~\ref{thm:sateds-joint-dist} in case $\bfX$ is sated with respect to enough difference classes.  In order to formulate this structure theorem, we first settle on some more bespoke notation.

In the following we shall make repeated reference to certain
factors assembled from the partially invariant factors of our $\bbZ^d$-action $T$, so we now give these factors their own names. They
will be indexed by subsets of $[d] := \{1,2,\ldots,d\}$, or more
generally by subfamilies of the collection $\binom{[d]}{\geq 2}$ of
all subsets of $[d]$ of size at least $2$.  On the whole, these
indexing subfamilies will be up-sets in $\binom{[d]}{\geq
2}$: $\I\subseteq \binom{[d]}{\geq 2}$ is an \textbf{up-set} if $u \in \I$ and
$[d]\supseteq v\supseteq u$ imply $v \in \I$. For example, given $e
\subseteq [d]$ we write $\langle e\rangle := \{u \in\binom{[d]}{\geq
2}:\ u\supseteq e\}$ (note the non-standard feature of our notation
that $e \in \langle e\rangle$ if and only if $|e| \geq 2$): up-sets
of this form are \textbf{principal}.  We will abbreviate
$\langle\{i\}\rangle$ to $\langle i\rangle$. It will also be helpful
to define the \textbf{depth} of a non-empty up-set $\I$ to be
$\min\{|e|:\ e\in \I\}$.

The corresponding factor for $e =
\{i_1,i_2,\ldots,i_k\}\subseteq [d]$ with $k\geq 2$ is
$\Phi_e := \S^{T^{\bf{e}_{i_1}} = T^{\bf{e}_{i_2}} = \ldots = T^{\bf{e}_{i_k}}}$, so this is the partially invariant factor for the $(k-1)$-dimensional subgroup
\[\bbZ(\bf{e}_{i_1} - \bf{e}_{i_2}) + \bbZ(\bf{e}_{i_1} - \bf{e}_{i_3}) + \cdots + \bbZ(\bf{e}_{i_1} - \bf{e}_{i_k}).\]
More generally, given a family
$\cal{A} \subseteq \binom{[d]}{\geq 2}$ we define
$\Phi_{\cal{A}} := \bigvee_{e\in\cal{A}}\Phi_e$.

From the ordering among the factors $\Phi_e$ it is clear that
$\Phi_{\cal{I}} = \Phi_{\cal{A}}$ whenever $\cal{A} \subseteq
\binom{[d]}{\geq 2}$ is a family that generates $\cal{I}$ as an
up-set, and in particular that $\Phi_e = \Phi_{\langle e\rangle}$ when $|e|\geq 2$.

We now return to the Furstenberg self-joining $\mu^\rm{F}$.  For $e = \{i_1 < i_2 < \ldots < i_k\} \subseteq [d]$ we write $\mu^\rm{F}_e$ for the Furstenberg self-joining of the transformations $T^{\bf{e}_{i_1}}$, $T^{\bf{e}_{i_2}}$, \ldots, $T^{\bf{e}_{i_k}}$:
\[\mu^{\rm{F}}_e(A_1\times \cdots\times A_k)
:= \lim_{N\to
\infty}\frac{1}{N}\sum_{n=1}^N\mu(T^{-n\bf{e}_{i_1}}(A_1)\cap \cdots\cap T^{-n\bf{e}_{i_k}}(A_k)),\]
so this clearly extends the definition of Section~\ref{sec:Fberg} in the sense that $\mu_{[d]}^\rm{F} = \mu^\rm{F}$.  Of course, we know the existence of each $\mu^\rm{F}_e$ by the results of the previous chapter.

We next record some simple properties of the family of self-joinings $\mu^\rm{F}_e$ for $e \subseteq [d]$.  Given subsets $e \subseteq e' \subseteq [d]$, in the following we write $\pi_e$ for the coordinate projection $X^{e'}\to X^e$, since the choice of $e'$ will always be clear from the context.

\begin{lem}\label{lem:Fberg-project}
If $e \subseteq e' \subseteq [d]$ then $(\pi_e)_{\#}\mu^\rm{F}_{e'} = \mu^\rm{F}_e$.
\end{lem}

\noindent\textbf{Proof}\quad This is immediate from the definition: if $e = \{i_1 < i_2 < \ldots < i_k\} \subseteq e' = \{j_1 < j_2 <
\ldots < j_l\}$ and
$A_{i_j} \in \S$ for each $j \leq k$ then
\[(\pi_e)_{\#}\mu^{\rm{F}}_{e'}(A_{i_1}\times \cdots\times A_{i_k})\\
:= \lim_{N\to
\infty}\frac{1}{N}\sum_{n=1}^N\mu(T^{-n\bf{e}_{j_1}}(B_{j_1})\cap
\cdots\cap T^{-n\bf{e}_{j_l}}(B_{j_l}))\]
where $B_j := A_j$ if $j \in e$ and $B_j := X$ otherwise; but then
this last average simplifies summand-by-summand directly to
\[\lim_{N\to \infty}\frac{1}{N}\sum_{n=1}^N\mu(T^{-n\bf{e}_{i_1}}(A_1)\cap \cap\cdots\cap T^{-n\bf{e}_{i_k}}(A_k))
=: \mu^{\rm{F}}_e(A_1\times \cdots\times A_k),\]
as required. \qed

\begin{lem}\label{lem:diag}
For any $e \subseteq [d]$ and $A \in \Phi_e$ we have
\[\mu^\rm{F}_e(\pi_i^{-1}(A)\triangle \pi_j^{-1}(A)) = 0\quad\quad\forall i,j \in e:\]
thus, the restriction
$\mu^{\rm{F}}_e\uhr_{\Phi_e^{\otimes e}}$ is just the diagonal
measure $(\mu\uhr_{\Phi_e})^{\Delta e}$.
\end{lem}

\noindent\textbf{Proof}\quad If $e = \{i_1 < i_2 < \ldots < i_k\}$ and $A_j
\in \Phi_e$ for each $j \leq k$ then by definition we have
\begin{eqnarray*}&&\mu^{\rm{F}}_e(A_1 \times A_2\times \cdots\times A_k)\\ &&\quad = \lim_{N\to \infty}\frac{1}{N}\sum_{n=1}^N\mu(T^{-n\bf{e}_{i_1}}(A_1)\cap T^{-n\bf{e}_{i_2}}(A_2)\cap\cdots\cap
T^{-n\bf{e}_{i_k}}(A_k))\\
&&\quad = \lim_{N\to
\infty}\frac{1}{N}\sum_{n=1}^N\mu(T^{-n\bf{e}_{i_1}}(A_1\cap
A_2\cap\cdots\cap A_k))\\
&& \quad = \mu(A_1\cap A_2\cap\cdots\cap A_k),
\end{eqnarray*}
as required. \qed

It follows from the last lemma that whenever $e \subseteq e'$ the
factors $\pi_i^{-1}(\Phi_e) \leq \S^{\otimes d}$ for $i\in e$ are
all equal up to $\mu^{\rm{F}}_{e'}$-negligible sets.  It will prove
helpful later to have a dedicated notation for these factors.

\begin{dfn}[Oblique copies]
For each $e \subseteq [d]$ we refer to the common
$\mu^{\rm{F}}_{[d]}$-completion of the $\s$-subalgebra
$\pi_i^{-1}(\Phi_e)$, $i \in e$, as the \textbf{oblique copy} of
$\Phi_e$, and denote it by $\Phi^{\rm{F}}_e$.  More generally we
shall refer to factors formed by repeatedly applying $\cap$ and
$\vee$ to such oblique copies as \textbf{oblique factors}.
\end{dfn}

We are now ready to derive the more nontrivial consequences we need from Theorem~\ref{thm:sateds-joint-dist}.  These will appear in two separate propositions.

\begin{prop}\label{prop:Fberg1}
For each pair $i\leq d$ let
\[\sfC_i := \bigvee_{j\leq d,\,j\neq i}\sfZ_0^{\bf{e}_i - \bf{e}_j}.\]
If $\bfX$ is $\sfC_i$-sated for each $i$ then the coordinate projections $\pi_i:X^d\to X$ are relatively independent under $\mu^\rm{F}$ over the further factors
\[\pi_i^{-1}\Big(\bigvee_{j\leq d,\,j\neq i}\S^{T^{\bf{e}_i} = T^{\bf{e}_j}}\Big) = \pi_i^{-1}(\Phi_{\langle i\rangle}).\]
\end{prop}

\noindent\textbf{Proof}\quad This follows by applying Theorem~\ref{thm:sateds-joint-dist} to the $\bbZ^{d+1}$-system $\vec{\bfX}$ introduced at the beginning of the previous section. Indeed, as explained there the coordinate projections $\pi_i:\vec{\bfX}\to \bfX_i$ witness that $\vec{\bfX}$ is a joining of the systems $\bfX_i \in \sfZ_0^{\bbZ(\bf{e}_{d+1} - \bf{e}_i)}$.

Let
\[\sfD_i := \bigvee_{j\leq d,\,j\neq i}\sfZ_0^{\bbZ(\bf{e}_i - \bf{e}_{d+1}) + \bbZ(\bf{e}_j - \bf{e}_{d+1})},\]
an idempotent class of $\bbZ^{d+1}$-systems.  Now the assumption that $\bfX$ is $\sfC_i$-sated as a $\bbZ^d$-system implies that $\bfX_i$ is $\sfD_i$-sated as a $\bbZ^{d+1}$-system.  Indeed, given any extension of $\bbZ^{d+1}$-systems $\pi:\bfY\to \bfX_i$ the subaction $(\sfD_i\bfY)^{\uhr(\bbZ^d\oplus\{\bs{0}\}}$ is clearly a member of the class $\sfC_i$, so the $\sfC_i$-satedness of $\bfX$ implies that $\pi$ is relatively independent from $\zeta_{\sfD_i}^\bfY:\bfY\to \sfD_i\bfY$ over its further factor map $\zeta_{\sfC_i}^{\bfX}$, which agrees with $\zeta_{\sfD_i}^{\bfX_i}$ because the whole of $\bfX_i$ is already $\bbZ(\bf{e}_{d+1} - \bf{e}_i)$-partially invariant.

Setting $\G_i := \bbZ(\bf{e}_i - \bf{e}_{d+1})$ for $i=1,2,\ldots,d$ and $\L := \bbZ\bf{e}_{d+1}$, these subgroups define a direct-sum decomposition of $\bbZ^{d+1}$. Therefore Theorem~\ref{thm:sateds-joint-dist} applies to tell us that the factors $\pi_i^{-1}(\S)$ are relatively independent under $\mu^\rm{F}$ over their further factors
\[\pi_i^{-1}\Big(\bigvee_{j\leq d,\,j\neq i}\S^{T_i\uhr(\bbZ(\bf{e}_i - \bf{e}_{d+1}) + \bbZ(\bf{e}_j - \bf{e}_{d+1}))}\Big) = \pi_i^{-1}(\Phi_{\langle i\rangle}),\]
as required. \qed

For our second application of Theorem~\ref{thm:sateds-joint-dist} we need a preparatory lemma.

\begin{lem}\label{lem:inherit-satedness}
If $\sfC \subseteq \sfD$ are idempotent classes of $\G$-systems for any discrete group $\G$ and $\bfX$ is $\sfC$-sated, then $\sfD\bfX$ is also $\sfC$-sated.
\end{lem}

\noindent\textbf{Proof}\quad If $\bfX$ is $\sfC$-sated and $\pi:\bfY \to \sfD\bfX$ is any extension, then the relatively independent product $\t{\bfX} := \bfX\times_{\{\zeta_\sfD^\bfX = \pi\}} \bfY$ is an extension of $\bfX$ through the first coordinate projection (it for the sake of using this relatively independent product that we need $\G$ to be a group).  Therefore by $\sfC$-satedness the factor map $\zeta_\sfC^{\t{\bfX}}$ is relatively independent from this coordinate projection over the further factor map $\zeta_\sfC^\bfX:\bfX\to\sfC\bfX$ of the latter, and so the same must be true of $\zeta_\sfC^\bfY$.  However, the factor map $\zeta_\sfC^\bfX$ is clearly contained in the factor map $\zeta_\sfD^\bfX$ since $\sfC \subseteq \sfD$, and so it must actually equal $\zeta_\sfC^{\sfD\bfX}\circ \zeta_\sfD^{\bfX}:\bfX\to \sfC(\sfD\bfX)$.  Hence $\pi$ is relatively independent from $\zeta_\sfC^\bfY$ over its further factor map $\zeta_\sfC^{\sfD\bfX}$, as required. \qed

\begin{prop}\label{prop:Fberg2}
For each subset $e = \{i_1,i_2,\ldots,i_k\}\subseteq [d]$ let
\[\sfC_e := \bigvee_{j\in [d]\setminus e}\sfZ_0^{\bbZ(\bf{e}_{i_1} - \bf{e}_{i_2}) + \cdots + \bbZ(\bf{e}_{i_1} - \bf{e}_{i_k}) + \bbZ(\bf{e}_{i_1} - \bf{e}_j)},\]
and suppose now that $\bfX$ is $\sfC_e$-sated for every $e$ (so this includes the assumption of the previous proposition when $e$ is a singleton).  Then under $\mu^\rm{F}$ the oblique factors have the property that $\Phi^\rm{F}_\I$ and $\Phi^\rm{F}_{\I'}$ are relatively independent over $\Phi^\rm{F}_{\I\cap \I'}$ for any up-sets $\I,\I' \subseteq \binom{[d]}{\geq 2}$.
\end{prop}

\noindent\textbf{Proof}\quad\textbf{Step 1}\quad First observe that
the result is trivial if $\I \supseteq \I'$, so now suppose that
$\I' = \langle e\rangle$ where $e$ is a maximal member of
$\binom{[d]}{\geq 2}\setminus \I$.  Let $\{a_1,a_2,\ldots,a_m\}$ be
the antichain of minimal elements of $\I$, so that $\Phi^\rm{F}_\I =
\bigvee_{l\leq m}\Phi^\rm{F}_{a_k}$. The maximality assumption on
$e$ implies that $e\cup \{j\}$ contains some $a_k$ for every $j \in
[d]\setminus e$, and so $\I\cap \I'$ is precisely the up-set
generated by these sets $e\cup \{j\}$ for $j\in[d]\setminus e$.  We
must therefore show that $\Phi^\rm{F}_e$ is relatively independent
from $\bigvee_{k\leq m}\Phi^\rm{F}_{a_k}$ under $\mu^\rm{F}$ over
the common factor $\bigvee_{j\in [d]\setminus
e}\Phi^\rm{F}_{e\cup\{j\}}$.

Observe also that since $e \not\in \I$ we can find some $j_k \in
a_k\setminus e$ for each $k\leq m$.  Moreover, each $j\in
[d]\setminus e$ must appear as some $j_k$ in this list, since it
appears at least for any $k$ for which $a_k \subseteq e \cup \{j\}$.

Now Lemma~\ref{lem:diag} implies that $\Phi^\rm{F}_{a_k}$ agrees with $\pi_{j_k}^{-1}(\Phi_{a_k})$ up to $\mu^\rm{F}$-negligible sets.  On the other hand, we clearly have $\pi_{j_k}^{-1}(\Phi_{a_k})\leq \pi_{j_k}^{-1}(\S)$, and so in fact it will suffice to show that $\Phi^\rm{F}_e$ is relatively independent from $\bigvee_{j\in [d]\setminus e}\pi_j^{-1}(\S)$ over $\bigvee_{j\in [d]\setminus e}\Phi^\rm{F}_{e\cup\{j\}}$.

This alteration of the problem is important because it provides the linear independence needed to apply Theorem~\ref{thm:sateds-joint-dist}.  Indeed, considering again the $\bbZ^{d+1}$-system $\vec{\bfX}$, in the present setting we see that the $\s$-subalgebras
\[\Phi^\rm{F}_e\quad\hbox{and}\quad \pi_j^{-1}(\S)\ \hbox{for}\ j\in[d]\setminus e\]
constitute a collection of factors of $\vec{\bfX}$ that are partially invariant under the subgroups
\[\G_e := \bbZ(\bf{e}_i - \bf{e}_{d+1}) + \sum_{\ell\in e\setminus \{i\}}\bbZ(\bf{e}_i - \bf{e}_\ell)\quad\hbox{and}\quad \G_j := \bbZ(\bf{e}_j - \bf{e}_{d+1})\ \hbox{for}\ j\in [d]\setminus e\]
respectively, where $i \in e$ is arbitrary.  On the one hand these subgroups can be inserted into a direct sum decomposition of $\bbZ^{d+1}$, and on the other we may argue just as in the proof of Proposition~\ref{prop:Fberg1} that the $\bbZ^{d+1}$-system defined by the factor $\Phi^\rm{F}_e$ is sated relative to the class $\bigvee_{j \in [d]\setminus e}\sfZ_0^{\G_e + \G_j}$, using our satedness assumption on $\bfX$ and Lemma~\ref{lem:inherit-satedness}. The conclusion therefore follows from Theorem~\ref{thm:sateds-joint-dist}.

\textbf{Step 2}\quad The general case can now be treated for fixed
$\I$ by induction on $\I'$. If $\cal{I}' \subseteq \cal{I}$ then the
result is clear, so now let $e$ be a minimal member of
$\cal{I}'\setminus\cal{I}$ of maximal size, and let $\cal{I}'' :=
\cal{I}'\setminus\{e\}$.  It will suffice to prove that if $F \in
L^\infty(\mu^{\rm{F}})$ is $\Phi^{\rm{F}}_{\cal{I}'}$-measurable
then
\[\sfE_{\mu^{\rm{F}}}(F\,|\,\Phi^{\rm{F}}_{\cal{I}}) = \sfE_{\mu^{\rm{F}}}(F\,|\,\Phi^{\rm{F}}_{\cal{I}\cap\cal{I}'}),\]
and furthermore, by an approximation in $\|\cdot\|_2$ by finite sums
of products, to do so only for $F$ that are of the form $F_1\cdot
F_2$ with $F_1$ and $F_2$ being bounded and respectively
$\Phi^{\rm{F}}_{\langle e \rangle}$- and
$\Phi^{\rm{F}}_{\cal{I}''}$-measurable. However, for such a product
we can write
\[\sfE_{\mu^{\rm{F}}}(F\,|\,\Phi^{\rm{F}}_{\cal{I}}) =
\sfE_{\mu^{\rm{F}}}\big(\sfE_{\mu^{\rm{F}}}(F\,|\,\Phi^{\rm{F}}_{\cal{I}\cup\cal{I}''})\,\big|\,\Phi^{\rm{F}}_{\cal{I}}\big)
=
\sfE_{\mu^{\rm{F}}}\big(\sfE_{\mu^{\rm{F}}}(F_1\,|\,\Phi^{\rm{F}}_{\cal{I}\cup\cal{I}''})\cdot
F_2\,\big|\,\Phi^{\rm{F}}_{\cal{I}}\big).\]
By Step 1 we have
\[\sfE_{\mu^{\rm{F}}}(F_1\,|\,\Phi^{\rm{F}}_{\cal{I}\cup\cal{I}''})
=
\sfE_{\mu^{\rm{F}}}(F_1\,|\,\Phi^{\rm{F}}_{(\cal{I}\cup\cal{I}'')\cap\langle
e\rangle}),\] and on the other hand $(\cal{I}\cup\cal{I}'')\cap\langle
e\rangle \subseteq \cal{I}''$ (because $\cal{I}''$ contains every
subset of $[d]$ that strictly includes $e$, since $\I'$ is an
up-set), so $(\cal{I}\cup\cal{I}'')\cap\langle
e\rangle = \I''\cap \langle e\rangle$ and therefore another appeal to Step 1 gives
\[\sfE_{\mu^{\rm{F}}}(F_1\,|\,\Phi^{\rm{F}}_{(\cal{I}\cup\cal{I}'')\cap\langle
e\rangle}) =
\sfE_{\mu^{\rm{F}}}(F_1\,|\,\Phi^{\rm{F}}_{\cal{I}''}).\]
Therefore the above expression for
$\sfE_{\mu^{\rm{F}}}(F_1F_2\,|\,\Phi^{\rm{F}}_{\cal{I}})$
simplifies to
\begin{multline*}
\sfE_{\mu^{\rm{F}}}\big(\sfE_{\mu^{\rm{F}}}(F_1\,|\,\Phi^{\rm{F}}_{\cal{I}''})\cdot
F_2\,\big|\,\Phi^{\rm{F}}_{\cal{I}}\big) =
\sfE_{\mu^{\rm{F}}}\big(\sfE_{\mu^{\rm{F}}}(F_1\cdot
F_2\,|\,\Phi^{\rm{F}}_{\cal{I}''})\,\big|\,\Phi^{\rm{F}}_{\cal{I}}\big)\\
=
\sfE_{\mu^{\rm{F}}}\big(\sfE_{\mu^{\rm{F}}}(F\,|\,\Phi^{\rm{F}}_{\cal{I}''})\,\big|\,\Phi^{\rm{F}}_{\cal{I}}\big)
= \sfE_{\mu^{\rm{F}}}(F\,|\,\Phi^{\rm{F}}_{\cal{I}\cap
\cal{I}''}) =
\sfE_{\mu^{\rm{F}}}(F\,|\,\Phi^{\rm{F}}_{\cal{I}\cap\cal{I}'}),
\end{multline*}
where the third equality follows by the inductive hypothesis applied to $\cal{I}''$ and $\cal{I}$. \qed

\section{Infinitary hypergraph removal and completion of the proof}

Propositions~\ref{prop:Fberg1} and~\ref{prop:Fberg2} tell us a great
deal about the structure of the probability measure $\mu^\rm{F}$ for
a system $\bfX$ that is sated relative to all the necessary classes
in terms of the partially-ordered family of factors
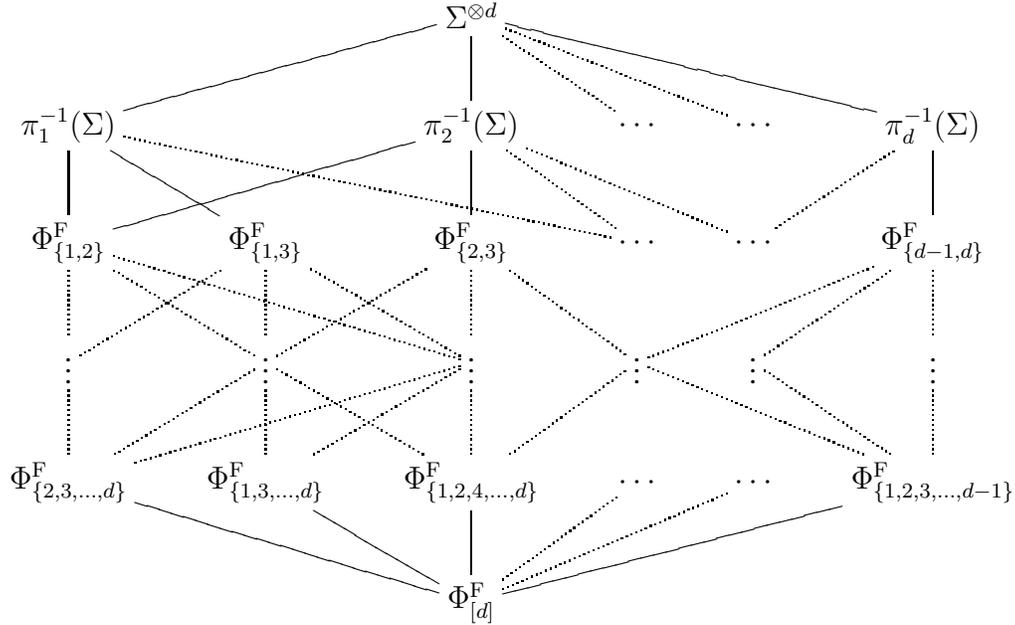
\begin{center}
$\phantom{i}$\xymatrix{ & & \S^{\otimes d}\ar@{-}[dll]\ar@{-}[d]\ar@{.}[dr]\ar@{.}[drr]\ar@{-}[drrr]\\
\pi_1^{-1}(\S)\ar@{-}[d]\ar@{-}[dr]\ar@{.}[drrr] & & \pi_2^{-1}(\S)\ar@{-}[dll]\ar@{-}[d]\ar@{.}[dr]\ar@{.}[drr] & \cdots & \cdots & \pi_d^{-1}(\S)\ar@{.}[dl]\ar@{-}[d]\\
\Phi^\rm{F}_{\{1,2\}}\ar@{.}[d]\ar@{.}[dr]\ar@{.}[drr] & \Phi^\rm{F}_{\{1,3\}}\ar@{.}[dr]\ar@{.}[d]\ar@{.}[dl] & \Phi^\rm{F}_{\{2,3\}}\ar@{.}[dr]\ar@{.}[d]\ar@{.}[dl] &\cdots &\cdots & \Phi^\rm{F}_{\{d-1,d\}}\ar@{.}[dll]\ar@{.}[dl]\ar@{.}[d]\\
\vdots\ar@{.}[d] & \vdots\ar@{.}[dl]\ar@{.}[d]\ar@{.}[dr]& \vdots\ar@{.}[dll]\ar@{.}[dl]\ar@{.}[d] & \vdots\ar@{.}[dl]\ar@{.}[drr] & \vdots\ar@{.}[dr]& \vdots\ar@{.}[d] \\
\Phi^\rm{F}_{\{2,3,\ldots,d\}}\ar@{-}[drr] & \Phi^\rm{F}_{\{1,3,\ldots,d\}}\ar@{-}[dr] & \Phi^\rm{F}_{\{1,2,4,\ldots,d\}}\ar@{-}[d] & \cdots\ar@{.}[dl] & \cdots\ar@{.}[dll] & \Phi^\rm{F}_{\{1,2,3,\ldots,d-1\}}\ar@{-}[dlll]\\
& & \Phi^\rm{F}_{[d]}
}
\end{center}
by showing that large collections of the $\s$-subalgebras appearing here are relatively independent over the collections of further $\s$-subalgebras that they have in common.

It is worth stressing at this point that we have \emph{not} proved any such assertion for the joint distribution of all the original factors $\Phi_e\leq \S$, but only for their oblique copies inside $\S^{\otimes d}$.  The problem of describing the joint distribution of the factors $\Phi_e$ themselves seems to be much harder, because it runs into precisely the difficulties with linear dependence discussed in Section~\ref{sec:meta}: for example, if $e_1,e_2,e_3 \subseteq [d]$ are three subsets that are pairwise non-disjoint, then we have $\Phi_{e_i} = \S^{T\uhr \G_{e_i}}$ for $\G_{e_i} = \sum_{j,j' \in e_i}\bbZ(\bf{e}_j - \bf{e}_{j'})$, and these three subgroups are now clearly not linearly independent.  In our analysis of the oblique factors $\Phi^\rm{F}_e$ we carefully avoided a similar problem during Step 1 of the proof of Proposition~\ref{prop:Fberg2}, where we exploited the fact that $\Phi^\rm{F}_e$ is contained modulo negligible sets in $\pi_j^{-1}(\S)$ for any choice of $j \in e$, so that by making careful choices of the coordinates with which to express these oblique copies we were able to reduce the joint distribution of interest to the case covered by Theorem~\ref{thm:sateds-joint-dist}, involving only linearly independent subgroups.  However, it seems clear that no similar trick will be available in the study of the factors $\Phi_e$.

Happily, however, we do not need any such more precise information to complete our proof of Theorem~\ref{thm:multirec2}: in the remainder of this chapter we show how the structure proved above for $\mu^\rm{F}$ suffices.  This will proceed through a slight modification of Tao's infinitary hypergraph removal lemma from~\cite{Tao07}, which first appeared in the form given below in~\cite{Aus--newmultiSzem}.

\begin{prop}\label{prop:infremoval}
Suppose that $(X,\S,\mu)$ is a standard Borel space and $\l$ is a $d$-fold coupling of $\mu$ on $(X^d,\S^{\otimes d})$ with coordinate projection maps $\pi_i:X^d\to X$, and that $(\Psi_e)_e$ is a collection of $\s$-subalgebras of $\S$ indexed by subsets $e \in \binom{[d]}{\geq 2}$ with the following properties:
\begin{itemize}
\item[{[i]}] if $e \subseteq e'$ then $\Psi_e \geq \Psi_{e'}$;
\item[{[ii]}] if $i,j \in e$ and $A \in \Psi_e$ then $\l(\pi_i^{-1}(A)\triangle \pi_j^{-1}(A)) = 0$, so that we may let $\Psi^\dag_e$ be the common $\l$-completion of the lifted $\s$-algebras $\pi_i^{-1}(\Psi_e)$ for $i\in e$;
\item[{[iii]}] if we define $\Psi^\dag_\I := \bigvee_{e\in \I}\Psi^\dag_e$ for each up-set $\I \in \binom{[d]}{\geq 2}$, then the $\s$-subalgebras $\Psi^\dag_\I$ and $\Psi^\dag_{\I'}$ are relatively independent under $\l$ over $\Psi^\dag_{\I\cap \I'}$.
\end{itemize}

In addition, suppose that $\I_{i,j}$ for $i=1,2,\ldots,d$ and $j =
1,2,\ldots,k_i$ are collections of up-sets in $\binom{[d]}{\geq 2}$
such that $[d] \in \I_{i,j} \subseteq \langle i\rangle$ for each
$i,j$, and that the sets $A_{i,j}\in
\Phi_{\I_{i,j}}$ are such that
\[\l\Big(\prod_{i=1}^d\Big(\bigcap_{j=1}^{k_i}A_{i,j}\Big)\Big) = 0.\]
Then we must also have
\[\mu\Big(\bigcap_{i=1}^d\bigcap_{j=1}^{k_i}A_{i,j}\Big) = 0.\]
\end{prop}

\noindent\textbf{Proof of Theorem~\ref{thm:multirec2} from Proposition~\ref{prop:infremoval}}\quad Clearly the conclusion holds for a system $\bfX$ if it holds for any extension of $\bfX$, so by Theorem~\ref{thm:sateds-exist} we may assume that $\bfX$ is $\sfC_e$-sated for every $e\subseteq [d]$.

Now suppose that $A_1,A_2,\ldots,A_d \in \S$ are
such that $\mu^{\rm{F}}(A_1\times A_2\times\cdots\times A_d) = 0$. Then by Proposition~\ref{prop:Fberg1} we have
\[\mu^{\rm{F}}(A_1\times A_2\times\cdots\times A_d) = \int_{X^d}\bigotimes_{i=1}^d \sfE_{\mu}(1_{A_i}\,|\,\Phi_{\langle i\rangle})\,\d\mu^{\rm{F}} = 0.\]
The level set $B_i := \{\sfE_{\mu}(1_A\,|\,\Phi_{\langle
i\rangle})
> 0\}$ (of course, this is unique only up to $\mu$-negligible sets) lies in $\Phi_{\langle i\rangle}$, and the above vanishing requires that also $\mu^{\rm{F}}(B_1\times B_2\times
\cdots\times B_d) = 0$.  Now setting $k_i=1$,
$\I_{i,1}:= \langle i\rangle$ and $A_{i,1} := B_i$ for each $i \leq
d$, Lemma~\ref{lem:diag} and Proposition~\ref{prop:Fberg2} imply that Proposition~\ref{prop:infremoval} applies to the partially invariant factors $\Phi_e$ and their oblique copies to give $\mu(B_1\cap
B_2\cap\cdots\cap B_d) = 0$. On the other hand we must have
$\mu(A\setminus B_i) = 0$ for each $i$, and so overall $\mu(A) \leq
\mu(B_1\cap B_2\cap\cdots\cap B_d) + \sum_{i=1}^d\mu(A\setminus B_i)
= 0$, as required. \qed

The remainder of this chapter is given to the proof of Proposition~\ref{prop:infremoval}. This proceeds by induction on a suitable ordering of the possible collections
of up-sets $(\I_{i,j})_{i,j}$, appealing to a handful of different
possible cases at different steps of the induction.  At the
outermost level, this induction will be organized according to the
depth of our up-sets.

The proof given below is taken essentially unchanged from~\cite{Aus--newmultiSzem}, where in turn the statement and proof were adopted with only slight modifications from~\cite{Tao07}.  The reader may consult~\cite{Aus--newmultiSzem} for an explanation of these modifications.

\begin{dfn}
A family $(\I_{i,j})_{i,j}$ has the property \textbf{P} if it
satisfies the conclusion of Proposition~\ref{prop:infremoval}.
\end{dfn}

We separate the various components of the induction into separate lemmas.

\begin{lem}[Lifting using relative independence]\label{lem:lift-rel-ind}
Suppose that all up-sets in the collection $(\I_{i,j})_{i,j}$ have
depth at least $k$, that all those with depth exactly $k$ are
principal, and that there are $\ell \geq 1$ of these.  Then if
property P holds for all similar collections having $\ell - 1$
up-sets of depth $k$, then it holds also for this collection.
\end{lem}

\noindent\textbf{Proof}\quad Let $\I_{i_1,j_1} = \langle e_1\rangle$,
$\I_{i_2,j_2} = \langle e_2\rangle$, \ldots, $\I_{i_\ell,j_\ell} =
\langle e_\ell\rangle$ be an enumeration of all the (principal)
up-sets of depth $k$ in our collection.  We will treat two separate
cases.

First suppose that two of the generating sets agree; by re-ordering
if necessary we may assume that $e_1 = e_2$.  Clearly we can assume
that there are no duplicates among the coordinate-collections
$(\I_{i,j})_{j=1}^{k_i}$ for each $i$ separately, so we must have
$i_1 \neq i_2$. However, if we now suppose that $A_{i,j}\in
\I_{i,j}$ for each $i$, $j$ are such that
\[\l\Big(\prod_{i=1}^d\Big(\bigcap_{j=1}^{k_i}A_{i,j}\Big)\Big) = 0,\]
then by assumption [ii] the same equality holds if we simply replace $A_{i_1,j_1} \in
\langle e_1\rangle$ with $A'_{i_1,j_1}:= A_{i_1,j_1}\cap
A_{i_2,j_2}$ and $A_{i_2,j_2}$ with $A'_{i_2,j_2} := X$. Now this
last set can simply be ignored to leave an instance of a
$\l$-negligible product for the same collection of
up-sets omitting $\I_{i_2,j_2}$, and so property P of this reduced
collection completes the proof.

On the other hand, if all the $e_i$ are distinct, we shall simplify
the last of the principal up-sets $\I_{i_\ell,j_\ell}$ by exploiting
the relative independence among the lifted $\s$-algebras $\Psi_e^\dag$. Assume for notational simplicity that $(i_\ell,j_\ell) =
(1,1)$; clearly this will not affect the proof.  We will reduce to
an instance of property P associated to the collection $(\I'_{i,j})$
defined by
\[\I'_{i,j} := \left\{\begin{array}{ll}\langle e_\ell\rangle\setminus\{e_\ell\}&\quad\hbox{if}\ (i,j) = (1,1)\\ \I_{i,j}&\quad\hbox{else,}\end{array}\right.\]
which has one fewer up-set of depth $k$ and so falls under the
inductive assumption.

Indeed, by property [iii] under
$\l$ the set $\pi_1^{-1}(A_{1,1})$ is relatively
independent from all the sets $\pi_i^{-1}(A_{i,j})$, $(i,j) \neq
(1,1)$, over the $\s$-algebra $\pi_1^{-1}(\Psi_{\langle
e_\ell\rangle\setminus\{e_\ell\}})$, which is dense inside $\Psi^\dag_{\langle
e_\ell\rangle\setminus\{e_\ell\}}$.  Therefore
\begin{multline*}
0 =
\l\Big(\prod_{i=1}^d\Big(\bigcap_{j=1}^{k_i}A_{i,j}\Big)\Big)\\
= \int_{X^d}\sfE_\mu(1_{A_{1,1}}\,|\,\Psi_{\langle
e_\ell\rangle\setminus\{e_\ell\}})\circ\pi_1\cdot\prod_{j=2}^{k_1}1_{\pi_1^{-1}(A_{1,j})}\cdot\prod_{i=2}^d\prod_{j=1}^{k_i}1_{\pi_i^{-1}(A_{i,j})}\,\d\l.
\end{multline*}
Setting $A'_{1,1}:= \{\sfE_\mu(1_{A_{1,1}}\,|\,\Psi_{\langle
e_\ell\rangle\setminus\{e_\ell\}}) > 0\} \in \Psi_{\langle
e_\ell\rangle\setminus\{e_\ell\}}$ and $A'_{i,j} := A_{i,j}$ for
$(i,j) \neq (1,1)$, we have that $\mu(A_{1,1}\setminus A'_{1,1}) =
0$ and it follows from the above equality that also
$\l\big(\prod_{i=1}^d\big(\bigcap_{j=1}^{k_i}A'_{i,j}\big)\big)
= 0$, so an appeal to property P for the reduced collection of
up-sets completes the proof. \qed

\begin{lem}[Lifting under finitary generation]\label{lem:lift-gen}
Suppose that all up-sets in the collection $(\I_{i,j})_{i,j}$ have
depth at least $k$ and that among those of depth $k$ there are $\ell
\geq 1$ that are non-principal.  Then if property P holds for all
similar collections having at most $\ell - 1$ non-principal up-sets
of depth $k$, then it also holds for this collection.
\end{lem}

\noindent\textbf{Proof}\quad Let $\I_{i_1,j_1}$, $\I_{i_2,j_2}$, \ldots,
$\I_{i_\ell,j_\ell}$ be the non-principal up-sets of depth $k$, and
now in addition let $e_1$, $e_2$, \ldots, $e_r$ be all the members
of $\I_{i_\ell,j_\ell}$ of size $k$ (so, of course, $r \leq
\binom{d}{k}$). Once again we will assume for simplicity that
$(i_\ell,j_\ell) = (1,1)$.  We break our work into two further
steps.

\textbf{Step 1}\quad First consider the case of a collection
$(A_{i,j})_{i,j}$ such that for the set $A_{1,1}$, we can actually
find \emph{finite} subalgebras of sets $\B_s\in \Psi_{\{e_s\}}$ for
$s = 1,2,\ldots,r$ such that $A_{i_\ell,j_\ell} \in \B_1\vee
\B_2\vee \cdots \vee \B_r\vee\Psi_{\I_{1,1}\cap\binom{[d]}{\geq
k+1}}$ (so $A_{1,1}$ lies in one of our non-principal up-sets of
depth $k$, but it fails to lie in an up-set of depth $k+1$ only `up
to' finitely many additional generating sets). Choose $M\geq
\max_{s\leq r}|\B_s|$, so that we can certainly express
\[A_{1,1} = \bigcup_{m=1}^{M^r} (B_{m,1}\cap B_{m,2}\cap\cdots\cap B_{m,r}\cap C_m)\]
with $B_{m,s}\in \B_s$ for each $s\leq r$ and $C_m \in
\Psi_{\I_{1,1}\cap\binom{[d]}{\geq k+1}}$.  Inserting this
expression into the equation
\[\l\Big(\prod_{i=1}^d\Big(\bigcap_{j=1}^{k_i}A_{i,j}\Big)\Big)
= 0\] now gives that each of the $M^r$ individual product sets
\[\Big((B_{m,1}\cap B_{m,2}\cap\cdots\cap B_{m,r}\cap C_m)\cap \bigcap_{j=2}^{k_1}A_{1,j}\Big)\times\prod_{i=2}^d \Big(\bigcap_{j=1}^{k_i}A_{i,j}\Big)\]
is $\l$-negligible.

Now consider the family of up-sets comprising the original
$\I_{i,j}$ if $i=2,3,\ldots,d$ and the collection $\langle
e_1\rangle$, $\langle e_2\rangle$, \ldots, $\langle e_r\rangle$,
$\I_{1,2}$, $\I_{1,3}$, \ldots, $\I_{1,k_1}$ corresponding to $i =
1$. We have broken the depth-$k$ non-principal up-set $\I_{1,1}$
into the higher-depth up-set $\I_{1,1}\cap\binom{[d]}{\geq k+1}$ and
the principal up-sets $\langle e_s\rangle$, and so there are only
$\ell-1$ minimal-depth non-principal up-sets in this new family.  It
is clear that for each $m\leq M^r$ the above product set is
associated to this family of up-sets, and so an inductive appeal to
property P for this family tells us that also
\[\mu\Big((B_{m,1}\cap B_{m,2}\cap\cdots\cap B_{m,r}\cap C_m)\cap \bigcap_{j=2}^{k_1}A_{1,j}\cap \bigcap_{i=2}^d \bigcap_{j=1}^{k_i}A_{i,j}\Big) = 0\]
for every $m\leq M^r$.  Since the union of these sets is just
$\bigcap_{i=1}^d\bigcap_{j=1}^{k_i} A_{i,j}$, this gives the desired
negligibility in this case.

\textbf{Step 2}\quad Now we return to the general case, which will
follow by a suitable limiting argument applied to the conclusion of
Step 1.  Since any $\Psi_e$ is countably generated modulo $\mu$, for each
$e$ with $|e| = k$ we can find an increasing sequence of finite
subalgebras $\B_{e,1} \subseteq \B_{e,2} \subseteq\ldots$ that
generates $\Psi_e$ up to $\mu$-negligible sets. In terms of
these define approximating sub-$\s$-algebras
\[\Xi^{(n)}_{i,j} := \Psi_{\I_{i,j}\cap \binom{[d]}{\geq k+1}}\vee \bigvee_{e \in \I_{i,j}\cap \binom{[d]}{k}} \B_{e,n},\]
so for each $\I_{i,j}$ these form an increasing family of
$\s$-algebras that generates $\Psi_{\I_{i,j}}$ up to
$\mu$-negligible sets (inded, if $\I_{i,j}$ does not contain any
sets of the minimal depth $k$ then we simply have $\Xi^{(n)}_{i,j} = \Psi_{\I_{i,j}}$ for all $n$). Now property [iii] implies for each $n$ that $\Psi^\dag_{\I_{1,1}}$ and $\bigvee_{(i,j)\neq
(1,1)}\pi_i^{-1}(\Xi^{(n)}_{i,j})$ are relatively independent over
$\pi_1^{-1}(\Xi^{(n)}_{1,1})$.

Given now a family of sets $(A_{i,j})_{i,j}$ associated to
$(\I_{i,j})_{i,j}$, for each $(i,j)$ the conditional expectations
$\sfE_\mu(1_{A_{i,j}}\,|\,\Xi^{(n)}_{i,j})$ form an almost surely
uniformly bounded martingale converging to $1_{A_{i,j}}$ in
$L^2(\mu)$. Letting \[B^{(n)}_{i,j}:=
\{\sfE_\mu(1_{A_{i,j}}\,|\,\Xi^{(n)}_{i,j})
> 1-\delta\}\] for some small $\delta > 0$ (to be specified
momentarily), it is clear that we also have $\mu(A_{i,j}\triangle
B_{i,j}^{(n)}) \to 0$ as $n\to\infty$.  Let
\[F := \prod_{i=1}^d\Big(\bigcap_{j=1}^{k_i}B^{(n)}_{i,j}\Big).\] We
now compute using the above-mentioned relative independence that
\begin{eqnarray*}
&&\l(F\setminus \pi_i^{-1}(A_{i,j}))\\ && \quad = \int_{X^d}
\Big(\prod_{(i',j')}1_{B^{(n)}_{i',j'}}\circ\pi_{i'}\Big) - 1_{A_{i,j}}\circ\pi_i\cdot\Big(\prod_{(i',j')}1_{B^{(n)}_{i',j'}}\circ\pi_{i'}\Big)\,\d\l\\
&&\quad = \int_{X^d} (1_{B^{(n)}_{i,j}\setminus A_{i,j}}\circ\pi_i)\cdot
\Big(\prod_{(i',j')\neq
(i,j)}1_{B^{(n)}_{i',j'}}\circ\pi_{i'}\Big)\,\d\l\\
&&\quad = \int_{X^d} (\sfE_\mu(1_{B^{(n)}_{i,j}\setminus
A_{i,j}}\,|\,\Xi^{(n)}_{i,j})\circ\pi_i)\cdot
\Big(\prod_{(i',j')\neq
(i,j)}1_{B^{(n)}_{i',j'}}\circ\pi_{i'}\Big)\,\d\l
\end{eqnarray*}
for each pair $(i,j)$.

However, from the definition of $B^{(n)}_{i,j}$ we must have
\[\sfE_\mu(1_{B^{(n)}_{i,j}\setminus A_{i,j}}\,|\,\Xi^{(n)}_{i,j}) \leq \delta 1_{B^{(n)}_{i,j}}\]
almost surely, and therefore the above integral inequality implies
that \[\l(F\setminus \pi_i^{-1}(A_{i,j}))\leq
\delta\int_{X^d} (1_{B^{(n)}_{i,j}}\circ\pi_i)\cdot
\Big(\prod_{(i',j')\neq
(i,j)}1_{B^{(n)}_{i',j'}}\circ\pi_{i'}\Big)\,\d\l =
\delta \l(F).\]

From this we can estimate as follows:
\[\l(F) \leq
\l\Big(\prod_{i=1}^d\Big(\bigcap_{j=1}^{k_i}A_{i,j}\Big)\Big)
+ \sum_{(i,j)}\l(F\setminus \pi_i^{-1}(A_{i,j})) \leq 0 + \Big(\sum_{i=1}^dk_i\Big)\delta \l(F),\]
and so provided we chose $\delta < \big(\sum_{i=1}^dk_i\big)^{-1}$
we must in fact have $\l(F) = 0$.

We have now obtained sets $(B^{(n)}_{i,j})_{i,j}$ that are
associated to the family $(\I_{i,j})_{i,j}$ and satisfy the property
of lying in finitely-generated extensions of the relevant factors
corresponding to the members of the $\I_{i,j}$ of minimal size, and
so we can apply the result of Step 1 to deduce that
$\mu\big(\bigcap_{i=1}^d\bigcap_{j=1}^{k_i}B^{(n)}_{i,j}\big) = 0$.
It follows that
\[\mu\Big(\bigcap_{i=1}^d\bigcap_{j=1}^{k_i}A_{i,j}\Big) \leq \sum_{i,j}\mu(A_{i,j}\setminus B^{(n)}_{i,j}) \to 0\quad\quad\hbox{as }n\to\infty,\]
as required. \qed

\noindent\textbf{Proof of Proposition~\ref{prop:infremoval}}\quad We first
take as our base case $k_i = 1$ and $\I_{i,1} = \{[d]\}$ for each
$i=1,2,\ldots,d$.  In this case we know from property [ii] that for any $A \in
\Psi_{[d]}$ the pre-images $\pi_i^{-1}(A)$ are all equal up to
negligible sets, and so given $A_1$, $A_2$, \ldots, $A_d \in
\Psi_{[d]}$ we have $0 = \l(A_1\times
A_2\times\cdots\times A_d) = \mu(A_1\cap A_2\cap\cdots\cap A_d)$.

The remainder of the proof now just requires putting the preceding
lemmas into order to form an induction with three layers: if our
collection has any non-principal up-sets of minimal depth, then
Lemma~\ref{lem:lift-gen} allows us to reduce their number at the
expense only of introducing new principal up-sets of the same depth;
and having removed all the non-principal minimal-depth up-sets,
Lemma~\ref{lem:lift-rel-ind} enables us to remove also the principal
ones until we are left only with up-sets of increased minimal depth.
This completes the proof. \qed

\chapter{The Density Hales-Jewett Theorem}\label{chap:DHJ}

Much as for Szemer\'edi's Theorem and its multidimensional
generalization, the Ergodic Ramsey Theory approach to Theorem B
begins by establishing its equivalence to a result about stochastic
processes.  We have deferred the introduction of the stochastic
processes analog of Theorem B until now because it involves a less
well-known family of processes than the tuples of commuting
transformations that appear in Theorem A, and these new stochastic
processes require a separate discussion.  The proof
from~\cite{FurKat91} of the correspondence between Theorem B and an
assertion about these processes is also less well-known, and so we
recall this in the first section below for completeness.

After formulating the stochastic processes result to which Theorem B is equivalent, we introduce an additional semigroup $\G$ of transformations on these processes and argue that we may reduce further to the case of processes whose distributions are invariant.  This leaves us with a class of $\G$-systems, on which we will bring a notion of satedness to bear.  However, as promised at the beginning of Chapter 2, this first requires some modifications to that notion, effectively by imposing additional restrictions on the factor maps we allow in our theory of a kind not involved heretofore.  With these modifications in place we will proceed to analogs of Propositions~\ref{prop:Fberg1} and~\ref{prop:Fberg2} and thence to the proof of Theorem B.

\section{The correspondence with a class of stationary processes}\label{sec:startDHJ}

\subsection*{Combinatorial notation}

In addition to the finite spaces $[k]^N$ appearing in the statement of Theorem B, we will work with their union
\[[k]^\ast := \bigcup_{N\geq 1}[k]^N.\]
The spaces $[k]^N$ and $[k]^\ast$ are referred to as the \textbf{$N$-dimensional} and \textbf{infinite-dimensional combinatorial spaces} over the alphabet $[k]$ respectively. Most of this chapter will consider probabilities on product spaces
indexed by $[k]^\ast$.

If $A\subseteq [k]^N$ then we denote its \textbf{density} by
\[\d(A) := \frac{|A|}{k^N};\]
thus the assumption of Theorem B is that $N$ is sufficiently large in terms of $k$ and $\d(A)$.

Given two finite words $u,v \in [k]^\ast$ we denote their
concatenation by either $uv$ or $u\oplus v$.  For any finite $n$ we
define an \textbf{$n$-dimensional subspace} of $[k]^\ast$ to be an
injection $\phi:[k]^n \into [k]^\ast$ specified as follows: for some
integers $0 = N_0 < N_1 < N_2 < \ldots < N_n$, nonempty subsets $I_1
\subseteq [N_1]$, $I_2 \subseteq [N_2]\setminus [N_1]$, \ldots, $I_n
\subseteq [N_n]\setminus [N_{n-1}]$ and a word $w \in [k]^{N_n}$ we
let $\phi(v_1v_2\cdots v_n)$ be the word in $[k]^\ast$ of length
$N_n$ given by
\[\phi(v_1v_2\cdots v_n)_m := \left\{\begin{array}{ll}w_m&\quad\quad\hbox{if }m\in [N_n]\setminus (I_1 \cup I_2\cup\cdots\cup I_n)\\ v_i&\quad\quad\hbox{if }m\in I_i.\end{array}\right.\]

In these terms a combinatorial line is simply a $1$-dimensional
combinatorial subspace.

Similarly, an \textbf{infinite-dimensional subspace} (or often just
\textbf{subspace}) of $[k]^\ast$ is an injection $\phi:[k]^\ast
\into [k]^\ast$ specified using some infinite sequence $0 = N_0 <
N_1 < N_2 < \ldots$, nonempty subsets $I_{i+1} \subseteq
[N_{i+1}]\setminus [N_i]$ and words $w_i \in [k]^{N_i}$, where for
any $v \in [k]^n$ its image $\phi(v)$ has length $N_n$ and is given
by the above formula with $w := w_n$. It is clear that the
collection of all subspaces of $[k]^\ast$ forms a semigroup $\G$
under composition.

Finally, let us define \textbf{letter-replacement maps}: give $i\in
[k]$ and $e \subseteq [k]$, for each $N\geq 1$ we define
$r^N_{e,i}:[k]^N \to [k]^N$ by
\[r^N_{e,i}(w)_m := \left\{\begin{array}{ll}i&\quad\quad\hbox{if }w_m \in e\\ w_m&\quad\quad\hbox{if }w_m \in [k]\setminus e\end{array}\right.\]
for $m\leq N$, and let \[r_{e,i} := \bigcup_{N\geq
1}r^N_{e,i}:[k]^\ast\to [k]^\ast\] (so clearly $r_{e,i}$ actually takes
values in the subset $(([k]\setminus e)\cup \{i\})^\ast\subseteq [k]^\ast$).

\subsection*{Reformulation in terms of stochastic processes}

The correspondence that Furstenberg and Katznelson establish for Theorem B is between dense subsets of the finite-dimensional combinatorial spaces $[k]^N$ and stochastic processes \emph{indexed} by the infinite-dimensional combinatorial space $[k]^\ast$.

\begin{thm}[Infinitary Density Hales-Jewett Theorem]\label{thm:inf-DHJ}
For any $\delta > 0$, if $\mu$ is a Borel probability measure on
$\{0,1\}^{[k]^\ast}$ with the property that
\[\mu\{\bf{x} \in \{0,1\}^{[k]^\ast}:\ x_w = 1\} \geq \delta\quad\quad\forall w \in [k]^\ast,\]
then there is a combinatorial line $\phi:[k]\into [k]^\ast$ such that
\[\mu\{\bf{x} \in \{0,1\}^{[k]^\ast}:\ x_{\phi(i)} = 1\ \forall i \in [k]\} > 0.\]
\end{thm}

\noindent\textbf{Proof of Theorem B from Theorem~\ref{thm:inf-DHJ}}\quad Clearly we may restrict our attention to $k\geq 2$. We will suppose that theorem B fails, and show that this would give rise to a counterexample to Theorem~\ref{thm:inf-DHJ}.
We break this into two steps.

\textbf{Step 1}\quad First observe that if $N\geq L\geq 1$ and $A\subseteq [k]^N$ has $\d(A) > 1 - \frac{1}{k^{2L}}$ then $A$ necessarily contains a whole $L$-dimensional combinatorial subspace.  Indeed, having density as high as this implies that each of the $k^L$ subsets
\[A_u := \{w \in [k]^{N-L}:\ u\oplus w \in A\}\quad\quad\hbox{for}\ u \in [k]^L\]
has density greater than $1 - \frac{1}{k^L}$, and so there must be some $w \in \bigcap_{u\in [k]^L}A_u$, implying that the subspace $[k]^L\into [k]^N:u\mapsto u\oplus w$ has image lying entirely in $A$.

In particular, letting $L = 1$, if we assume that Theorem B fails then we may let
\[\delta_0 := \sup\{\delta > 0:\ \hbox{Theorem B fails for subsets of density $\delta$}\}\]
and deduce that $0 < \delta_0 < 1$.

\textbf{Step 2}\quad Now fix some integer $L \geq 1$ and let $A \subseteq [k]^{N}$ be a subset of density $\d(A) = \delta > (1 + \frac{1}{2k^{L+1}})^{-1}\delta_0$ for some $N \geq L$ such that $A$ contains no combinatorial lines.

Let $N = L + M$ and decompose $[k]^{N}$ as $[k]^L \oplus [k]^{M}$. For each $w \in [k]^L$ let
\[A_w := \{v \in [k]^{M}:\ w\oplus v \in A\}.\]
Clearly
\[\frac{1}{k^L}\sum_{w\in [k]^L}\d(A_w) = \d(A) = \delta,\]
and on the other hand $\d(A_w) < (1 + \frac{1}{2k^{L+1}})\delta$ for
each $w$ once $N$ is sufficiently large, for otherwise $A_w$ would
contain a combinatorial line by the definition of $\delta_0$.
Therefore the above equation between densities and Chebyshev's
inequality require that in fact \emph{every} $w \in [k]^L$ have
$\d(A_w) > \delta/2$.

Now defining the probability measure $\mu_L$ on $\{0,1\}^{[k]^L}$ by
\[\mu_L\{(x_w)_{w\in [k]^L}\} := \d(\{v \in [k]^M:\ x_w = 1_{A_w}(v)\ \forall w\in [k]^L\})\]
for each $(x_w)_{w\in [k]^L} \in \{0,1\}^{[k]^L}$, we see that for each $L$ we have produced a probability $\mu_L$ on $\{0,1\}^{[k]^L}$ such that
\[\mu_L\{\bf{x} \in \{0,1\}^{[k]^L}:\ x_w = 1\} = \d(A_w) \geq \delta/2 \geq \delta_0/4\quad\quad\forall w \in [k]^L\]
but
\[\mu_L\{\bf{x} \in \{0,1\}^{[k]^L}:\ x_{\phi(i)} = 1\ \forall i \in [k]\} = 0\]
for any combinatorial line $\phi:[k]\into [k]^L$.  Finally defining $\mu := \bigotimes_{L\geq 1}\mu_L$, we obtain a measure that contradicts Theorem~\ref{thm:inf-DHJ} with density $\delta_0/4$. \qed

\noindent\textbf{Remark}\quad The above proof is essentially taken from Proposition 2.1 of~\cite{FurKat91}, where the reverse implication is also proved. \fin

\section{Strongly stationary processes}

After introducing Theorem~\ref{thm:inf-DHJ}, Furstenberg and Katznelson make a further reduction to a special subclass of measures.

\begin{dfn}[Semigroup action of combinatorial subspaces]
If $\phi:[k]^N\into [k]^\ast$ is a combinatorial subspace then for any product space $K^{[k]^\ast}$ we define the corresponding map $T_\phi:K^{[k]^\ast}\to K^{[k]^N}$ by
\[(T_\phi(\bf{x}))_w := x_{\phi(w)}\quad\quad\hbox{for}\ w \in [k]^N\ \hbox{and}\ \bf{x} = (x_u)_{u\in [k]^\ast} \in K^{[k]^\ast},\]
and similarly define $T_\phi:K^{[k]^\ast}\to K^{[k]^\ast}$ in case $\phi:[k]^\ast\into [k]^\ast$.  In the latter case this specifies an action $\G\actson K^{[k]^\ast}$.
\end{dfn}

\begin{dfn}[Strongly stationary laws]
A probability measure $\mu$ on the product $(K^{[k]^\ast},\Psi^{\otimes [k]^\ast})$ for some standard Borel space $(K,\Psi)$ is \textbf{strongly stationary} if $T_{\phi\#}\mu = \mu$ for all subspaces $\phi \in \G$.  In this case the transformations $T_\phi$ give to $(K^{[k]^\ast},\Psi^{\otimes [k]^\ast},\mu)$ the structure of a probability-preserving $\G$-system.
\end{dfn}

\begin{lem}\label{lem:sssuffice}
If Theorem~\ref{thm:inf-DHJ} holds for all strongly stationary measures for any $\delta > 0$ then it holds for all measures satisfying the conditions of that theorem for any $\delta > 0$.
\end{lem}

\noindent\textbf{Proof}\quad This argument is again lifted directly from~\cite{FurKat91}, and we only sketch the details.  Given a measure $\mu$ satisfying the conditions of Theorem~\ref{thm:inf-DHJ} for some $\delta > 0$, by applying the Carlson-Simpson Theorem~\cite{Car88} to arbitrarily
fine finite open coverings of the finite-dimensional spaces of
probability distributions on $\{0,1\}^{[k]^n}$ for increasingly
large $n$, we obtain a subspace $\psi:[k]^\ast \into [k]^\ast$ and an
infinite word $w = w_1w_2\cdots \in [k]^\bbN$ such that the restricted laws
\[T_{\psi(w_1w_2\cdots w_m\oplus \ \cdot\ )\#}\mu\]
converge to a strongly stationary law as $m\to\infty$, and since
\emph{all} one-dimensional marginals of the input law gave
probability at least $\delta$ to $\{1\}$, the same is true of the
limit measure.  Finally, the subset of probability measures
\[\big\{\nu \in \Pr\{0,1\}^{[k]^\ast}:\ \nu\{\bf{x} \in \{0,1\}^{[k]^\ast}:\ x_{\phi(i)} = 1\ \forall i\leq k\} > 0\big\}\]
is finite-dimensional and open for any given line $\phi:[k]\into
[k]^\ast$, so if the limit measure is in this set the so is some
image of the original measure. \qed

An immediate consequence of the strong stationarity of a measure $\mu$ is that for any two $N$-dimensional subspaces $\phi,\psi:[k]^N\into [k]^\ast$ we have $T_{\phi\#}\mu = T_{\psi\#}\mu$.  In case $N = 0$ we refer to this common image measure as the \textbf{point marginal} $\mu$ and denote it by $\mu^{\rm{pt}}$, and similarly in case $N=1$ it is the \textbf{line marginal} of $\mu$ and is denoted by $\mu^{\rm{line}}$.  In these terms it is possible to give another, more convenient reformulation of Theorem~\ref{thm:inf-DHJ}.

\begin{thm}\label{thm:inf-DHJ2}
If $(K,\Psi)$ is a standard Borel space and $\mu$ is a strongly stationary law on $(K^{[k]^\ast},\Psi^{\otimes [k]^\ast})$ then for any $A_1,A_2,\ldots,A_k \in \Psi$ we have
\[\mu^{\rm{line}}(A_1\times A_2\times \cdots \times A_k) = 0\quad\quad\Rightarrow\quad\quad \mu^{\rm{pt}}(A_1\cap A_2\cap\cdots \cap A_k) = 0.\]
\end{thm}

The resemblance to Theorem~\ref{thm:multirec2} is far from accidental!

The proof of Theorem~\ref{thm:inf-DHJ2} will involve a version of satedness for our systems of interest; however, here a slight subtlety creeps in.  In the following we will need to work with only those $\G$-systems that are of the form $(K^{[k]^\ast},\Psi^{[k]^\ast},\mu,T)$ for some strongly stationary measure $\mu$ (of course, the huge semigroup $\G$ could also have invariant measures for all sorts of other Borel actions, not of this form).  On the other hand, the conclusion of Theorem~\ref{thm:inf-DHJ2} is not about the joint distribution of several copies of whole $\G$-systems under some self-joining.  Rather, it is about the joint distribution of some copies of just the `one-dimensional' point marginal $(K,\Psi,\mu^{\rm{pt}})$ under the line marginal: this is only a tiny fragment of the whole system $(K^{[k]^\ast},\Psi^{\otimes [k]^\ast},\mu,T)$.

The way we can keep track of the structure of point and line marginals between different such systems is by restricting the kinds of factor map we allow.

\begin{dfn}
Let $\sfA$ be the class of $\G$-systems given by strongly stationary measures on product spaces indexed by $[k]^\ast$, as above.

A \textbf{coordinatewise factor} (or \textbf{cw-factor}) of $\bfX = (K^{[k]^\ast},\Psi^{\otimes [k]^\ast},\mu,T) \in \sfA$ is a $\s$-subalgebra of the form $\Phi^{\otimes [k]^\ast} \leq \Psi^{\otimes [k]^\ast}$ for some $\Phi \leq \Psi$.  Slightly abusively, we will sometimes refer instead to the single-coordinate $\s$-subalgebra $\Psi$ as a cw-factor. Likewise, a \textbf{cw-factor map} is a map of the form
\[f^\ast:(K^{[k]^\ast},\Psi^{\otimes [k]^\ast},\mu,T)\to (L^{[k]^\ast},\Xi^{\otimes [k]^\ast},\nu,T):(x_w)_w \mapsto (f(x_w))_w\]
for some Borel map $f:(K,\Psi)\to (L,\Xi)$, and $f^\ast$ is a \textbf{cw-isomorphism} if $f$ is measurably invertible away from some $\mu^{\rm{pt}}$- and $\nu^{\rm{pt}}$-negligible sets (this is clearly equivalent to its being an isomorphism in the usual sense).

With $f^\ast$ as above we shall sometimes refer to $f$ as its corresponding \textbf{single-coordinate map}.
\end{dfn}

It is now easy to see that the class $\sfA$ is closed under joinings
and inverse limits, provided that we interpret a joining of two
systems $(K^{[k]^\ast},\Psi^{\otimes [k]^\ast},\mu,T)$ and
$(L^{[k]^\ast},\Xi^{\otimes [k]^\ast},\nu,T)$ as a strongly
stationary measure on $(K\times L)^{[k]^\ast}$ and that we restrict
our attention to inverse sequences whose connecting maps are all
cw-factor maps.  We will henceforth refer to a subclass
$\sfC\subseteq \sfA$ as \textbf{cw-idempotent} if it is closed under
cw-isomorphisms, joinings and inverse limits involving cw-factor
maps, and now observe that all of the definitions and lemmas of
Section~\ref{sec:idem} have direct analogs for cw-idempotent classes
obtained simply by insisting that all morphisms be given by
cw-factor maps. In particular, if $\sfC$ is a cw-idempotent class
and $\bfX = (K^{[k]^\ast},\Psi^{\otimes [k]^\ast},\mu,T)\in \sfA$
then the maximal cw-$\sfC$-factor of $\bfX$ is given by
$\Phi^{\otimes [k]^\ast}$ where $\Phi$ is the maximal $\s$-algebra
in the family
\begin{multline*}
\big\{\Xi\leq \Psi:\ \Xi\ \hbox{is generated by some Borel map}\ f:(K,\Psi)\to (K_1,\Psi_1)\\
 \hbox{such that}\ (K_1^{[k]^\ast},\Psi_1^{\otimes [k]^\ast},f^\ast_\#\mu,T) \in \sfC\big\}.
\end{multline*}
We will write a cw-factor map coordinatizing this maximal $\sfC$-factor as $\zeta^\ast_\sfC$ for some map $\zeta_\sfC:K\to \sfC K$ of single-coordinate spaces.

Given these observations we can make our analog of Definition~\ref{dfn:sated}.

\begin{dfn}[CW-sated systems]\label{dfn:Csated}
For a cw-idempotent class $\sfC\subseteq \sfA$, a system $\bfX \in \sfA$ is \textbf{cw-$\sfC$-sated} if for any cw-extension
\[\pi^\ast:\t{\bfX} = (\t{K}^{[k]^\ast},\t{\Psi}^{\otimes[k]^\ast},\t{\mu},T)\to \bfX\]
the single-coordinate maps $\pi:\t{K}\to K$ and $\t{\zeta}_\sfC:\t{K}\to \sfC\t{K}$ are relatively independent under $\t{\mu}^{\rm{pt}}$ over $\zeta_\sfC\circ\pi:\t{K}\to K\to\sfC K$, where $\t{\zeta}^\ast_\sfC$ and $\zeta_\sfC^\ast$ coordinatize the maximal $\sfC$-factors of $\t{\bfX}$ and $\bfX$ respectively.
\end{dfn}

\begin{thm}\label{thm:Csateds-exist}
If $(\sfC_i)_{i\in I}$ is a
countable family of cw-idempotent classes then any system $\bfX_0 \in \sfA$
admits a cw-extension $\pi:\bfX \to \bfX_0$ that is cw-$\sfC_i$-sated for every $i\in I$.
\end{thm}

\noindent\textbf{Proof outline}\quad This proceeds in exact analogy with the proof of Theorem~\ref{thm:sateds-exist}.  First, applying the argument for Lemma~\ref{lem:inverse-lim-sated} to a bounded measurable function $f$ on the single-coordinate space of an inverse limit shows that an inverse limit of cw-$\sfC$-sated systems through cw-factor maps is cw-$\sfC$-sated.

Next, given a system
\[\bfX = (K^{[k]^\ast},\Psi^{\otimes [k]^\ast},\mu,T)\in \sfA,\]
we show how to produce a cw-sated extension for a single cw-idempotent class $\sfC$: first enumerate an $L^2$-dense sequence $(f_r)_{r\geq 1}$ in the unit ball of $L^\infty(\mu^{\rm{pt}})$; then apply the same exhaustion argument as in Step 1 of Theorem~\ref{thm:sateds-exist} to produce an inverse sequence of cw-extensions
\begin{multline*}
\ldots\stackrel{(\zeta^{n+2}_{n+1})^\ast}{\to} \bfX_{n+1} = (K_{n+1}^{[k]^\ast},\Psi_{n+1}^{\otimes [n]^\ast},\mu_{n+1},T)\\ \stackrel{(\zeta^{n+1}_n)^\ast}{\to}\bfX_n = (K_n^{[k]^\ast},\Psi_n^{\otimes [n]^\ast},\mu_n,T)\stackrel{(\zeta^n_{n-1})^\ast}{\to}\ldots\to \bfX
\end{multline*}
such that for each $r$ it happens cofinally often that this extension is within a factor of $2$ of achieving the optimal increase in the $L^2$-norm of the conditional expectation $\sfE_{\mu^{\rm{pt}}_n}(f_r\circ\psi^n_0\,|\,\zeta_\sfC^{(n)})$ (where $\zeta_\sfC^{(n)}$ is the single-coordinate map coordinatizing $\sfC\bfX_n$); and finally take the inverse limit of this sequence.  Just as in the proof of Theorem~\ref{thm:sateds-exist}, if this inverse limit were not cw-$\sfC$-sated then this would lead to a contradiction with our assumption on the increase of $\|\sfE_{\mu^{\rm{pt}}_n}(f_r\circ\psi^n_0\,|\,\zeta_\sfC^{(n)})\|_2$ for some finite $n$.

Finally the proof is completed by arguing that given a countable collection of cw-idempotent classes $\sfC_i$, we can produce one long inverse sequence of extensions in which for each $i$ there is a cofinal subsequence of cw-$\sfC_i$-sated systems, so that the inverse limit is cw-$\sfC_i$-sated for every $i$. \qed

This completes the modifications we need for our approach to Theorem B.  Note that detailed proofs of the above results written in the setting of strongly stationary laws are given in~\cite{Aus--DHJ}.

\vspace{7pt}

\noindent\textbf{Remark}\quad In principle one could give a complete unification of Chapter~\ref{chap:basics} with the above modifications to it by phrasing all of these results in terms of a general (not necessarily full) subcategory $\bf{Cat}$ of the category $\G\hbox{-\textbf{Sys}}$ of all $\G$-systems, and adopting a flexible meaning for the term `relatively independent'.  In this work we have preferred to draw a more informal parallel between our two settings of interest, but it may be instructive to deduce from the proofs of Chapter 2 what basic properties we really need for the existence of sated extensions and the various lemmas that support it.  Although we leave the proof to the reader, it turns out that $\bf{Cat}$ must admit two basic constructions:
\begin{itemize}
\item it must have inverse limits;
\item it must have generated factors: that is, if
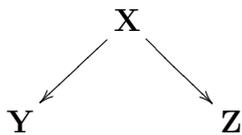
\begin{center}
$\phantom{i}$\xymatrix{& \bfX\ar[dl]\ar[dr]\\ \bfY & & \bfZ
}
\end{center}
is a diagram in $\bf{Cat}$, then there is an essentially unique minimal system $\bfW$ that may be inserted into this diagram as
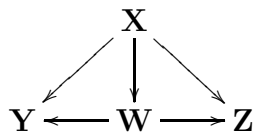
\begin{center}
$\phantom{i}$\xymatrix{& \bfX\ar[dl]\ar[d]\ar[dr]\\ \bfY & \bfW\ar[l]\ar[r] & \bfZ
}
\end{center}
\end{itemize}

Note, interestingly, that it does \emph{not} seem to be essential that any diagram such as
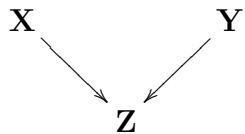
\begin{center}
$\phantom{i}$\xymatrix{\bfX\ar[dr] & & \bfY\ar[dl]\\ & \bfZ
}
\end{center}
have a common extension of $\bfX$ and $\bfY$ that can be inserted above it (of course, working in the whole of $\G\hbox{-\textbf{Sys}}$ when $\G$ is a group such a common extension is provided by the relatively independent product).

While these assumptions on $\bf{Cat}$ are relatively innocuous, more drastic steps are needed if we are to accommodate the instances of relative independence appearing in both Theorem~\ref{thm:sateds-exist} and Theorem~\ref{thm:Csateds-exist}.  The former of these asserts the relative independence of two whole factors of some extended system, whereas the latter concerns only the relative independent of functions of a single fixed coordinate within each of those factors (that is, relative independence under $\mu^{\rm{pt}}$ rather than $\mu$).  In order to treat these together, one could for example augment the category $\bf{Cat}$ by attaching to each system some distinguished subalgebra of bounded measurable functions (the whole of $L^\infty$ in the first case, and the subalgebra of functions of $x_w$ for some distinguished $w \in [k]^\ast$ in the second), and then re-defining conditional expectation as an operator acting only between these subalgebras for different systems and satisfying the usual conditions of idempotence and agreement of integrals against functions in the target subalgebra.

Altogether these very abstract considerations seem more demanding than worthwhile, and I know of few other situations in which a non-standard example of an abstract category of systems having these properties has been useful in ergodic theory. One related area which could fit into this mould is the study of partial exchangeability in probability theory, for which we refer the reader to Kallenberg's book~\cite{Kal02}, the survey papers~\cite{Aus--ERH,Ald10} and the references given there. \fin

\section{Another appeal to the infinitary hypergraph removal lemma}

The cw-idempotent classes for which we will apply Theorem~\ref{thm:Csateds-exist} are as follows.

\begin{dfn}[Partially insensitive processes]
Given a subset $e \subseteq [k]$, a process $(K^{[k]^\ast},\Psi^{\otimes [k]^\ast},\mu,T) \in \sfA$ is \textbf{$e$-insensitive} if its line marginal satisfies
\[x_i = x_j\quad\quad\hbox{for}\ \mu^{\rm{line}}\hbox{-a.e.}\ (x_1,x_2,\ldots,x_k)\in K^k\ \hbox{for all}\ i,j \in e.\]
We write $\sfA_e \subseteq \sfA$ for the subclass of all $e$-insensitive processes.
\end{dfn}

The persistence of $e$-insensitivity under inverse limits and joinings is immediate, and so we have:

\begin{lem}
The class $\sfA_e$ is cw-idempotent for each $e\subseteq [d]$. \qed
\end{lem}

In parallel with the developments of Section~\ref{sec:moreFberg}, given an arbitrary process
\[\bfX = (K^{[k]^\ast},\Psi^{\otimes [k]^\ast},\mu,T)\in \sfA,\]
for each $e \subseteq [d]$ we let $\Phi_e$ denote the \textbf{$e$-insensitive $\s$-subalgebra} of $\Psi$, consisting of those $A \in \Psi$ such that $\mu^{\rm{line}}(\pi_i^{-1}(A)\triangle \pi_j^{-1}(A)) = 0$ for all $i,j \in e$, where $\pi_i:K^k \to K$ is the coordinate projection.  Letting $\zeta_e:(K,\Psi) \to (K_e,\Psi_e)$ be some map of standard Borel spaces such that $\Phi_e$ agrees with $\{\zeta_e^{-1}(E):\ E \in \Psi_e\}$ modulo $\mu^{\rm{pt}}$-negligible sets, it follows that
\[\zeta^\ast_e:K^{[k]^\ast}\to K_e^{[k]^\ast}:(x_w)_w\mapsto (\zeta_e(x_w))_w\]
is a cw-factor map that coordinatizes $\bfX \to \sfA_e\bfX$.

Directly from the definition of $\Phi_e$ we observe that if $i,j \in e$ then $\pi_i^{-1}(\Phi_e)$ and $\pi_j^{-1}(\Phi_e)$ differ only by $\mu^{\rm{line}}$-negligible sets, and we denote their common $\mu^{\rm{line}}$-completion by $\Phi_e^\dag$.  If now $\I \subseteq \binom{[k]}{\geq 2}$ is an up-set, then similarly to the setup of Section~\ref{sec:moreFberg} we define $\Phi_\I := \bigvee_{e \in \I}\Phi_e$ and $\Phi^\dag_\I := \bigvee_{e \in \I}\Phi^\dag_e$.

In terms of these definitions, the consequences of cw-satedness that we need are now essentially parallel to Propositions~\ref{prop:Fberg1} and~\ref{prop:Fberg2}.

\begin{prop}\label{prop:line1}
For each $i\leq k$ let
\[\sfC_i := \bigvee_{j\leq k,\,j\neq i}\sfA_{\{i,j\}}.\]
If a system $\bfX$ with strongly stationary measure $\mu$ is cw-$\sfC_i$-sated for each $i$ then the $\s$-algebras $\pi_i^{-1}(\Psi) \leq \Psi^{\otimes k}$ are relatively independent under $\mu^{\rm{line}}$ over the further factors
\[\pi_i^{-1}\Big(\bigvee_{j\leq k,\,j\neq i}\Phi_{\{i,j\}}\Big).\]
\end{prop}

\noindent\textbf{Proof}\quad Clearly it will suffice to prove that $\pi_1^{-1}(\Psi)$ is relatively independent from $\pi_2^{-1}(\Psi)\vee\cdots\vee\pi_d^{-1}(\Psi)$ under $\mu^{\rm{line}}$ over
\[\Xi:= \bigvee_{j=2}^k\Phi_{\{1,j\}},\]
since the cases of the other coordinates under $\mu^{\rm{line}}$ then follow by symmetry.

We prove this by contradiction, so suppose that $f_1$, $f_2$, \ldots, $f_d \in L^\infty(\mu^{\rm{pt}})$ are such that
\[\int_{K^k}f_1\otimes f_2\otimes \cdots\otimes f_k\,\d\mu^{\rm{line}}\neq \int_{K^k}\sfE(f_1\,|\,\Xi)\otimes f_2\otimes \cdots\otimes f_k\,\d\mu^{\rm{line}}.\]
We will deduce from this a contradiction with the cw-satedness of $\mu$.  By replacing $f_1$ with $f_1 - \sfE(f_1\,|\,\Xi)$ it suffices to assume that $\sfE(f_1\,|\,\Xi) = 0$ but that the left-hand integral above does not vanish.

For each $j=2,3,\ldots,k$ recall the
letter-replacement map $r_{\{1,j\},j}:[k]^\ast\to[k]^\ast$ defined in Section~\ref{sec:startDHJ}. In view of the strong stationarity of $\mu$, we may transport the above non-vanishing integral to any combinatorial line in $[k]^\ast$: in particular, picking some $w \in [k]^\ast$ for which $w^{-1}\{j\} \neq \emptyset$ for every $j$, the points $\{w,r_{\{1,2\},2}(w),r_{\{1,3\},3}(w),\ldots,r_{\{1,k\},k}(w)\}$ form such a line, and so we have
\[\int_{X^{[k]^\ast}} f_1(x_w)\cdot f_2(x_{r_{\{1,2\},2}(w)})\cdot\cdots\cdot f_k(x_{r_{\{1,k\},k}(w)})\,\mu(\d\bf{x}) = \kappa\neq 0.\]

Now define the probability measure $\l$ on
$(K\times K^{\{2,3,\ldots,k\}})^{[k]^\ast}$ to be the joint law under
$\mu$ of
\[(x_w)_w \mapsto \big(x_w,x_{r_{\{1,2\},2}(w)},x_{r_{\{1,3\},3}(w)},\ldots,x_{r_{\{1,k\},k}(w)}\big)_w.\]
We see that all of its coordinate projections onto individual copies of $K$
are still just $\mu^{\rm{pt}}$, the cw-factor map
\[\phi^\ast_1:(x_w,y_{2,w},y_{3,w},\ldots,y_{k,w})_w\mapsto (x_w)_w\]
has $\phi^\ast_{1\#}\l = \mu$, and the cw-factor map
\[\phi^\ast_j:(x_w,y_{2,w},y_{3,w},\ldots,y_{k,w})_w\mapsto (y_{j,w})_w\]
for $j=2,3,\ldots,k$ is $\l$-almost surely
$\{1,j\}$-insensitive. Therefore through the cw-factor map $\phi^\ast_1$ the law $\l$ defines an extension of $\mu$ as a measure space.

This new measure $\l$ may not be strongly stationary, so may not define an extension of members of $\sfA$.  However, we can now repeat the trick of Lemma~\ref{lem:sssuffice}. By the Carlson-Simpson Theorem there are a subspace $\psi:[k]^\ast\into [k]^\ast$ and an infinite word $w
\in [k]^\bbN$ such that the pulled-back measures
\[T_{\psi(w_1w_2\cdots w_n\oplus\ \cdot\ )\#}\l\]
converge in the coupling topology on $(K\times K^{\{2,3,\ldots,k\}})^{[k]^\ast}$ (recall that for couplings of fixed marginal measures this is
compact; see Theorem 6.2 in~\cite{Gla03}) to a strongly stationary measure $\tilde{\mu}$. Since $\mu$ was already
strongly stationary, we must still have $\phi^\ast_{1\#}\tilde{\mu}= \mu$,
and by the definition of the coupling topology as the weakest for
which integration of fixed product functions is continuous it
follows that we must still have, firstly, that
\[\int_{(K\times K^{\{2,3,\ldots,k\}})^{[k]^\ast}}(f\circ\pi_u\circ\phi^\ast_1)\cdot\prod_{j \in [k]\setminus e}(h_j\circ\pi_u\circ\phi^\ast_j)\,\d\tilde{\mu}  =\kappa \neq 0\]
for each $u\in [k]^\ast$ (where now we may omit the assumption that
$u$ contains every letter at least once, by strong stationarity),
and secondly that the cw-factors generated by the maps $\phi^\ast_j$
are $\{1,j\}$-insensitive under $\tilde{\mu}$, since this is
equivalent to the assertion that for any $A \in \Psi$ and line
$\ell:[k]\into [k]^\ast$ we have
\[\int_{(K\times K^{\{2,3,\ldots,k\}})^{[k]^\ast}} 1_A(\phi_j(z_{\ell(1)}))\cdot 1_{K\setminus A}(\phi_j(z_{\ell(j)}))\,\tilde{\mu}(\d\bf{z}) = 0\]
and this is clearly a closed condition in the coupling topology.

It follows that this strongly stationary measure $\t{\mu}$ gives a genuine cw-extension $\phi^\ast_1:\t{\bfX}\to \bfX$ such that the lift of $f_1\circ\pi_1$ as a function of any one coordinate must have a nontrivial inner product with some pointwise product of $\{1,j\}$-insensitive
functions under $\t{\mu}$ over $j=2,3,\ldots,k$.  Hence this lift has nonzero conditional expectation onto a $\s$-subalgebra of $\Psi\otimes \Psi^{\otimes \{2,3,\ldots,k\}}$ coordinatizing a cw-factor in the class $\sfC_1$, but recalling our assumption that $\sfE(f_1\,|\,\Xi) = 0$, this provides the desired contradiction with cw-$\sfC_1$-satedness. \qed

\begin{prop}\label{prop:line2}
For each $e \subseteq [k]$ let
\[\sfC_e:= \bigvee_{j\in [d]\setminus e}\sfA_{e\cup \{j\}}.\]
If $\bfX$ is cw-$\sfC_e$-sated for every $e$ then for any up-sets
$\I,\I'\subseteq \binom{[k]}{\geq 2}$ the $\s$-subalgebras
$\Phi^\dag_\I$ and $\Phi^\dag_{\I'}$ are relatively independent
under $\mu^{\rm{line}}$ over $\Phi^\dag_{\I\cap \I'}$.
\end{prop}

\noindent\textbf{Proof}\quad As for Proposition~\ref{prop:Fberg2} we start with the case in which $\I' = \langle e\rangle$ for $e$ a member of $\binom{[d]}{\geq 2}\setminus \I$ of maximal size, and again just as for that proposition it suffices to show that $\Phi^\dag_e$ is relatively independent from $\bigvee_{j \in [k]\setminus e}\pi_j^{-1}(\Psi)$ over $\bigvee_{j \in [k]\setminus e}\Phi_{e\cup \{j\}}^\dag$ under $\mu^{\rm{line}}$.

Again this is best proved by deriving a contradiction with cw-satedness.  Pick some $i \in e$, so $\Phi^\dag_e$ agrees with $\pi_i^{-1}(\Phi_e)$ up to negligible sets, let
\[\Xi := \bigvee_{j \in [k]\setminus e}\Phi_{e\cup\{j\}},\]
and suppose we have some $f \in L^\infty(\mu^{\rm{pt}})$ that is $\Phi_e$-measurable and such that $\sfE(f\,|\,\Xi) = 0$, and also $h_j \in L^\infty(\mu^{\rm{pt}})$ for each $j\in [k]\setminus e$ such that
\[\int_{K^k} (f\circ\pi_i)\cdot \prod_{j\in [k]\setminus e}(h_j\circ\pi_j)\,\d\mu^{\rm{line}} = \kappa \neq 0.\]

Arguing as for the preceding proposition, this nonvanishing can be transported to any combinatorial line in $[k]^\ast$, including to a line such as $\{r_{e,1}(w)$, $r_{e,2}(w)$, $r_{e,3}(w)$, \ldots, $r_{e,k}(w)\}$ for any $w$ that contains every letter at least once.  This gives
\[\int_{K^{[k]^\ast}} f(x_{r_{e,i}(w)})\cdot \prod_{j\in [k]\setminus e}h_j(x_{r_{e,j}(w)})\,\mu(\d\bf{x}) = \kappa\]
for any such $w$, but since $f$ is $e$-insensitive we may replace the first factor in this integrand simply by $f(x_w)$.

It follows that if we define the probability measure $\l$ on
$(K\times K^{[k]\setminus e})^{[k]^\ast}$ to be the joint law under
$\mu$ of
\[(x_w)_w \mapsto \big(x_w,(x_{r_{e,j}(w)})_{j\in [k]\setminus e}\big)_w\]
then all of its coordinate projections onto individual copies of $K$
are still just $\mu^{\rm{pt}}$, the cw-factor map
\[\phi^\ast:\big(x_w,(y_{j,w})_{j\in [k]\setminus e}\big)_w\mapsto (x_w)_w\]
has $\phi^\ast_\#\l = \mu$ and the cw-factor maps
\[\phi^\ast_j:\big(x_w,(y_{j,w})_{j\in [k]\setminus e}\big)_w\mapsto (y_{j,w})_w\]
are $\l$-almost surely $(e\cup\{j\})$-insensitive.  Therefore
through $\phi^\ast$ the measure $\l$ is
an extension of the measure $\mu$, and the above inequality gives a non-zero
inner product for the lift of $f\circ\pi_u$ through $\phi^\ast$ with some product over $j\in
[k]\setminus e$ of
$(e\cup\{j\})$-insensitive functions under $\l$, which we can express as
\[\int_{K^{[k]^\ast}} (f\circ\pi_u\circ\phi^\ast)\cdot\prod_{j \in [k]\setminus e}(h_j\circ\pi_u\circ\phi^\ast_j)\,\d\l  = \kappa\]
for any $u\in [k]^\ast$ that contains each letter at least once.

To complete the proof, we may argue exactly as for
Proposition~\ref{prop:line1} that within the
not-necessarily-strongly-stationary law $\l$ we can find
infinite-dimensional subspaces $\psi$ for which the corresponding
image measures under $T_\psi$ converge in the coupling topology to a
strongly stationary extension $\t{\mu}$ of $\mu$, and such that this
extension preserves the feature that the lift of $f\circ\pi_u$ has a
nontrivial inner product with a pointwise product of
$(e\cup\{j\})$-insensitive functions under $\t{\mu}^{\rm{pt}}$ for
any word $u$.  By our assumption that $\sfE(f\,|\,\Xi) = 0$ this
gives a contradiction with cw-$\sfC_e$-satedness, as required.

The general case can now follows by induction on $\I'$ for each fixed $\I$ exactly as for Proposition~\ref{prop:Fberg2}. \qed

\noindent\textbf{Proof of Theorem~\ref{thm:inf-DHJ2}}\quad An initial application of Theorem~\ref{thm:Csateds-exist} allows us to assume that $\bfX$ is cw-sated for all the classes involved in Propositions~\ref{prop:line1} and~\ref{prop:line2}.

Next, exactly as for the proof of Theorem~\ref{thm:multirec2}, applying Proposition~\ref{prop:line1} shows that it suffices to prove Theorem~\ref{thm:inf-DHJ2} in case the sets $A_i$ lie in the $\s$-subalgebras $\Phi_{\langle i\rangle} = \bigvee_{j \in [d]\setminus \{i\}}\Phi_{\{i,j\}} \leq \Psi$.

Finally, it follows from the definitions and Proposition~\ref{prop:line2} that the probability space $(K,\Psi,\mu^{\rm{pt}})$, its self-coupling $\mu^{\rm{line}}$ and the $\s$-subalgebras $\Phi_e$ and their lifts $\Phi^\dag_e$ for $e \subseteq [d]$ satisfy all the conditions of the `infinitary removal result' Proposition~\ref{prop:infremoval}, so another appeal to that proposition completes the proof. \qed

\subsection*{Postscript to the above proof}

After the appearance of Furstenberg and Katznelson's original, technically rather demanding proof of Theorem B in~\cite{FurKat91}, considerable efforts were made to provide firstly a simpler proof, and more importantly one that could be made effective to deduce some quantitative bound on the necessary dependence of $N_0$ and $\delta$ and $k$.

Both of these goals were recently achieved by a large online collaboration, instigated by Tim Gowers and involving several other mathematicians, called Polymath1.  Importantly, their new proof does give a dependence of $N_0$ on $\delta$ and $k$ similar to the dependence obtained for the Multidimensional Szemer\'edi Theorem by using the hypergraph regularity and removal lemmas.  All these developments can be found online~(\cite{Gow(online)}) and in the preprint~\cite{DHJ09}.

Importantly, the infinitary proof of Theorem B that we have reported
above relies on an observation that was originally taken from their
work.  I will not attempt an exact translation here since the
lexicons of these two approaches are very different, but the outcome
for stochastic processes is essentially the observation that an
initially-given system $\bfX \in \sfA$ can be combined in a strongly
stationary joining with some $\{1,j\}$-insensitive systems as in our
proof of Proposition~\ref{prop:line1}, which then gives some
information on the structure of the original process $\bfX$ (in our
case by an appeal to cw-satedness).

\chapter{Coda: a general structural conjecture}\label{chap:spec}

It seems inadequate to finish this dissertation without discussing
at least some of the issues obviously left open by the preceding
chapters.  Perhaps most interesting for ergodic theory is the
meta-question introduced in Section~\ref{sec:meta}, and in this last
chapter I offer a few further speculations on what additional
answers to it we might hope for.

Our first clue in this direction is offered by the works~\cite{HosKra05} of Host and Kra and~\cite{Zie07} of Ziegler, establishing the special case of Theorem C corresponding to different powers of a single ergodic transformation: that is, the result that if $T:\bbZ\actson (X,\S,\mu)$ is ergodic and $f_1$, $f_2$, \ldots, $f_d \in L^\infty(\mu)$ then the averages
\[S_N(f_1,f_2,\ldots,f_d) := \frac{1}{N}\sum_{n=1}^N(f_1\circ T^n)\cdot (f_2\circ T^{2n})\cdot\cdots\cdot (f_d\circ T^{dn})\]
converge in $L^2(\mu)$ as $N\to\infty$.  Importantly, those two works both rest on a quite detailed result about `characteristic factors' for these averages:

\begin{thm}[Host-Kra Theorem]\label{thm:HK}
If $\bfX = (X,\S,\mu,T)$ is as above then there is a factor $\Phi \leq \S$ that is \textbf{characteristic} for the averages $S_N$ in the sense that
\[S_N(f_1,f_2,\ldots,f_d)\sim S_N(\sfE(f_1\,|\,\Phi),\sfE(f_2\,|\,\Phi),\ldots,\sfE(f_d\,|\,\Phi))\]
in $L^2(\mu)$ for any $f_1$, $f_2$, \ldots $f_d \in L^\infty(\mu)$ as $N\to\infty$, and which can be generated by a factor map to a $(d-1)$-step pro-nilsystem: that is, it can be generated by some increasing sequence of factor maps
\[\pi_n:(X,\S,\mu,T)\to (G_n/\G_n,\rm{Borel},m_{G_n/\G_n},R_{g_n})\]
to systems that are given by rotations by elements $g_n$ on compact $(d-1)$-step nilmanifolds $G_n/\G_n$.
\end{thm}

\noindent\textbf{Remark}\quad This notion of a characteristic factor is just a slight modification to that of a partially characteristic factor that we met in Proposition~\ref{prop:sated-implies-pleasant}. In fact, Ziegler proves in~\cite{Zie07} that there is a unique minimal factor with the properties given by the above theorem, and in Leibman's later treatment of these two proofs in ~\cite{Lei05(HKvsZ)} it is shown that the pro-nilsystem characteristic factors constructed by Host and Kra are precisely these minimal factors. \fin

This very surprising theorem asserts that for a completely arbitrary ergodic $\bbZ$-system $\bfX$, its nonconventional averages $S_N$ are entirely controlled by some highly-structured factor of $\bfX$, which can be expressed in terms of the very concrete data of rotations on nilmanifolds.  In this informal discussion we will assume familiarity with the definition and basic properties of such `nilsystems' here; they are treated thoroughly in~\cite{HosKra05} and~\cite{Zie07} and the references given there.

Host and Kra and Ziegler's proofs of the one-dimensional case of
Theorem C proceed via two different approaches to
Theorem~\ref{thm:HK}.  They are both rather longer than the proof in
our Chapter~\ref{chap:conv}, but using Theorem~\ref{thm:HK} they
give a much more precise picture of the limit.  On the other hand,
the strategy used in our Chapter~\ref{chap:conv} simply cannot be
specialized to the one-dimensional setting: it is essential for our
approach that the result be formulated for the linearly independent
directions $\bf{e}_1$, $\bf{e}_2$, \ldots, $\bf{e}_d \in \bbZ^d$.
This is because even if we are initially given a $\bbZ$-system
$(X,\S,\mu,T)$, we must re-interpret it as a $\bbZ^d$-system in
order to pass to an extension that is sated in the way required by
Proposition~\ref{prop:sated-implies-pleasant}.  To do this we define
a new $\bbZ^d$-action $T'$ on $X$ by $(T')^{\bf{e}_i} := T^i$, but
once we ascend to our sated extension this special structure of a
collection of powers of a single transformation will be lost, and so
we can no longer focus on the special, one-dimensional case of
convergence.  In a sense, this quiet assumption of linear
independence was a precursor to the discussion of
Section~\ref{sec:meta}: we need the linear independence of the
subgroups $\bbZ\bf{e}_i\leq \bbZ^d$ in order that a corresponding
notion of satedness has useful consequences.

However, these two very different approaches to different cases of
Theorem C do suggest a reconciliation of the issue raised at the end
of Section~\ref{sec:meta}: what becomes of our meta-question on the
possibly joinings of $\bbZ^D$-systems $\bfX_i \in \sfZ_0^{\G_i}$ if
the subgroups $\G_i$ are not linearly independent?  The centrepiece
of this final chapter is a conjectural answer to this question.  If
true, it would offer the first step in a complete `interpolation'
between the structural result~\ref{thm:HK} of Host and Kra and our
much softer result~\ref{thm:sateds-joint-dist}.

In order to formulate our conjecture, we first need some more notation.  The notion of an isometric extension of ergodic probability-preserving systems and the fact that any such can be coordinatized as a skew-product extension over the base system by some compact homogeneous space are very classical; see, for instance, Glasner's book~\cite{Gla03}.  Here we will also assume familiarity with a natural but less common generalization of this theory to the case in which the base system is not necessarily ergodic, in which the fibres of our extension must be allowed to vary in a suitable `measurable' way over the ergodic components of the base system.  This theory is set up generally in~\cite{Aus--ergdirint}, where the lengthy but routine work of re-establishing all the well-known theorems from the ergodic case is carried out in full, and we will also adopt the basic notations of that paper.

\begin{dfn}[Direct integral of pro-nilsystems]
If $\G$ is a discrete Abelian group then a $\G$-system $\bfX = (X,\S,\mu,T)$ is a \textbf{direct integral of $k$-step pro-nilsystems} if it admits a tower of factors
\[\bfX = \bfX_k \to \bfX_{k-1}\to \ldots\to \bfX_1 \to \bfX_0\]
in which the action of $\G$ on $\bfX_0$ is trivial, each extension $\bfX_i\to \bfX_{i-1}$ for $i\geq 1$ can be coordinatized as a relatively ergodic extension by measurably-varying compact metrizable Abelian group data
\begin{center}
$\phantom{i}$\xymatrix{\bfX_i\ar[dr]\ar@{<->}[rr]^-\cong & & \bfX_{i-1}\ltimes (A_{i,\bullet},m_{A_{i,\bullet}},\s_i)\ar[dl]^{\rm{canonical}}\\
& \bfX_{i-1},}
\end{center}
(so the measurable group data $A_{i,\bullet}$ really varies only over the base system $\bfX_0$) and for each ergodic component $\mu_s$ of $\mu$ the resulting $k$-step Abelian distal ergodic $\G$-system
\[(X,\S,\mu_s,T) \cong (A_{1,s}\times A_{2,s}\times\cdots\times A_{k,s},\rm{Borel},\rm{Haar},\s_1\ltimes \s_2\ltimes \cdots\ltimes \s_k)\]
is measure-theoretically isomorphic to an inverse limit of actions of $\G$ by commuting rotations on $k$-step nilmanifolds.
\end{dfn}

\noindent\textbf{Remark}\quad In fact it seems likely that the above class of systems can be set up in several different ways, which will presumably turn out to be equivalent.  I haven chosen the above definition here because I suspect it will ultimately prove relatively convenient for establishing the necessary properties of these systems, but an alternative has already appeared in the literature in the paper~\cite{ChuFraHos09} of Chu, Frantzikinakis and Host. \fin

The following lemma is now routine, given the ergodic case which is classical (it follows from the nilmanifold case of Ratner's Theorem: see, for instance,~\cite{Lei07,Lei10}).

\begin{dfn}
If $\L\leq \G$ is an inclusion of discrete Abelian groups, then the class $\sfZ_{\nil,k}^\L$ of those $\G$-systems whose $\L$-subactions are direct integral of $k$-step pro-nilsystems is an idempotent class of $\G$-systems. We refer to it as the class of \textbf{$\L$-partially $k$-step pro-nilsystems}. \qed
\end{dfn}

We are now ready to offer our conjectural strengthening of Theorem~\ref{thm:sateds-joint-dist} to the case of linearly dependent subgroups $\G_i$.

\begin{conj}[General Structural Conjecture]\label{conj}
Suppose that $\G_i \leq \bbZ^D$ for $i=1,2,\ldots,r$ are subgroups among which there are no pairwise inclusions and $n_1$, $n_2$, \ldots, $n_r \geq 0$ are integers. Then depending on these data there are finite families of pairs
\[(\L_{i,1},m_{i,1}),(\L_{i,2},m_{i,2}),\ldots,(\L_{i,k_i},m_{i,k_i})\quad\quad\hbox{for}\ i=1,2,\ldots,r\]
such that each $m_{i,j}\geq 0$ is an integer and $\L_{i,j}\leq \bbZ^d$ is a subgroup properly containing $\G_i$ for each $i,j$, and for which the following holds.

If $\bfX_i \in \sfZ_{\nil,n_i}^{\G_i}$ for each $i=1,2,\ldots,r$ and each $\bfX_i$ is sated with respect to all possible joins of classes of the form $\sfZ_{\nil,n}^\G$ for $\G \leq \bbZ^D$ and $n\geq 0$,
then for any joining $\pi_i:\bfY \to \bfX_i$, $i=1,2,\ldots,r$, the factors $\pi_i^{-1}(\S_i)$ are relatively independent over their further factors
\[\pi_i^{-1}\Big(\bigvee_{j=1}^{k_i}\Phi_{i,j}\Big)\]
where $\Phi_{i,j}$ is the factor of $\bfX_i$ generated by the factor map to $(\sfZ_0^{\G_i}\cap\sfZ_{\nil,m_{i,j}}^{\L_{i,j}})\bfX_i$.
\end{conj}

\noindent\textbf{Remark}\quad We invoke the `no-inclusions' condition on the subgroups $\G_i$ in order to avoid degenerate cases. Without it, we might for example be asking for the collection of all possible joinings between two systems $\bfX_i \in \sfZ_0^{\G_i}$ for $i=1,2$ with $\G_1 \geq \G_2$, and in this case Lemma~\ref{lem:two-fold-joinings} tells us something about the less constrained system $\bfX_2$, but on the side of the more constrained system $\bfX_1$ the joining may clearly be completely arbitrary. \fin

In particular, the case in which $\bfX_i$ has trivial
$\G_i$-subaction corresponds to $n_i = 0$, and in this case the
above conjecture asserts that given enough satedness, the factors
$\pi_i^{-1}(\S_i)$ of the joining system are relatively independent
over some further factors, each of which is assembled as a join of
systems from the classes
$\sfZ_0^{\G_i}\cap\sfZ_{\nil,m_{i,j}}^{\L_{i,j}}$. In particular,
while each of these ingredients may not be partially invariant under
any subgroup of $\bbZ^D$ strictly larger than $\G_i$, for each them
we do know \emph{something} quite concrete (in terms of
pro-nilsystems) about the subaction of some properly larger subgroup
$\L_{i,j}\gneqq\G_i$.

Of course, the above conjecture does not strictly cover
Theorem~\ref{thm:sateds-joint-dist}, since that gives much more
precise information on the pairs $(\L_{i,j},m_{i,j})$ in case the
$\G_i$ are linearly independent: to wit, the $\L_{i,j}$ are the sums
$\G_i + \G_\ell$ for $\ell\neq i$, and $m_{i,j} = 0$.  While a final
understanding of Conjecture~\ref{conj} would presumably also give a
recipe for producing these pairs in the general case (and so would
recover the exact details of our known special cases), the slightly
incomplete formulation of Conjecture~\ref{conj} seems ample for our
present discussion, and as I write this any sensible guess as to its
completion appears beyond reach.

Indeed, by itself Conjecture~\ref{conj} seems very optimistic, so it is worth mentioning some special cases of it beyond Theorem~\ref{thm:sateds-joint-dist} for which we have some supplementary evidence.

Firstly, if $D = 2$, each $n_i = 0$ and the $\G_i$ are pairwise
linearly-independent one-dimensional subgroups $\bbZ\bf{v}_i\leq
\bbZ^2$, then we can take a sensible guess at a more precise version
of the above conjecture: that any joining of systems $\bfX_i \in
\sfZ_0^{\G_i}$ should be relatively independent over the maximal
$(r-1)$-step pro-nilsystem factors $\bfX_i \to \sfZ_{\nil,r}\bfX_i$.
Indeed, this would simply correspond to the Host-Kra Theorem in the
case of the $\bbZ^2$-system
\[\vec{\bfX} := (X^k,\S^{\otimes k},\mu^\rm{F},\vec{T})\]
with $\vec{T}^{\bf{e}_1} := T\times T\times\cdots\times T$ and
$\vec{T}^{\bf{e}_2} := T\times T^2\times \cdots\times T^k$, where
now the subgroups are $\G_i = \bbZ(\bf{e}_2 - i\bf{e}_1)$ and the
coordinate projections $\pi_i:X^k\to X$ define factor maps to
suitable $\G_i$-partially-invariant $\bbZ^2$-systems $\bfX_i$,
constructed from $\bfX$ as in the proof of
Proposition~\ref{prop:Fberg1}.  In fact, I strongly suspect that the
methods of either~\cite{HosKra05} or~\cite{Zie07} could be adapted
directly to proving this more general result on the possible
joinings of such partially-invariant systems.  Other, similar
results on possible joinings of partially-invariant systems that do
not require any extensions but would correspond to further special
cases of Conjecture~\ref{conj} have appeared in Frantzikinakis and
Kra~\cite{FraKra05} (where nonconventional averages such as in our
Theorem C are studied, but subject to some additional hypotheses on
the individual ergodicity of several one-dimensional subactions), in
Chu~\cite{Chu09} and in Chu, Frantzikinakis and
Host~\cite{ChuFraHos09}.  In each of these cases, the joining in
question has been either the Furstenberg self-joining of some tuple
of commuting transformations, or the related Host-Kra self-joining
(originally defined in~\cite{HosKra05} for the case of powers of a
single transformation, and since adapted to the multi-directional
case in~\cite{Hos09,Chu09,ChuFraHos09}).  However, in each of these
cases it seems likely that the methods employed could be adapted to
proving a corresponding instance of Conjecture~\ref{conj}.

Another special case of Conjecture~\ref{conj}, the first beyond Theorem~\ref{thm:sateds-joint-dist} that \emph{does} require an ascent to sated extensions, appears in~\cite{Aus--lindeppleasant1,Aus--lindeppleasant2}. Indeed, the principal structural result of~\cite{Aus--lindeppleasant2} can be phrased as asserting that if $\bf{p}_1$, $\bf{p}_2$ and $\bf{p}_3 \in \bbZ^2$ are three directions which together with the origin $\bs{0}\in \bbZ^2$ lie in general position, then for a sufficiently sated system $\bfX$ the Furstenberg self-joining $\mu^\rm{F}$ of the quadruple of transformations $\id,T^{\bf{p}_1}$, $T^{\bf{p}_2}$, $T^{\bf{p}_3}$ is such that the coordinate projections $\pi_i:X^{\{0,1,2,3\}} \to X$ are relatively independent over their further factors
\[\pi_0^{-1}(\S^{T^{\bf{p}_1} = T^{\bf{p}_2}}\vee \S^{T^{\bf{p}_1} = T^{\bf{p}_3}}\vee \S^{T^{\bf{p}_2} = T^{\bf{p}_3}}\vee \S_{\nil,2}^T)\]
and
\[\pi_i^{-1}(\S^{T^{\bf{p}_i}}\vee \S^{T^{\bf{p}_i} = T^{\bf{p}_j}}\vee \S^{T^{\bf{p}_i} = T^{\bf{p}_k}}\vee \S_{\nil,2}^T)\quad\quad\hbox{for}\ \{i,j,k\} = \{1,2,3\}.\]
Arguing again as for Proposition~\ref{prop:Fberg1}, this would follow from a special case of Conjecture~\ref{conj} (again with some more precise information on the pairs $(\L_{i,j},m_{i,j})$) when $D=3$, $r=4$ and $\G_0,\G_1,\G_2,\G_3$ are four one-dimensional subgroups of $\bbZ^3$ any three of which are linearly independent.

At present no proof (or disproof) of Conjecture~\ref{conj} seems to be at hand.  Nevertheless, the various cases mentioned above do give me hope for it, and I strongly suspect that any result as powerful as this would constitute a major addition to our toolkit for approaching questions of multiple recurrence. For example, I would expect it to shed considerable new light on the Bergelson-Leibman conjecture on the convergence of `polynomial' nonconventional averages~\cite{BerLei02}.  For a recent discussion of these latter question see~\cite{Aus--lindeppleasant1,Aus--lindeppleasant2}, where the proof of an instance of this latter conjecture was the original motivation for the result on joint distributions mentioned above.

\bibliographystyle{uclathes}
\bibliography{bibfile}

\begin{thebibliography}{LRR03}

\bibitem[Ald]{Ald10}
David~J. Aldous.
\newblock ``More uses of exchangeability: representations of complex random
  structures.''
\newblock to appear in \emph{{P}robability and {M}athematical {G}enetics:
  {P}apers in {H}onour of {S}ir {J}ohn {K}ingman}.

\bibitem[Ausa]{Aus--DHJ}
Tim Austin.
\newblock ``Deducing the {D}ensity {H}ales-{J}ewett {T}heorem from an
  infinitary removal lemma.''
\newblock Preprint, available online at \verb|arXiv.org|: 0903.1633.

\bibitem[Ausb]{Aus--newmultiSzem}
Tim Austin.
\newblock ``Deducing the multidimensional {S}zemer\'edi {T}heorem from an
  infinitary removal lemma.''
\newblock To appear, \emph{{J}. d'{A}nalyse {M}ath.}

\bibitem[Ausc]{Aus--ergdirint}
Tim Austin.
\newblock ``Extensions of probability-preserving systems by measurably-varying
  homogeneous spaces and applications.''
\newblock Preprint, available online at \verb|arXiv.org|: 0905.0516.

\bibitem[Ausd]{Aus--lindeppleasant1}
Tim Austin.
\newblock ``Pleasant extensions retaining algebraic structure, {I}.''
\newblock Preprint, available online at \verb|arXiv.org|: 0905.0518.

\bibitem[Ause]{Aus--lindeppleasant2}
Tim Austin.
\newblock ``Pleasant extensions retaining algebraic structure, {II}.''
\newblock Preprint, available online at \verb|arXiv.org|: 0910.0907.

\bibitem[Aus08]{Aus--ERH}
Tim Austin.
\newblock ``On exchangeable random variables and the statistics of large graphs
  and hypergraphs.''
\newblock {\em Probability Surveys}, (5):80--145, 2008.

\bibitem[Aus09]{Aus--nonconv}
Tim Austin.
\newblock ``On the norm convergence of nonconventional ergodic averages.''
\newblock {\em Ergodic Theory Dynam. Systems}, {\bfseries 30}(2):321--338,
  2009.

\bibitem[Ber96]{Ber96}
Vitaly Bergelson.
\newblock ``Ergodic {R}amsey {T}heory -- an {U}pdate.''
\newblock In M.~Pollicott and K.~Schmidt, editors, {\em Ergodic Theory of
  $\mathbb{Z}^d$-actions: Proceedings of the Warwick Symposium 1993-4}, pp.
  1--61. Cambridge University Press, Cambridge, 1996.

\bibitem[BL02]{BerLei02}
V.~Bergelson and A.~Leibman.
\newblock ``A nilpotent {R}oth theorem.''
\newblock {\em Invent. Math.}, {\bfseries 147}(2):429--470, 2002.

\bibitem[Car88]{Car88}
Timothy~J. Carlson.
\newblock ``Some unifying principles in {R}amsey theory.''
\newblock {\em Discrete Math.}, {\bfseries 68}:117--169, 1988.

\bibitem[CFH]{ChuFraHos09}
Qing Chu, Nikos Frantzikinakis, and Bernard Host.
\newblock ``Commuting averages with polynomial iterates of distinct degrees.''
\newblock Preprint, available online at \verb|arXiv.org|: 0912.2641.

\bibitem[Chu09]{Chu09}
Qing Chu.
\newblock ``Convergence of weighted polynomial multiple ergodic averages.''
\newblock {\em Proc. Amer. Math. Soc.}, {\bfseries 137}:1363--1369, 2009.

\bibitem[CL84]{ConLes84}
Jean-Pierre Conze and Emmanuel Lesigne.
\newblock ``Th\'eor\`emes ergodiques pour des mesures diagonales.''
\newblock {\em Bull. Soc. Math. France}, {\bfseries 112}(2):143--175, 1984.

\bibitem[CL88a]{ConLes88.1}
Jean-Pierre Conze and Emmanuel Lesigne.
\newblock ``Sur un th\'eor\`eme ergodique pour des mesures diagonales.''
\newblock In {\em Probabilit\'es}, volume 1987 of {\em Publ. Inst. Rech. Math.
  Rennes}, pp. 1--31. Univ. Rennes I, Rennes, 1988.

\bibitem[CL88b]{ConLes88.2}
Jean-Pierre Conze and Emmanuel Lesigne.
\newblock ``Sur un th\'eor\`eme ergodique pour des mesures diagonales.''
\newblock {\em C. R. Acad. Sci. Paris S\'er. I Math.}, {\bfseries
  306}(12):491--493, 1988.

\bibitem[ET36]{ErdTur36}
P.~Erd\H{o}s and P.~Tur\'an.
\newblock ``On some sequences of integers.''
\newblock {\em J. London Math. Soc.}, {\bfseries 11}:261--264, 1936.

\bibitem[FK78]{FurKat78}
Hillel Furstenberg and Yitzhak Katznelson.
\newblock ``An ergodic {S}zemer\'edi {T}heorem for commuting transformations.''
\newblock {\em J. d'Analyse Math.}, {\bfseries 34}:275--291, 1978.

\bibitem[FK91]{FurKat91}
Hillel Furstenberg and Yitzhak Katznelson.
\newblock ``A {D}ensity {V}ersion of the {H}ales-{J}ewett {T}heorem.''
\newblock {\em J. d'Analyse Math.}, {\bfseries 57}:64--119, 1991.

\bibitem[FK05]{FraKra05}
Nikos Frantzikinakis and Bryna Kra.
\newblock ``Convergence of multiple ergodic averages for some commuting
  transformations.''
\newblock {\em Ergodic Theory Dynam. Systems}, {\bfseries 25}(3):799--809,
  2005.

\bibitem[Fur67]{Fur67}
Harry Furstenberg.
\newblock ``Disjointness in ergodic theory, minimal sets, and a problem in
  {D}iophantine approximation.''
\newblock {\em Math. Systems Theory}, {\bfseries 1}:1--49, 1967.

\bibitem[Fur77]{Fur77}
Hillel Furstenberg.
\newblock ``Ergodic behaviour of diagonal measures and a theorem of
  {S}zemer\'edi on arithmetic progressions.''
\newblock {\em J. d'Analyse Math.}, {\bfseries 31}:204--256, 1977.

\bibitem[Fur81]{Fur81}
Hillel Furstenberg.
\newblock {\em Recurrence in Ergodic Theory and Combinatorial Number Theory}.
\newblock Princeton University Press, Princeton, 1981.

\bibitem[FW96]{FurWei96}
Hillel Furstenberg and Benjamin Weiss.
\newblock ``A mean ergodic theorem for
  $\frac{1}{N}\sum_{n=1}^{N}f({T}^nx)g({T}^{n^2}x)$.''
\newblock In Vitaly Bergleson, ABC March, and Joseph Rosenblatt, editors, {\em
  Convergence in Ergodic Theory and Probability}, pp. 193--227. De Gruyter,
  Berlin, 1996.

\bibitem[Gla03]{Gla03}
Eli Glasner.
\newblock {\em Ergodic {T}heory via {J}oinings}.
\newblock American {M}athematical {S}ociety, {P}rovidence, 2003.

\bibitem[Gow98]{Gow98}
W.~T. Gowers.
\newblock ``A new proof of {S}zemer\'edi's theorem for arithmetic progressions
  of length four.''
\newblock {\em Geom. Funct. Anal.}, {\bfseries 8}(3):529--551, 1998.

\bibitem[Gow01]{Gow01}
W.~T. Gowers.
\newblock ``A new proof of {S}zemer\'edi's theorem.''
\newblock {\em Geom. Funct. Anal.}, {\bfseries 11}(3):465--588, 2001.

\bibitem[Gow07]{Gow07}
W.~T. Gowers.
\newblock ``Hypergraph regularity and the multidimensional {S}zemer\'edi
  theorem.''
\newblock {\em Ann. of Math. (2)}, {\bfseries 166}(3):897--946, 2007.

\bibitem[GRS90]{GraRotSpe90}
Robert~L. Graham, Bruce~L. Rothschild, and Joel~H. Spencer.
\newblock {\em Ramsey {T}heory}.
\newblock John Wiley \& Sons, New York, 1990.

\bibitem[HJ63]{HalJew63}
A.~W. Hales and R.~I. Jewett.
\newblock ``Regularity and positional games.''
\newblock {\em Trans. Amer. Math. Soc.}, {\bfseries 106}:222--229, 1963.

\bibitem[HK]{HosKra07}
Bernard Host and Bryna Kra.
\newblock ``Uniformity seminorms on $\ell^\infty$ and applications.''
\newblock Preprint, available online at \verb|arXiv.org|: 0711.3637.

\bibitem[HK01]{HosKra01}
Bernard Host and Bryna Kra.
\newblock ``Convergence of {C}onze-{L}esigne averages.''
\newblock {\em Ergodic Theory Dynam. Systems}, {\bfseries 21}(2):493--509,
  2001.

\bibitem[HK05]{HosKra05}
Bernard Host and Bryna Kra.
\newblock ``Nonconventional ergodic averages and nilmanifolds.''
\newblock {\em Ann. Math.}, {\bfseries 161}(1):397--488, 2005.

\bibitem[Hos09]{Hos09}
Bernard Host.
\newblock ``Ergodic seminorms for commuting transformations and applications.''
\newblock {\em Studia Math.}, {\bfseries 195}(1):31--49, 2009.

\bibitem[Kal02]{Kal02}
Olav Kallenberg.
\newblock {\em Foundations of modern probability}.
\newblock Probability and its Applications (New York). Springer-Verlag, New
  York, second edition, 2002.

\bibitem[Lei05]{Lei05(HKvsZ)}
A.~Leibman.
\newblock ``Host-Kra and Ziegler factors and convergence of multiple
  averages.''
\newblock In B.~Hasselblatt and A.~Katok, editors, {\em Handbook of Dynamical
  Systems}, volume~1B, pp. 841--853. Elsevier, 2005.

\bibitem[Lei07]{Lei07}
A.~Leibman.
\newblock ``Orbits on a nilmanifold under the action of a polynomial sequence
  of translations.''
\newblock {\em Ergodic Theory Dynam. Systems}, {\bfseries 27}(4):1239--1252,
  2007.

\bibitem[Lei10]{Lei10}
A.~Leibman.
\newblock ``Orbit of the diagonal in the power of a nilmanifold.''
\newblock {\em Trans. Amer. Math. Soc.}, {\bfseries 362}(3):1619--1658, 2010.

\bibitem[LPT00]{LemParTho00}
M.~Lema{\'n}czyk, F.~Parreau, and J.-P. Thouvenot.
\newblock ``Gaussian automorphisms whose ergodic self-joinings are
  {G}aussian.''
\newblock {\em Fund. Math.}, {\bfseries 164}(3):253--293, 2000.

\bibitem[LRR03]{LesRitdelaRue03}
E.~Lesigne, B.~Rittaud, and T.~de~la Rue.
\newblock ``Weak disjointness of measure-preserving dynamical systems.''
\newblock {\em Ergodic Theory Dynam. Systems}, {\bfseries 23}(4):1173--1198,
  2003.

\bibitem[NRS06]{NagRodSch06}
Brendan Nagle, Vojt{\v{e}}ch R{\"o}dl, and Mathias Schacht.
\newblock ``The counting lemma for regular {$k$}-uniform hypergraphs.''
\newblock {\em Random Structures Algorithms}, {\bfseries 28}(2):113--179, 2006.

\bibitem[Pola]{DHJ09}
D.~H.~J. Polymath.
\newblock ``A new proof of the density {H}ales-{J}ewett theorem.''
\newblock Preprint, available online at \verb|arXiv.org|: 0910.3926.

\bibitem[Polb]{Gow(online)}
Polymath1.
\newblock ``A combinatorial approach to {D}ensity {H}ales-{J}ewett.''
\newblock Online project, viewable at http://gowers.wordpress.com/.

\bibitem[Rot53]{Rot53}
K.~F. Roth.
\newblock ``On certain sets of integers.''
\newblock {\em J. London Math. Soc.}, {\bfseries 28}:104--109, 1953.

\bibitem[RS04]{RodSko04}
Vojt{\v{e}}ch R{\"o}dl and Jozef Skokan.
\newblock ``Regularity lemma for {$k$}-uniform hypergraphs.''
\newblock {\em Random Structures Algorithms}, {\bfseries 25}(1):1--42, 2004.

\bibitem[Rue]{delaRue09}
Thierry de~la Rue.
\newblock ``Notes on {A}ustin's multiple ergodic theorem.''
\newblock Unpublished, available online at \verb|arXiv.org|: 0907.0538.

\bibitem[Shk05]{Shk05}
Ilya~D. Shkredov.
\newblock ``On a problem of Gowers.''
\newblock {\em Dokl. Akad. Nauk}, {\bfseries 400}(2):169--172, 2005.
\newblock (Russian).

\bibitem[Sze75]{Sze75}
Endre Szemer\'edi.
\newblock ``On sets of integers containing no $k$ elements in arithmetic
  progression.''
\newblock {\em Acta Arith.}, {\bfseries 27}:199--245, 1975.

\bibitem[Tao07]{Tao07}
Terence Tao.
\newblock ``A correspondence principle between (hyper)graph theory and
  probability theory, and the (hyper)graph removal lemma.''
\newblock {\em J. d'Analyse Math.}, {\bfseries 103}:1--45, 2007.

\bibitem[Tao08]{Tao08(nonconv)}
Terence Tao.
\newblock ``Norm convergence of multiple ergodic averages for commuting
  transformations.''
\newblock {\em Ergodic Theory and Dynamical Systems}, {\bfseries 28}:657--688,
  2008.

\bibitem[Tow09]{Tow09}
Henry~P. Townser.
\newblock ``Convergence of {D}iagonal {E}rgodic {A}verages.''
\newblock {\em Ergodic Theory Dynam. Systems}, {\bfseries 29}:1309--1326, 2009.

\bibitem[TV06]{TaoVu06}
Terence Tao and Van Vu.
\newblock {\em Additive combinatorics}.
\newblock Cambridge University Press, Cambridge, 2006.

\bibitem[Wae27]{vdW27}
B.~L. van~der Waerden.
\newblock ``Beweis einer Baudetschen Vermutung.''
\newblock {\em Nieuw. Arch. Wisk.}, {\bfseries 15}:212--216, 1927.

\bibitem[Zha96]{Zha96}
Qing Zhang.
\newblock ``On convergence of the averages {$(1/N)\sum\sp N\sb {n=1}f\sb 1(R\sp
  nx)f\sb 2(S\sp nx)f\sb 3(T\sp nx)$}.''
\newblock {\em Monatsh. Math.}, {\bfseries 122}(3):275--300, 1996.

\bibitem[Zie07]{Zie07}
Tamar Ziegler.
\newblock ``Universal characteristic factors and {F}urstenberg averages.''
\newblock {\em J. Amer. Math. Soc.}, {\bfseries 20}(1):53--97 (electronic),
  2007.

\end{thebibliography}

\end{document}